\documentclass[12pt,a4paper]{article}

\usepackage{amsfonts,amsmath,amsthm,amssymb,amsopn}
\usepackage{graphicx}
\usepackage{epstopdf}
\usepackage{algorithmic}
\usepackage{enumitem}

\usepackage{bm,color}
\usepackage{stmaryrd}
\usepackage[T1]{fontenc}
\usepackage{float,subcaption}
\usepackage{booktabs}
\usepackage{multirow}
\usepackage{siunitx}
\usepackage{pifont}
\usepackage{hyperref}
\usepackage{cleveref}
\usepackage{lscape}

\usepackage{bbold}
\newcommand*{\I}{\mathbb{I}}

\newcommand\bu{\bm{u}}

\newcommand\bx{\bm{x}}

\newcommand\dd{\mathrm{d}}
\newcommand\pd[2]{\frac{\partial {#1}}{\partial {#2}}}

\newcommand\abs[1]{\lvert #1 \rvert}
\newcommand\norm[1]{\left\lVert #1 \right\rVert}

\newcommand\xr{i+\frac12}
\newcommand\xl{i-\frac12}
\newcommand\yr{j+\frac12}
\newcommand\yl{j-\frac12}
\newcommand\zr{k+\frac12}
\newcommand\zl{k-\frac12}

\theoremstyle{plain}
\newtheorem{definition}{Definition}[section]

\theoremstyle{definition}
\newtheorem{example}{Example}[section]
\newtheorem{proposition}{Proposition}[section]
\newtheorem{remark}{Remark}[section]
\newtheorem{lemma}{Lemma}[section]

\crefname{equation}{Equation}{Equations}
\crefname{figure}{Figure}{Figures}
\crefname{table}{Table}{Tables}
\crefname{example}{Example}{Examples}
\crefname{section}{Section}{Sections}

\renewcommand{\title}{A compact fourth-order doubly conservative active flux method with maximum-principle-preserving limiting for degenerate convection--diffusion equations}

\newcommand{\authorOne}{Junming Duan\footnote{School of Science and Engineering, The Chinese University of Hong Kong (Shenzhen), Longgang, 518172 Shenzhen, Guangdong, duanjunming@cuhk.edu.cn}}

\begin{document}

\begin{center} \Large
\title

\vspace{1cm}

\date{}
\normalsize

\authorOne
\end{center}

\begin{abstract}

The active flux (AF) method is a compact high-order finite volume method originally proposed for solving hyperbolic conservation laws, that evolves cell averages together with point values shared by neighboring cells.
This paper develops a compact fourth-order AF method for scalar degenerate convection--diffusion equations on Cartesian meshes.
Fourth-order accuracy is achieved without introducing additional unknowns relative to the standard third-order AF method.
For convection, a downwind point value is incorporated into the biased point-value stencil.
For diffusion, the proposed method directly discretizes the diffusion potential using compact fourth-order operators, thereby avoiding auxiliary gradient variables, hyperbolic reformulations, pseudo-time stepping, or inner iterations in the existing high-order AF methods.
To deal with degeneracy, the diffusion contribution to the cell-average evolution is written in a doubly conservative form, while the point values are evolved using compact finite-difference operators, thanks to the flexibility of point values.
To enforce the discrete maximum principle, a monolithic convex limiting is constructed, which simultaneously limits convection fluxes and diffusion potentials by blending them with their low-order maximum-principle-preserving counterparts, thus the conservation for convection and double conservation for diffusion are maintained.
The point values are projected onto the admissible interval, and an optional smoothness-indicator-based blending is employed to control oscillations.
For the one-dimensional linear problem, Fourier analysis establishes fourth-order spatial accuracy, while von Neumann analysis yields larger stable CFL limits than the fourth-order $P^3$ local discontinuous Galerkin method for the explicit Runge--Kutta methods considered.
Numerical experiments demonstrate the fourth-order convergence, preservation of prescribed bounds, and accurate resolution of shocks and sharp wave fronts in strongly degenerate problems.



\end{abstract}


\section{Introduction}\label{sec:introduction}
Convection--diffusion equations frequently arise in various physical and engineering contexts, such as heat conduction, diffusion processes, and fluid flow in porous media encountered in reservoir engineering.
This paper is concerned with solving the following convection--diffusion equation
\begin{equation}\label{eq:cde_nonconservative}
	u_t+\nabla\cdot\bm{f}(u)
	=\nabla\cdot\bigl(\kappa(u)\nabla u\bigr),
	\quad (\bx,t)\in\mathbb{R}^d\times\mathbb{R}^{+},
\end{equation}
where $u(\bx,t)\in\mathbb{R}$ is the scalar unknown, 
$\bm{f}=(f_1,\ldots,f_d)$ is the convection flux, and
$\kappa(u)\geqslant 0$ is the diffusion coefficient.
Introducing the diffusion potential
\begin{equation*}
	a(u)=\int^u\kappa(s)\,\dd s,
\end{equation*}
equation~\eqref{eq:cde_nonconservative} can equivalently be written as
\begin{equation}\label{eq:cde_conservative}
	u_t+\nabla\cdot\bm{f}(u)=\Delta a(u).
\end{equation}
When $\kappa(u)$ vanishes for some $u$, the equation becomes degenerate and exhibits both hyperbolic and parabolic behavior. 
Its solution may contain shock waves and sharp wave fronts with finite velocity.
These features make it important for a numerical method not only to attain high-order accuracy for smooth solutions, but also to approximate the relevant weak or ultraweak solution in degenerate regimes
\cite{Evje_2000_Monotone_SJNA,Jerez_2017_Entropy_SJoNA,Yu_2026_Lax_MC}.

For the periodic domain, the solution satisfies the maximum principle
\begin{equation*}
	u(\bx,t)\in[u_{\min},u_{\max}],
	\quad
	u_{\min}=\min_{\bx}u_0(\bx),
	\quad
	u_{\max}=\max_{\bx}u_0(\bx),
\end{equation*}
where $u_0$ denotes the initial data.
Preserving such bounds at the discrete level, referred to as the maximum-principle-preserving (MPP) property, is important for robustness and for suppressing nonphysical overshoots and undershoots near shocks and moving interfaces.
Low-order monotone schemes provide a robust foundation for degenerate convection--diffusion equations
\cite{Evje_2000_Monotone_SJNA}, while high-order WENO and discontinuous Galerkin discretizations have been developed for nonlinear degenerate parabolic problems
\cite{Zhang_2009_Numerical_JSC,Liu_2011_High_SJSC,Jiang_2021_High_JSC,Xu_2025_High_JCP}.
MPP high-order methods for convection--diffusion equations have also been studied in
\cite{Zhang_2012_Maximum_SJoSC,Yu_2025_High_JCAM}.
To deal with degeneracy, the Lax--Wendroff-type result in \cite{Yu_2026_Lax_MC} shows that, under the given convergence, consistency, and bounded-variation assumptions, a doubly conservative discretization permits the two discrete summations by parts needed to recover an ultraweak solution.
This provides an important criterion for the numerical methods for degenerate cases.

The active flux (AF) method is a compact finite volume method originally developed for hyperbolic conservation laws
\cite{Eymann_2011_Active_InCollection,Eymann_2011_Active_InProceedings,Eymann_2013_Multidimensional_InCollection,Roe_2017_Is_JSC}.
In addition to cell averages, AF methods evolve point values shared by neighboring cells.
The cell averages satisfy a conservative finite volume update, whereas the point values provide a continuous representation at cell interfaces and make the physical flux directly available there.
The AF methods evolve the point values using exact or approximate evolution operators \cite{Eymann_2013_Active,Fan_2015_Investigations_InCollection,Fan_2017_Acoustic,Barsukow_2021_active_JSC,Helzel_2019_new_JSC,Chudzik_2024_Active_JSC,Barsukow_2025_Active,Chudzik_2025_Fully}.

A method-of-lines formulation was subsequently introduced in
\cite{Abgrall_2023_combination_CAMC}, in which both cell averages and point values are discretized in space first and advanced by Runge--Kutta methods.
Upwind information in the point-value update can be incorporated through Jacobian or flux-vector splitting
\cite{Abgrall_2023_combination_CAMC,Abgrall_2023_Extensions_EMMNA,Duan_2025_Active_SJSC}.
This formulation has facilitated extensions to higher-order accuracy and more general meshes
\cite{Abgrall_2023_Active_BJPAM,Barsukow_2026_generalized_CF}.

The compactness and flexibility of the AF degrees of freedom (DoFs) make the method attractive for convection--diffusion problems, but comparatively few AF discretizations of diffusion are available.
One approach replaces the parabolic part by a first-order hyperbolic system with auxiliary gradient variables
\cite{Nishikawa_2016_Third_CF,He_2019_novel_InProceedings,Duan_2025_asymptotic}.
Depending on the formulation and the target problem, this approach may require pseudo-time integration, inner iterations, or the solution of an additional steady relaxation problem.
Another approach introduces a degenerate first-order system without time derivatives for the auxiliary variables \cite{He_2021_Towards}, which is related to local discontinuous Galerkin (LDG) methods
\cite{Cockburn_1998_local_SJNA}.
The diffusion discretization in \cite{He_2021_Towards}, however, is only second-order accurate.
These developments do not yet provide a compact high-order AF method that retains the appropriate doubly conservative formulation for strongly degenerate diffusion and enforces a discrete maximum principle.

In this work, we develop a compact doubly conservative AF method for scalar degenerate convection--diffusion equations on Cartesian meshes with suitable MPP limiting.
The main contributions and findings are summarized as follows.
\begin{enumerate}
	\item[(i)]
	We construct a compact fourth-order AF discretization using the same set of evolved unknowns as the standard third-order AF method.
	For convection, fourth-order accuracy is obtained by incorporating a downwind point value into the biased point-value stencil.
	For diffusion, the potential $a(u)$ is discretized directly using compact fourth-order operators.
	The resulting method introduces no auxiliary gradient variables and requires neither a hyperbolic reformulation nor pseudo-time stepping or inner iterations as needed in the existing AF methods.
	The construction is first presented in one dimension and then extended to Cartesian meshes in multiple dimensions.
	
	\item[(ii)]
	The diffusion contribution to the cell-average evolution is written in a doubly conservative form.
	Double conservation is imposed on the cell-average update, where it supplies the discrete form required by the Lax--Wendroff-type result in \cite{Yu_2026_Lax_MC}, while the flexibility of the point-value evolution is used to retain a compact fourth-order stencil.
	Numerical evidence from strongly degenerate examples shows that an otherwise consistent AF discretization without double conservation may fail to approximate the ultraweak solution, whereas the proposed doubly conservative method captures the correct wave fronts.
	
	\item[(iii)]
	We adapt the monolithic convex limiting of \cite{Kuzmin_2020_Monolithic_CMAME} to both convection fluxes and diffusion potentials.
	The high-order corrections are limited simultaneously using the same low-order MPP intermediate states.
	The resulting cell-average update preserves both local conservation of the convection part and double conservation of the diffusion part.
	Point values are projected onto the admissible interval, while an optional smoothness-indicator-based blending is employed separately to control oscillations near discontinuities.
	
	\item[(iv)]
	For the 1D linear problem, discrete Fourier analysis establishes fourth-order spatial accuracy.
	Von Neumann analysis is used to determine the CFL condition of the fully discrete method coupled with the Runge--Kutta methods \cite{Gottlieb_2001_Strong_SR,Spiteri_2002_New_SJNA}.
	The maximum stable CFL numbers are $1.049$ for pure convection and $0.33$ for pure diffusion.
	In the corresponding 1D comparison, these values are larger than those of the fourth-order $P^3$ LDG method.
	
	\item[(v)]
	Numerical experiments in one, two, and three dimensions demonstrate the observed fourth-order convergence for smooth solutions, preservation of prescribed bounds, and accurate resolution of shocks and sharp fronts in degenerate problems.
\end{enumerate}

The remainder of this paper is organized as follows.
Section~\ref{sec:1d_formulation} presents the 1D fourth-order doubly conservative AF discretization.
Section~\ref{sec:1d_fourier} analyzes its spatial accuracy and fully discrete linear stability.
The extension to Cartesian meshes in multiple dimensions is given in Section~\ref{sec:multid_formulation}.
Section~\ref{sec:limiting} introduces the maximum-principle-preserving limiting and the optional smoothness-indicator-based blending.
Numerical experiments are presented in Section~\ref{sec:results}, and concluding remarks are given in Section~\ref{sec:conclusion}.

\section{1D formulation}\label{sec:1d_formulation}
This section presents an AF method for solving
\begin{equation}\label{eq:1d_cde_dc}
	u_t + f(u)_x = a(u)_{xx}.
\end{equation}
Numerical evidence shows that, for strongly degenerate cases, an AF method without double conservation, presented in \ref{app:nodc} and based on \eqref{eq:1d_cde_nodc}, may fail to approximate the ultraweak solution, see also \cite{Yu_2026_Lax_MC}.
This motivates the doubly conservative AF discretization of \eqref{eq:1d_cde_dc} developed below.
This work focuses on the spatial discretization in the AF method.
The fully-discrete method is obtained using the explicit SSP-RK54 \cite{Spiteri_2002_New_SJNA}.

Assume that a 1D computational domain is divided into $N$ uniform cells
$I_i = [x_{\xl}, x_{\xr}]$ with the cell size $\Delta x = x_{\xr}-x_{\xl}$, $i=1,\cdots,N$.
The cell centers are denoted by $x_i = (x_{\xl} + x_{\xr})/2$.
The DoFs of the AF method consist of approximations to cell averages of the unknowns as well as point values at the cell interfaces,
\begin{equation*}
	\bar{u}_i(t) = \frac{1}{\Delta x} \int_{I_i} u_h(x, t) \dd x,\quad
	u_{\xr}(t) = u_h(x_{\xr}, t),
\end{equation*}
where $u_h(x,t)$ is the numerical solution.

\subsection{A fourth-order discretization for the convection part}
Following the finite volume method, one has
\begin{subequations}\label{eq:1d_conv}
	\begin{align}
		\frac{\dd \bar{u}_i }{\dd t} &+ \frac{1}{\Delta x}\left( f_{\xr} - f_{\xl} \right) = 0, \label{eq:1d_conv_av} \\
		\frac{\dd u_{\xr} }{\dd t} &+ D^{\texttt{4th}}_{+} (f^{+})_{\xr} + D^{\texttt{4th}}_{-}(f^{-})_{\xr} = 0, \label{eq:1d_conv_pnt}
	\end{align}
\end{subequations}
with the flux $f_{\xr} = f(u_{\xr})$.
Thanks to the point value DoFs in the AF method, the flux is available as soon as the point values $u_{\xr}$ are obtained.

Three approaches for constructing higher-order AF methods have been discussed for the 1D case \cite{Abgrall_2023_Extensions_EMMNA}.
To obtain a $4$th-order approximation to the derivative of some variable $z$ at $x_{\xr}$, we use the four-point biased stencil $z_{\xl}$, $z_{i}$, $z_{\xr}$, $z_{i+1}$, which is mentioned in \cite{Abgrall_2023_combination_CAMC}, and a similar idea is used in \cite{Zeng_2024_explicit_NMPDE}.
A cubic polynomial is reconstructed on this stencil, and its derivative is evaluated at $x_{\xr}$,
\begin{equation}\label{eq:1d_upwind_ops_4p}
	D^{\texttt{4th}}_{+}(z)_{\xr} = \frac{1}{3\Delta x}\left( z_{\xl} - 6z_{i} + 3z_{\xr} + 2z_{i+1} \right).
\end{equation}
With a mirror-symmetric stencil, one also has
\begin{equation}\label{eq:1d_upwind_ops_4m}
	D^{\texttt{4th}}_{-}(z)_{\xr} = \frac{1}{3\Delta x}\left( -2z_{i} - 3z_{\xr} + 6z_{i+1} - z_{i+\frac32} \right).
\end{equation}
Note that the cell-centered value is obtained through Simpson's rule
\begin{equation*}
	z_i = \tfrac{1}{4}( 6\bar{z}_i - (z_{\xl} + z_{\xr}) ),
\end{equation*}
which is a $4$th-order approximation.
The numerical flux is obtained by the local Lax--Friedrichs flux vector splitting (FVS) introduced in \cite{Duan_2025_Active_SJSC}, where the FVS gives
\begin{equation*}
	f^{\pm} = \frac{1}{2}\left(f(u) \pm \alpha_{\xr} u\right),
	\quad
	(f^{+}(u))' \geqslant 0,
	\quad
	(f^{-}(u))' \leqslant 0,
\end{equation*}
with the viscosity coefficient $\alpha_{\xr}$ determined by the maximum absolute value of the wave speed across the spatial stencil,
\begin{equation*}
	\alpha_{\xr} = \max\{ \abs{f'(u_{\xl})}, ~\abs{f'(u_{i})}, ~\abs{f'(u_{\xr})}, ~\abs{f'(u_{i+1})}, ~\abs{f'(u_{i+\frac32})} \}.
\end{equation*}

\subsection{A compact doubly conservative discretization for the diffusion part}\label{subsec:dc}
Consider Evje and Karlsen's example \cite{Evje_2000_Monotone_SJNA}, presented in Example~\ref{ex:1d_degenerate_Evje_Karlsen}.
The diffusion coefficient is zero on a set of nonzero measure, which is a strongly degenerate case.
Figure~\ref{fig:1d_degenerate_Evje_Karlsen} shows the solutions of a doubly conservative AF method to be constructed, and also of an AF method \eqref{eq:1d_diff_nodc} without double conservation in \ref{app:nodc}.
One observes that the double conservation is important for capturing the correct wave front in this strongly degenerate case.
The difficulty due to degeneracy will be overcome by considering the double conservation proposed in \cite{Yu_2026_Lax_MC}.
Thanks to the flexibility of the point-value evolution, the double conservation is only enforced where it is essential, i.e., the cell-average evolution, whereas a compact $4$th-order discretization is exploited for the point value.

Following \cite{Yu_2026_Lax_MC}, the relevant doubly conservative property can be given as follows for the present uniform mesh.
For a forward Euler step, a scheme for the cell average is doubly conservative if it can be written in the form
\begin{equation}\label{eq:1d_doc_definition}
	\frac{\widetilde{u}_i^{n+1}-\widetilde{u}_i^n}{\Delta t^n}
	+ \frac{\widetilde{f}_{i+\frac12}^n-\widetilde{f}_{i-\frac12}^n}{\Delta x}
	= \frac{\widetilde{a}_{i+1}^n-2\widetilde{a}_i^n+\widetilde{a}_{i-1}^n}{\Delta x^2},
\end{equation}
where $\widetilde{u}_i$, $\widetilde{f}_{i+\frac12}$, and $\widetilde{a}_i$ are approximations to the conservative variable, the convection flux $f(u)$, and the diffusion potential $a(u)$, respectively.
For a constant state $u$, consistency requires $\widetilde{u}_i=u$, $\widetilde{f}_{i+\frac12}=f(u)$, and $\widetilde{a}_i=a(u)$.
The general definition in \cite{Yu_2026_Lax_MC} additionally requires local Lipschitz continuity of these numerical approximations and recovery of global conservation from the locally conserved variables.
The term ``doubly conservative'' refers to the fact that the diffusion contribution is a difference of diffusion fluxes, with each diffusion flux itself represented as a difference of numerical diffusion potentials.
Under the additional convergence, consistency, and bounded-variation assumptions given in \cite{Yu_2026_Lax_MC}, this property permits two discrete summations by parts and yields an ultraweak solution in the limit.

The AF method can be built on \eqref{eq:1d_cde_dc} as follows,
\begin{subequations}\label{eq:1d_af_dc}
	\begin{align}
		&\frac{\dd \bar{u}_i }{\dd t} + \frac{1}{\Delta x}\left( f_{\xr} - f_{\xl} \right) = \frac{1}{\Delta x}\left( (a_x)_{\xr} - (a_x)_{\xl} \right), \label{eq:1d_af_av_dc} \\
		&(a_x)_{\xr} = D^{\texttt{4th}}_{c} a(u)_{\xr}, \nonumber \\
		&\frac{\dd u_{\xr} }{\dd t} + D^{\texttt{4th}}_{+} (f^{+})_{\xr} + D^{\texttt{4th}}_{-}(f^{-})_{\xr} = D^{2,\texttt{4th}} a(u)_{\xr}, \nonumber 
	\end{align}
\end{subequations}
where the central finite difference is obtained by averaging the biased operators \eqref{eq:1d_upwind_ops_4p} and \eqref{eq:1d_upwind_ops_4m} as
\begin{equation}\label{eq:1d_diff_central_4th}
	D^{\texttt{4th}}_{c} (z)_{\xr} = \frac12\left( D^{\texttt{4th}}_{+}(z) + D^{\texttt{4th}}_{-}(z) \right)_{\xr} =  \frac{1}{6\Delta x}\left( z_{\xl} - 8z_i + 8z_{i+1} - z_{i+\frac32} \right),
\end{equation}
and
\begin{equation}\label{eq:1d_diff_compact_a}
	D^{2,\texttt{4th}} a(u)_{\xr} = \frac{1}{\Delta x^2}\left(- \frac{1}{3}a_{i-\frac12} + \frac{16}{3}a_{i} - 10 a_{i+\frac12} + \frac{16}{3}a_{i+1} - \frac{1}{3}a_{i+\frac32}\right)
\end{equation}
is a $4$th-order approximation to $a(u)_{xx}$ at $x_{\xr}$.

\begin{proposition}
	With a forward Euler step, the cell-average update in \eqref{eq:1d_af_av_dc} has the doubly conservative form \eqref{eq:1d_doc_definition}, with
	$\widetilde{u}_i=\bar{u}_i$, $\widetilde{f}_{i+\frac12}=f(u_{i+\frac12})$, and $\widetilde{a}_i=\hat{a}_i$, where
	$\hat{a}_i = \dfrac{1}{6} \left( - a_{\xl} + 8a_{i} - a_{i+\frac12} \right)$.
\end{proposition}

\begin{proof}
	The central finite difference \eqref{eq:1d_diff_central_4th} can be rewritten as
	\begin{align*}
		D^{\texttt{4th}}_{c} (z)_{\xr} &= \frac{1}{6\Delta x}\left( z_{\xl} - 8z_i + 8z_{i+1} - z_{i+\frac32} \right) =: \dfrac{1}{\Delta x}\left( \hat{z}_{i+1} - \hat{z}_{i} \right),
	\end{align*}
	with $\hat{z}_i = \frac{1}{6} \left( - z_{\xl} + 8z_{i} - z_{i+\frac12} \right)$.
	Then \eqref{eq:1d_af_av_dc} can be simplified as 
	\begin{equation*}
		\frac{\dd \bar{u}_i }{\dd t}
		+ \frac{1}{\Delta x}\left( f(u_{\xr}) - f(u_{\xl}) \right)
		= \frac{1}{\Delta x^2}\left( \hat{a}_{i+1} - 2\hat{a}_{i} + \hat{a}_{i-1}\right)
	\end{equation*}
	as required.
\end{proof}

For completeness and clarity, a compact discretization \eqref{eq:1d_diff_nodc} for the diffusion part without double conservation is also presented in \ref{app:nodc}, which is used to generate the right panel in Figure~\ref{fig:1d_degenerate_Evje_Karlsen}.

\section{Fourier analysis}\label{sec:1d_fourier}
In this section, the accuracy and stability of the proposed AF method are studied by considering the following linear problem,
\begin{subequations}\label{eq:1d_linear_cde}
	\begin{align}
		&u_t + c u_x = \kappa_0 u_{xx}, \quad x\in[0,2\pi], ~t>0, \\
		&u(x,0) = \hat{u}_0\exp(\I\omega x), \label{eq:1d_heat_init}
	\end{align}
\end{subequations}
with constant $c, \kappa_0 > 0$, $\I=\sqrt{-1}$.
Its exact solution is $u(x,t) = \hat{u}_0\exp(\lambda t+\I\omega x)$ with $\lambda=-\I c\omega-\kappa_0\omega^2$.

For the linear case, the AF method reduces to the following form,
\begin{subequations}\label{eq:1d_af_constant}
	\begin{align}
		\frac{\dd \bar{u}_i }{\dd t}
		=&- \frac{c}{\Delta x} (u_{\xr} - u_{\xl}) \nonumber\\
		&+ \frac{\kappa_0}{2\Delta x^2}\left( -u_{i-\frac32} + 4\bar{u}_{i-1} + u_{\xl} - 8\bar{u}_{i} + u_{\xr} + 4\bar{u}_{i+1} - u_{i+\frac32} \right), \\
		\frac{\dd u_{\xr} }{\dd t}
		=&- \frac{c}{3\Delta x}\left( u_{\xl} - 6u_{i} + 3u_{\xr} + 2u_{i+1} \right) \nonumber\\
		&+ \frac{\kappa_0}{3\Delta x^2}\left( -5u_{i-\frac12} + 24\bar{u}_{i} - 38u_{\xr} + 24\bar{u}_{i+1} - 5u_{i+\frac32}\right).
	\end{align}
\end{subequations}

\subsection{Semi-discrete AF method}
Make the Fourier ansatz for the DoFs as
\begin{equation*}
	[\bar{u}_i(t), u_{\xr}(t)]^\top = \hat{\bm{u}}(t)\exp(\I\omega x_i) = [\hat{u}_1(t), \hat{u}_2(t)]^\top\exp(\I\omega x_i).
\end{equation*}
Let $t_x = \exp(\I\xi)$ with $\xi=\omega \Delta x\in[-\pi,\pi]$.
Substituting this ansatz into \eqref{eq:1d_af_constant} yields the following evolution equation for the Fourier coefficients,
\begin{align}\label{eq:1d_fourier_ode}
	&\frac{\dd}{\dd t} \hat{\bm{u}}
	= \left(\frac{c}{\Delta x}\bm{G}_c + \frac{\kappa_0}{\Delta x^2}\bm{G}_d\right) \hat{\bm{u}}
	=: \bm{G} \hat{\bm{u}}, \\
	&\bm{G}_c = 
	\begin{bmatrix}
		0 & t_x^{-1} -1 \\
		3 - t_x & \frac{1}{6}\left( -5t_x^{-1} - 8 + t_x \right)
		\\
	\end{bmatrix}, \nonumber\\
	&\bm{G}_d = 
	\begin{bmatrix}
		2(t_x^{-1} - 2 + t_x) & \frac12\left(-t_x^{-2} + t_x^{-1} + 1 - t_x\right) \\
		8(1 + t_x) & -\frac{1}{3}\left( 5t_x^{-1} + 38 + 5t_x \right) \\
	\end{bmatrix}.\nonumber
\end{align}

\begin{proposition}\label{prop:1d_eigval}
	Consider the AF method~\eqref{eq:1d_af_constant} for the linear convection--diffusion equation~\eqref{eq:1d_linear_cde} with a uniform mesh and periodic boundary conditions.
	Using the Fourier analysis, the semi-discrete Fourier matrix $\bm{G}$ has two eigenvalues.
	The physically relevant one, denoted by $\lambda_1$, approximates the analytical value $\lambda=-\I c\omega-\kappa_0\omega^2$ with $4$th-order accuracy.
	The other spurious one is of $\mathcal{O}(c/\Delta x + \kappa_0/\Delta x^2)$ with a negative real part.
	Thus the semi-discrete evolution equation~\eqref{eq:1d_fourier_ode} is stable.
\end{proposition}

\begin{proof}
	The calculations in this section can be obtained with the help of Mathematica.
	The characteristic polynomial of $\bm{G}$ is
	\begin{equation*}
		\det(\lambda \bm{I} - \bm{G})
		= \lambda^2 + c_1\lambda + c_2, 
	\end{equation*}
	with
	\begin{align*}
		c_1 = & - \frac{c (-5 + (-8 + t_x) t_x)}{6 \Delta x t_x} 
		- \frac{2\kappa_0 (1 + (-50+t_x)t_x)}{6 \Delta x^2 t_x}, \\
		c_2 = &\ \frac{(-1+t_x)}{6 \Delta x^{4} t_x^{2}}(
		6 c^{2} \Delta x^{2}(3-t_x)t_x
		-4 \kappa_0^{2}(1-t_x)(1+(-26+t_x)t_x) \\
		&+\kappa_0 c \Delta x(1+t_x(57-(-39+t_x)t_x))
		).
	\end{align*}
	The two eigenvalues are $\frac12\left(-c_1 \pm \sqrt{c_1^2 - 4c_2}\right)$.
	By Taylor expansion, one obtains
	\begin{align*}
		\lambda_1 &= - \left( \I c\omega + \kappa_0\omega^2\right) 
		+ \frac{1}{2880} (\I c + 7\kappa_0 \omega) \omega^5\Delta x^4 + \mathcal{O}(\Delta x^5) = \lambda + \mathcal{O}(\Delta x^4), \\
		\lambda_2 &= -\frac{2c}{\Delta x} - \frac{16\kappa_0}{\Delta x^2}
		+ \frac23\omega(3\I c + \kappa_0\omega) + \mathcal{O}(\Delta x).
	\end{align*}
	Therefore, $\lambda_1$ is physical with the $4$th-order accuracy, while $\lambda_2$ is spurious.
\end{proof}

To study the accuracy of the proposed method, consider the initial data $\hat{\bm{u}}_0\exp(\I\omega x_i)$ for the numerical solution in the cell $I_i=[x_{i-\frac12}, x_{i+\frac12}]$.
The coefficient $\hat{\bm{u}}_0 = [\hat u_1(0), \hat u_2(0)]^\top$ is obtained based on \eqref{eq:1d_heat_init},
\begin{equation*}
	\hat{\bm{u}}_0 = \exp(-\I \omega x_i) \left(\frac{1}{\Delta x}\int_{I_i} u(x,0) \dd x, ~ u(x_{\xr},0)\right)^\top= \hat u_0 \left(\tfrac{2}{\xi}\sin\left(\tfrac{\xi}{2}\right) , ~e^{\frac{\I\xi}{2}} \right)^\top.
\end{equation*}

\begin{proposition}\label{prop:1d_eigvec}
	With the same assumption as Proposition~\ref{prop:1d_eigval}, the eigenmode $\bm{V}_1$ corresponding to the physically relevant eigenvalue $\lambda_1$ approximates $\hat{\bm{u}}_0$ with the $4$th-order accuracy, while the contribution of the spurious mode to the projected initial data is of order $\mathcal{O}(\Delta x^4)$. 
\end{proposition}

\begin{proof}
	The corresponding eigenvectors of $\bm{G}$ are
	\begin{equation*}
		\bm{v}_1 = (\tilde{b}_1, 1)^\top,~
		\bm{v}_2 = (\tilde{b}_2, 1)^\top,
	\end{equation*}
	where
	\begin{align*}
		\tilde{b}_1 + \tilde{b}_2 &= \frac{2 (11 \kappa_0+c \Delta x) \cos \xi+26 \kappa_0-3 \I c \Delta x \sin
			\xi+4 c \Delta x}{3 e^{\I\xi} (8 \kappa_0-c \Delta x)+24 \kappa_0+9
			c \Delta x}, \\
		\tilde{b}_1\tilde{b}_2 &= \frac{\I e^{-\I \xi} (c
			\Delta x+ \I \kappa_0 \sin \xi)}{(8
			\kappa_0+c \Delta x) \cot \left(\frac{\xi }{2}\right)-2 \I c \Delta x}.
	\end{align*}
	The discrete initial data $\hat{\bm{u}}_0$ can be decomposed as $\hat{\bm{u}}_0 = \hat{u}_0^{(1)}\bm{v}_1 + \hat{u}_0^{(2)}\bm{v}_2 =: \bm{V}_1 + \bm{V}_2$, where
	\begin{equation*}
		\left(\hat{u}_0^{(1)}, \hat{u}_0^{(2)}\right) = \frac{\hat u_0}{\tilde{b}_1 - \tilde{b}_2}\left( \tfrac{2}{\xi}\sin\left(\tfrac{\xi}{2}\right) - e^{\frac{\I\xi}{2}}\tilde{b}_2, ~e^{\frac{\I\xi}{2}}\tilde{b}_1 - \tfrac{2}{\xi}\sin\left(\tfrac{\xi}{2}\right) \right).
	\end{equation*}
	By Taylor expansion, one has
	\begin{align*}
		&\bm{V}_1
		= \hat{\bm{u}}_0 + \left(0, -\dfrac{1}{2880}\hat{u}_0\omega^4\Delta x^4\right)^\top + \mathcal{O}(\Delta x^5),\\
		&\bm{V}_2 = \left(0, \dfrac{1}{2880}\hat{u}_0\omega^4\Delta x^4\right)^\top + \mathcal{O}(\Delta x^5).
	\end{align*}
\end{proof}

\begin{proposition}\label{prop:1d_accuracy}
	With the same assumption as Proposition~\ref{prop:1d_eigval}, let $\bm{u}(t) = \hat{\bm{u}}_0 \exp(\lambda t + \I\omega x_i)$ and $\bm{u}_h(t) = \hat{\bm{u}} \exp(\I\omega x_i)$ be the exact solution and numerical solution in the cell $I_i$.
	Then the error estimate for $t>0$ is
	\begin{equation*}
		\left\| \bm{e}(t) \right\|_{L^2} = \left\| \bm{u}(t) - \bm{u}_h(t) \right\|_{L^2}
		\leqslant C_1 t\Delta x^4 + C_2\Delta x^4,
	\end{equation*}
	where $C_1, C_2$ are positive constants independent of $\Delta x$.
	Thus the AF method is $4$th-order accurate in space.
\end{proposition}

\begin{proof}
	Following \cite{Guo_2013_Superconvergence_JCP}, the error can be decomposed as follows,
	\begin{align*}
		\left\| \bm{e}(t) \right\|_{L^2}
		&= \left \| \exp(\lambda t)\hat{\bm{u}}_0 - (\exp(\lambda_1t)\bm{V}_1 + \exp(\lambda_2t)\bm{V}_2) \right \|_{L^2} \\
		&\leqslant  | \exp(\lambda t) - \exp(\lambda_1 t)| \| \bm{V}_1 \|_{L^2}
		+ \| \hat{\bm{u}}_0 - \bm{V}_1 \|_{L^2}  \abs{\exp(\lambda t)}
		+ \abs{\exp(\lambda_2 t)} \| \bm{V}_2 \|_{L^2}.
	\end{align*}
	The second and third terms are $4$th-order accurate, according to Proposition~\ref{prop:1d_eigvec}, and $\operatorname{Re}(\lambda), \operatorname{Re}(\lambda_2) \leqslant 0$.
	In the first term, $\| \bm{V}_1 \|_{L^2}$ is of $\mathcal{O}(1)$, so that
	\begin{equation*}
		| \exp(\lambda t) - \exp(\lambda_1 t)| \| \bm{V}_1 \|_{L^2}
		\leqslant \widetilde{C}_1 t |\lambda - \lambda_1|
		\leqslant C_1 t \Delta x^4,
	\end{equation*}
	due to Proposition~\ref{prop:1d_eigval}.
\end{proof}

\subsection{Fully-discrete AF method}\label{sec:1d_cfl}
This section studies the stability of the fully-discrete method and gives a CFL condition.
When an SSP Runge--Kutta method \cite{Gottlieb_2001_Strong_SR, Spiteri_2002_New_SJNA} is used for time integration, the amplification matrix of the fully-discrete AF method is
\begin{equation}\label{eq:amplification_rk}
	\bm{\mathcal{G}} = R(\Delta t\bm{G}),
\end{equation}
where $\Delta t$ is the time step size, with $\bm{G}$ defined in \eqref{eq:1d_fourier_ode}.
The stability polynomials $R(z)$ for SSP-RK3, RK4, and SSP-RK54 are, respectively,
$R(z) = 1 + z + \frac12z^2 + \frac16z^3$, $R(z) = 1 + z + \frac12z^2 + \frac16z^3 + \frac{1}{24}z^4$, and
$R(z) = 1 + z + \frac12z^2 + \frac16z^3 + \frac{1}{24}z^4 + \eta z^5$ with $\eta=0.004477718303076$.
Define the P\'eclet number $\texttt{Pe} = \frac{c\Delta x}{\kappa_0} \in [0,+\infty)$, and consider the following CFL number
\begin{equation*}
	C_{\texttt{CFL}} = {\Delta t}\left( \frac{c}{\Delta x} + \frac{\kappa_0}{\Delta x^2} \right).
\end{equation*}
With this definition, $\Delta t\bm{G}$ can be parametrized as
\begin{equation*}
	\Delta t\bm{G} = C_{\texttt{CFL}}\left[ r \bm{G}_c + (1-r) \bm{G}_d \right],~
	r = \dfrac{\texttt{Pe}}{1 + \texttt{Pe}} \in [0, 1].
\end{equation*}
The two eigenvalues of $\bm{\mathcal{G}}$ can be expressed as a function of $C_{\texttt{CFL}}$, $r$, and $\xi = \omega \Delta x \in [-\pi, \pi]$.
The numerically estimated $C_{\texttt{CFL}}$ for stability is computed by grid search with different values of $r$ and $\xi$.
The results are summarized in Table~\ref{tab:1d_stability}, with a comparison with the LDG method \cite{Cockburn_1998_local_SJNA}.

When using the RK4 method, our results are consistent with the study in \cite{Watkins_2016_numerical_CF} for the LDG method.
For the same order of accuracy, one observes that the proposed $4$th-order AF method has a larger stable CFL number than the $P^3$ LDG method.

\begin{table}[htpb!]
	\centering
	\begin{tabular}{ccccccc}
		\toprule
		& \multicolumn{3}{c}{Pure advection $C_{\texttt{CFL}}=c\Delta t/\Delta x$}
		& \multicolumn{3}{c}{} \\
		\cmidrule(lr){2-4}
		Scheme
		& SSP-RK3
		& SSP-RK54
		& RK4
		& 
		& 
		&  \\
		\midrule
		AF4
		& 0.773
		& 1.049
		& 0.912
		& 
		& 
		&  \\
		$P^1$ LDG
		& 0.409
		& 0.662
		& 0.464
		& 
		& 
		&  \\
		$P^2$ LDG
		& 0.209
		& 0.343
		& 0.235
		& 
		& 
		&  \\
		$P^3$ LDG
		& 0.130
		& 0.215
		& 0.145
		& 
		& 
		&  \\
		$P^4$ LDG
		& 0.089
		& 0.149
		& 0.100
		& 
		& 
		&  \\
		\midrule
		&&
		\multicolumn{3}{c}{Pure diffusion $C_{\texttt{CFL}}=\kappa_0\Delta t/\Delta x^2$}
		&& \\
		\midrule
		& \multicolumn{3}{c}{Alternating flux for LDG}
		& \multicolumn{3}{c}{Central flux for LDG} \\
		\cmidrule(lr){2-4}
		\cmidrule(lr){5-7}
		Scheme
		& SSP-RK3
		& SSP-RK54
		& RK4
		& SSP-RK3
		& SSP-RK54
		& RK4 \\
		\midrule
		AF4
		& 0.157
		& 0.333
		& 0.174
		& 0.157
		& 0.333
		& 0.174 \\
		$P^1$ LDG
		& 0.069
		& 0.148
		& 0.077
		& 0.157
		& 0.333
		& 0.174 \\
		$P^2$ LDG
		& 0.016
		& 0.035
		& 0.018
		& 0.038
		& 0.081
		& 0.042 \\
		$P^3$ LDG
		& 0.005
		& 0.012
		& 0.006
		& 0.014
		& 0.030
		& 0.015 \\
		$P^4$ LDG
		& 0.002
		& 0.005
		& 0.002
		& 0.006
		& 0.013
		& 0.007 \\
		\bottomrule
	\end{tabular}
	\caption{Maximum stable CFL numbers of the $4$th-order AF method (AF4) and the LDG methods with different polynomial degrees.
		The LDG methods use an upwind flux for the convection part and alternating/central flux for the diffusion part.
	}
	\label{tab:1d_stability}
\end{table}

We are also interested in the total number of floating-point operations (FLOPs).
For the linear convection--diffusion equation, Table~\ref{tab:af4_ldg_p3_flops} lists the estimated FLOP counts required for one right-hand-side evaluation per cell by the AF and LDG methods.

\begin{table}[htb!]
	\centering
	\begin{tabular}{ccccc}
		\toprule
		Method
		&
		Convection
		&
		Diffusion
		&
		Coupling
		&
		Total
		\\
		\midrule
		AF4
		& 15
		& 20
		& 2
		& 37
		\\
		$P^3$ LDG
		& 24
		& 44
		& 4
		& 72
		\\
		\bottomrule
	\end{tabular}
	\caption{
		Estimated FLOP counts for one evaluation of the
		1D AF and $P^3$ LDG (with Legendre modal basis) spatial
		right-hand sides per cell.
		Coupling means the addition of the convection and diffusion parts of all DoFs.
	}
	\label{tab:af4_ldg_p3_flops}
\end{table}


\section{Multi-D formulation}\label{sec:multid_formulation}
\begin{figure}[hptb!]
	\centering
	\includegraphics[width=0.35\linewidth]{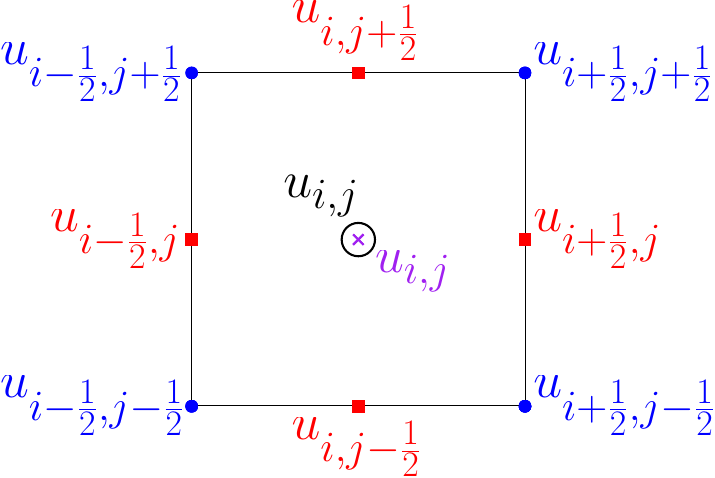}
	\caption{The DoFs of $u$: cell average (circle), face-centered values (squares), values at corners (dots). Note that the cell-centered point value $u_{i,j}$ (cross) is used in constructing the scheme, but does not belong to the DoFs.}
	\label{fig:2d_af_dofs}
\end{figure}

The 1D AF method \eqref{eq:1d_af_dc} can be extended to solve \eqref{eq:cde_conservative} on Cartesian meshes as follows.
The 3D formulation is presented first, followed by a detailed description of the 2D case.
Consider a uniform Cartesian mesh with $N_1\times N_2\times N_3$ cells
$I_{ijk} = [x_{\xl}, x_{\xr}]\times [y_{\yl}, y_{\yr}]\times [z_{\zl}, z_{\zr}]$ with cell centers $((x_{\xl} + x_{\xr})/2, (y_{\yl} + y_{\yr})/2, (z_{\zl} + z_{\zr})/2)$ and cell sizes $\Delta x, \Delta y, \Delta z$.
The DoFs consist of the cell averages, face-centered, edge-centered, and corner point values of the numerical solution $u_h(\bx, t)$, defined as
\begin{equation*}
	\overline{u}_{\bm{i}}(t) = \frac{1}{\Delta x\Delta y\Delta z}\int_{I_{ijk}} u_h(\bx,t) ~\dd \bx,\quad
	u_{\bm{\zeta}}(t) = u_h(\bx_{\bm{\zeta}}, t),
\end{equation*}
with $\bx=(x,y,z)=(x_1,x_2,x_3)$ and $\bm{\zeta} \in \{ (i+\frac12,j,k), (i,j+\frac12,k), (i,j,k+\frac12), (i+\frac12,j+\frac12,k), (i+\frac12,j,k+\frac12), (i,j+\frac12,k+\frac12), (\xr,\yr,\zr) \}$ being the locations of the point values.
Figure~\ref{fig:2d_af_dofs} shows the locations of the DoFs in the 2D case.

Define $\bm{i}=(i,j,k)$, $\bm{e}_1=(1,0,0)$, $\bm{e}_2=(0,1,0)$, $\bm{e}_3=(0,0,1)$.
The AF method can be written as follows,
\begin{align}
	&\frac{\dd \bar{u}_{\bm{i}}}{\dd t} + \sum_{\ell=1}^3\frac{1}{\Delta x_{\ell}}\left( \hat{f}_{\ell,\bm{i}+\frac12\bm{e}_\ell} - \hat{f}_{\ell,\bm{i}-\frac12\bm{e}_\ell} \right) = \sum_{\ell=1}^3\frac{1}{\Delta x_{\ell}} \left( \hat{a}_{\ell,\bm{i}+\frac12\bm{e}_\ell} - \hat{a}_{\ell,\bm{i}-\frac12\bm{e}_\ell} \right), \nonumber\\
	&\hat{f}_{\ell,\bm{i}+\frac12\bm{e}_\ell} = \overline{ f_{\ell}(u) }_{\bm{i}+\frac12\bm{e}_\ell},
	\quad \hat{a}_{\ell,\bm{i}+\frac12\bm{e}_\ell} = \overline{\left(\pd{a(u)}{x_\ell}\right)}_{\bm{i}+\frac12\bm{e}_\ell}, \nonumber\\
	&\left(\pd{a(u)}{x_\ell}\right)_{\bm{i}+\frac12\bm{e}_\ell} = D^{(\ell),\texttt{4th}}_{c} (a(u))_{\bm{i}+\frac12\bm{e}_\ell}, \nonumber\\
	&\frac{\dd u_{\bm{\zeta}} }{\dd t} + \sum_{\ell=1}^3 \left(D^{(\ell), \texttt{4th}}_{+} (f^{+}_{\ell})_{\bm{\zeta}} + D^{(\ell), \texttt{4th}}_{-}(f^{-}_{\ell})_{\bm{\zeta}} \right) = \sum_{\ell=1}^3 D^{(\ell),2,\texttt{4th}} (a(u))_{\bm{\zeta}}, \label{eq:2d_af_dc}
\end{align}
where the face average is obtained by the tensor-product Simpson's rule, e.g.,
\begin{align*}
	\overline{z}_{\xr,j,k} =&\ \frac{1}{36}\big( (z_{\xr,\yl,\zl} + z_{\xr,\yr,\zl} + z_{\xr,\yl,\zr} + z_{\xr,\yr,\zr}) \\
	&+ 4(z_{\xr,\yl,k} + z_{\xr,\yr,k} + z_{\xr,j,\zl} + z_{\xr,j,\zr}) \\
	&+ 16z_{\xr,j,k}\big).
\end{align*}
The finite difference operators $D^{(\ell),\texttt{4th}}_{\pm}$, $D^{(\ell),\texttt{4th}}_{c}$, and $D^{(\ell),2,\texttt{4th}}$ are applied along the $x_{\ell}$-direction, based on \eqref{eq:1d_upwind_ops_4p}-\eqref{eq:1d_upwind_ops_4m}, \eqref{eq:1d_diff_central_4th}, and \eqref{eq:1d_diff_compact_a}, respectively.

The 2D case is detailed below.
The $4$th-order biased finite difference operators $D^{(1),\texttt{4th}}_{\pm}$ are given by
\begin{align*}
	(D_{+}^{(1),\texttt{4th}}z)_{i+\frac12,j} &= \frac{1}{3\Delta x}\left( z_{i-\frac12,j} - 6z_{i,j} + 3z_{i+\frac12,j} + 2z_{i+1,j} \right), \\
	(D_{-}^{(1),\texttt{4th}}z)_{i+\frac12,j} &= \frac{1}{3\Delta x}\left( -2z_{i,j} - 3z_{i+\frac12,j} + 6z_{i+1,j} - z_{i+\frac32,j} \right), \\
	(D_{+}^{(2),\texttt{4th}}z)_{i,j+\frac12} &= \frac{1}{3\Delta y}\left( z_{i,j-\frac12} - 6z_{i,j} + 3z_{i,j+\frac12} + 2z_{i,j+1} \right), \\
	(D_{-}^{(2),\texttt{4th}}z)_{i,j+\frac12} &= \frac{1}{3\Delta y}\left( -2z_{i,j} - 3z_{i,j+\frac12} + 6z_{i,j+1} - z_{i,j+\frac32} \right).
\end{align*}
Note that the cell-centered value is obtained through
\begin{align*}
	z_{i,j} = \frac{1}{16}\Big[
	36\overline{z}_{i,j} &- 4\left(z_{\xl,j}+z_{\xr,j}+z_{i,\yl}+z_{i,\yr}\right)\nonumber \\
	&- \left(z_{\xl,\yl} + z_{\xr,\yl} + z_{\xl,\yr} + z_{\xr,\yr}\right)
	\Big].
\end{align*}
The central finite difference operator is
\begin{equation*}
	(D_{c}^{(\ell),\texttt{4th}} z)_{i,j} = \frac12\left( (D_{+}^{(\ell),\texttt{4th}}z)_{i,j} + (D_{-}^{(\ell),\texttt{4th}}z)_{i,j} \right).
\end{equation*}
The compact discretization for the pure second partial derivative in the $x$-direction can be expressed as
\begin{align*}
	\left(\dfrac{\partial^2 a(u)}{\partial x^2}  \right)_{i+\frac12,j+\frac12} =&\
	\frac{1}{\Delta x^2}\Big(- \frac{1}{3}a_{i-\frac12,j+\frac12} + \frac{16}{3}a_{i,j+\frac12} \\
	& - 10 a_{i+\frac12,j+\frac12} + \frac{16}{3}a_{i+1,j+\frac12} - \frac{1}{3}a_{i+\frac32,j+\frac12}\Big).
\end{align*}
The remaining formulas follow by half-index shifts in $i$ and/or $j$.

\section{Maximum-principle-preserving limiting and oscillation control}\label{sec:limiting}
This section constructs an MPP limiting for the proposed AF method, which aims at preserving the maximum principle at the discrete level.
The smoothness-indicator-based blending is also introduced for suppressing oscillations.

\begin{definition}
	An AF method is MPP if starting from cell averages and point values satisfying the bound $[u_{\min},u_{\max}]$, the cell averages and point values still stay in the bound at the next time step.
\end{definition}

For the cell average, the basic idea is to blend the high-order AF method with a low-order MPP scheme, similar to the idea in \cite{Duan_2025_Active_SJSC,Duan_2025_Active},
which amounts to using blended numerical fluxes, such that the local conservation is maintained.
However, a different convex limiting method is proposed, by adapting the monolithic convex limiting \cite{Kuzmin_2020_Monolithic_CMAME} to convection fluxes and diffusion potentials.
For the point value, this paper uses a direct post-processing as no conservation is required.
Since the SSP-RK54 is used for the time integration, which is a convex combination of the forward-Euler stages, it is sufficient to consider the forward Euler step and show the MPP property.

\subsection{1D MPP AF method}\label{sec:1d_pp}
\subsubsection{Monolithic limiting for cell averages}\label{sec:1d_limiting_average}

This section constructs monolithic convex limiting of convection fluxes and diffusion potentials for the cell-average update. The high-order correction of the convection part is associated with each cell interface, whereas that of the diffusion part is associated with each cell-centered numerical potential. The two families of limiting parameters are evaluated simultaneously from the same low-order states. In this way, local conservation of the convection part and double conservation of the diffusion part are both maintained.

Consider the following low-order finite volume scheme
\begin{align}\label{eq:1d_av_pp}
	&\bar{u}_i^{\texttt{L}}
	=
	\bar{u}_i^n
	-\frac{\Delta t^n}{\Delta x}
	\left(
	\hat{f}^{\texttt{L}}_{\xr}
	-\hat{f}^{\texttt{L}}_{\xl}
	\right)
	+\frac{\Delta t^n}{\Delta x}
	\left(
	\hat{a}^{\texttt{L}}_{\xr}
	-\hat{a}^{\texttt{L}}_{\xl}
	\right),
	\\
	&\hat{f}^{\texttt{L}}_{\xr}
	=
	f^{\texttt{LLF}}
	\left(
	\bar{u}_{i}^n,\bar{u}_{i+1}^n
	\right)
	=
	\frac12
	\left(
	f(\bar{u}_{i}^n)
	+f(\bar{u}_{i+1}^n)
	\right)
	-\frac{\alpha_{\xr}^n}{2}
	\left(
	\bar{u}_{i+1}^n-\bar{u}_{i}^n
	\right),
	\nonumber\\
	&\hat{a}^{\texttt{L}}_{\xr}
	=
	\frac{1}{\Delta x}
	\left(
	a(\bar{u}_{i+1}^n)-a(\bar{u}_{i}^n)
	\right).
	\nonumber
\end{align}
Here the viscosity coefficient is defined by
\begin{equation}\label{eq:1d_local_wave_speed}
	\alpha_{\xr}^n
	=
	\max_{
		u\in
		\left[
		\min\{\bar{u}_{i}^n,\bar{u}_{i+1}^n\},~
		\max\{\bar{u}_{i}^n,\bar{u}_{i+1}^n\}
		\right]
	}
	\abs{f'(u)}.
\end{equation}
For the diffusion part, define $\kappa_{\max} = \max_{u\in[u_{\min},u_{\max}]}\kappa(u)$.
Let $\mu=\frac{\Delta t^n}{\Delta x}$, $\nu=\frac{\Delta t^n}{\Delta x^2}$.
The low-order scheme \eqref{eq:1d_av_pp} can be rewritten as
\begin{equation}\label{eq:1d_low_directional_decomposition}
	\bar{u}_i^{\texttt{L}}
	=
	\frac12 u_i^{\texttt{L},+}
	+
	\frac12 u_i^{\texttt{L},-},
\end{equation}
where the intermediate states are
\begin{subequations}
	\begin{align}
		u_i^{\texttt{L},+}
		={}&
		\bar{u}_i^n
		-2\mu
		\left(
		\hat{f}_{\xr}^{\texttt{L}}
		-f(\bar{u}_i^n)
		\right)
		+2\frac{\Delta t^n}{\Delta x}
		\hat{a}_{\xr}^{\texttt{L}},
		\label{eq:1d_low_directional_plus}
		\\
		u_i^{\texttt{L},-}
		={}&
		\bar{u}_i^n
		-2\mu
		\left(
		f(\bar{u}_i^n)
		-\hat{f}_{\xl}^{\texttt{L}}
		\right)
		-2\frac{\Delta t^n}{\Delta x}
		\hat{a}_{\xl}^{\texttt{L}}.
		\label{eq:1d_low_directional_minus}
	\end{align}
\end{subequations}

\begin{lemma}\label{lem:1d_mpp_scheme}
	If $\bar{u}_i^n \in [u_{\min},u_{\max}]$ for all $i$,
	then the intermediate states satisfy $u_i^{\texttt{L},\pm}\in[u_{\min},u_{\max}]$,
	under the following CFL condition
	\begin{equation}\label{eq:1d_av_pp_dt}
		\Delta t^n
		\leqslant
		\frac{C_{\texttt{CFL}}}{
			\displaystyle
			\max_i
			\left(
			\frac{\alpha_{\xr}^n}{\Delta x}
			+
			\frac{\kappa_{\max}}{\Delta x^2}
			\right)
		},
		\quad
		C_{\texttt{CFL}}\leqslant\frac12.
	\end{equation}
	
\end{lemma}

\begin{proof}
	We only prove that
	$u_i^{\texttt{L},+}\in[u_{\min},u_{\max}]$ as the proof for
	$u_i^{\texttt{L},-}$ is analogous.
	If $\bar{u}_{i+1}^n=\bar{u}_i^n$, consistency immediately gives $u_i^{\texttt{L},+}=\bar{u}_i^n\in[u_{\min},u_{\max}]$.
	Otherwise, introduce
	\begin{equation*}
		\sigma_{\xr}^{\texttt{c}}
		=
		\frac{
			f(\bar{u}_{i+1}^n)-f(\bar{u}_i^n)
		}{
			\bar{u}_{i+1}^n-\bar{u}_i^n
		},\quad
		\sigma_{\xr}^{\texttt{d}}
		=
		\frac{
			a(\bar{u}_{i+1}^n)-a(\bar{u}_i^n)
		}{
			\bar{u}_{i+1}^n-\bar{u}_i^n
		}.
	\end{equation*}
	By the definition in \eqref{eq:1d_local_wave_speed} and the mean value theorem, one has $\abs{\sigma_{\xr}^{\texttt{c}}} \leqslant \alpha_{\xr}^n$.
	Since $\kappa(u) \in [0,\kappa_{\max}]$ on $u \in [u_{\min},u_{\max}]$, one obtains $\sigma_{\xr}^{\texttt{d}} \in [0, \kappa_{\max}]$.	
	Rewrite \eqref{eq:1d_low_directional_plus} as
	\begin{align*}
		u_i^{\texttt{L},+}
		={}&
		\bar{u}_i^n
		+
		\left[
		\mu
		\left(
		\alpha_{\xr}^n-\sigma_{\xr}^{\texttt{c}}
		\right)
		+
		2\nu\sigma_{\xr}^{\texttt{d}}
		\right]
		\left(
		\bar{u}_{i+1}^n-\bar{u}_i^n
		\right)
		\\
		={}&
		(1-r)\bar{u}_i^n
		+
		r\bar{u}_{i+1}^n,
	\end{align*}
	where $r =	\mu \left(\alpha_{\xr}^n-\sigma_{\xr}^{\texttt{c}}\right) + 2\nu\sigma_{\xr}^{\texttt{d}}$.
	Thanks to the bound of $\sigma_{\xr}^{\texttt{c}}$ and $\sigma_{\xr}^{\texttt{d}}$, one has $r\geqslant 0$ and
	\begin{equation*}
		r
		\leqslant
		2\mu\alpha_{\xr}^n
		+
		2\nu\kappa_{\max}
		=
		2\Delta t^n
		\left(
		\frac{\alpha_{\xr}^n}{\Delta x}
		+
		\frac{\kappa_{\max}}{\Delta x^2}
		\right)
		\leqslant
		2C_{\texttt{CFL}}
		\leqslant1.
	\end{equation*}
	Consequently,
	$u_i^{\texttt{L},+}$ is a convex combination of
	$\bar{u}_i^n$ and $\bar{u}_{i+1}^n$,
	and it follows that $u_i^{\texttt{L},+}\in[u_{\min},u_{\max}]$.
\end{proof}

To preserve the discrete maximum principle, while simultaneously guaranteeing the conservation of the convection part and also doubly conservative structure of the diffusion discretization simultaneously,
the limited scheme is defined as
\begin{align}
	\bar{u}_i^{\texttt{Lim}} &= \bar{u}_i^n - \mu \left( \hat{f}^{\texttt{Lim}}_{\xr}  - \hat{f}^{\texttt{Lim}}_{\xl} \right) + \frac{\Delta t^n}{\Delta x} \left( \hat{a}^{\texttt{Lim}}_{\xr}  - \hat{a}^{\texttt{Lim}}_{\xl} \right) \nonumber \\
	&= \bar{u}_i^n - \mu \left( \hat{f}^{\texttt{Lim}}_{\xr}  - \hat{f}^{\texttt{Lim}}_{\xl} \right) + \nu \left( \hat{a}^{\texttt{Lim}}_{i+1} - 2\hat{a}^{\texttt{Lim}}_{i} + \hat{a}^{\texttt{Lim}}_{i-1} \right), \label{eq:1d_limited_average}
\end{align}
where the limited flux and numerical diffusion potential are given by
\begin{align*}
	\hat{f}_{\xr}^{\texttt{Lim}}
	&=
	\hat{f}_{\xr}^{\texttt{L}}
	+
	\theta_{\xr}^{\texttt{c}}\Delta\hat{f}_{\xr},
	\qquad
	\\
	\hat{a}_i^{\texttt{Lim}}
	&=
	\hat{a}_i^{\texttt{L}}
	+
	\theta_{i}^{\texttt{d}}\Delta\hat{a}_i,
\end{align*}
with $\theta_{\xr}^{\texttt{c}}, \theta_{i}^{\texttt{d}} \in [0,1]$.
The corrections for the convection flux and diffusion potential are
\begin{equation*}
	\Delta\hat{f}_{\xr}
	=
	\hat{f}_{\xr}^{\texttt{H}}
	-\hat{f}_{\xr}^{\texttt{L}},
	\quad
	\Delta\hat{a}_i
	=
	\hat{a}_i^{\texttt{H}}
	-\hat{a}_i^{\texttt{L}},
\end{equation*}
with the high-order numerical flux in the AF method $\hat{f}_{\xr}^{\texttt{H}} = f(u_{\xr})$, and the low- and high-order numerical potentials
\begin{equation*}
	\hat{a}_i^{\texttt{L}}
	=
	a(\bar{u}_i^n),\quad
	\hat{a}_i^{\texttt{H}}
	=
	\frac16
	\left(
	-a(u_{\xl}^n)
	+8a(u_i^n)
	-a(u_{\xr}^n)
	\right).
\end{equation*}
%
The two limited states can be expressed as
\begin{align*}
	&\bar{u}_i^{\texttt{Lim}} = \frac12 u_i^{\texttt{Lim},+} + \frac12 u_i^{\texttt{Lim},-}, \\
	&u_i^{\texttt{Lim},+}
	=
	\frac13\Big[
	\left(u_i^{\texttt{L},+} - 6\mu\theta^{\texttt{c}}_{\xr}\Delta\hat{f}_{\xr} \right)
	+
	\left( u_i^{\texttt{L},+} + 6\nu\theta^{\texttt{d}}_{i+1}\Delta\hat{a}_{i+1} \right)
	+
	\left(u_i^{\texttt{L},+} - 6\nu\theta^{\texttt{d}}_i\Delta\hat{a}_i \right)
	\Big],
	\\
	&u_i^{\texttt{Lim},-}
	=
	\frac13\Big[
	\left( u_i^{\texttt{L},-}
	+6\mu\theta^{\texttt{c}}_{\xl}\Delta\hat{f}_{\xl} \right)
	+
	\left( u_i^{\texttt{L},-}
	+6\nu\theta^{\texttt{d}}_{i-1}\Delta\hat{a}_{i-1} \right)
	+
	\left( u_i^{\texttt{L},-}
	-6\nu\theta^{\texttt{d}}_i\Delta\hat{a}_i \right)
	\Big].
\end{align*}
To preserve the MPP property, it suffices to have each part of the decomposition of $u_i^{\texttt{Lim},\pm}$ stay in $[u_{\min},u_{\max}]$.
For $z^{\texttt{L}}\in[u_{\min},u_{\max}]$, define
\begin{equation*}
	\Theta
	\left(
	z^{\texttt{L}},z^{\texttt{H}}
	\right)
	=
	\begin{cases}
		\dfrac{z^{\texttt{L}}-u_{\min}}
		{z^{\texttt{L}}-z^{\texttt{H}}},
		&
		z^{\texttt{H}}<u_{\min},
		\\[1.1em]
		\dfrac{z^{\texttt{L}}-u_{\max}}
		{z^{\texttt{L}}-z^{\texttt{H}}},
		&
		z^{\texttt{H}}>u_{\max},
		\\[1.1em]
		1,
		&
		z^{\texttt{H}}\in[u_{\min},u_{\max}].
	\end{cases}
\end{equation*}
By construction, one has $(1-\vartheta)z^{\texttt{L}} + \vartheta z^{\texttt{H}} \in [u_{\min},u_{\max}]$, for all $\vartheta \in [0, \Theta(z^{\texttt{L}},z^{\texttt{H}})]$.
The limiting parameter for the convection flux is chosen as
\begin{align}
	\theta^{\texttt{c}}_{\xr}
	=
	\min\Big\{
	&
	\Theta\left(
	u_i^{\texttt{L},+},
	u_i^{\texttt{L},+}
	-6\mu\Delta\hat{f}_{\xr}
	\right),
	\Theta\left(
	u_{i+1}^{\texttt{L},-},
	u_{i+1}^{\texttt{L},-}
	+6\mu\Delta\hat{f}_{\xr}
	\right)
	\Big\},
	\label{eq:1d_convection_limiting_parameter}
\end{align}
which considers the two cells sharing the same interface.
The limiting parameter for the numerical potential in the cell $I_i$ is chosen as
\begin{align}
	\theta^{\texttt{d}}_i
	=
	\min\Big\{
	&
	\Theta\left(
	u_{i-1}^{\texttt{L},+},
	u_{i-1}^{\texttt{L},+}
	+6\nu\Delta\hat{a}_i
	\right),
	\Theta\left(
	u_i^{\texttt{L},+},
	u_i^{\texttt{L},+}
	-6\nu\Delta\hat{a}_i
	\right),
	\nonumber\\
	&
	\Theta\left(
	u_i^{\texttt{L},-},
	u_i^{\texttt{L},-}
	-6\nu\Delta\hat{a}_i
	\right),
	\Theta\left(
	u_{i+1}^{\texttt{L},-},
	u_{i+1}^{\texttt{L},-}
	+6\nu\Delta\hat{a}_i
	\right)
	\Big\},
	\label{eq:1d_diffusion_limiting_parameter}
\end{align}
which considers the four intermediate states affected by $\Delta\hat{a}_i$.
Note that $\theta^{\texttt{c}}_{\xr}$ and $\theta^{\texttt{d}}_i$ are evaluated from the same low-order states in one monolithic step.
The main result of the MPP limiting for the cell average is summarized as follows.
\begin{proposition}\label{prop:1d_joint_mpp_limiter}
	Suppose that the low-order intermediate states defined by
	\eqref{eq:1d_low_directional_plus}-\eqref{eq:1d_low_directional_minus} satisfy
	$u_i^{\texttt{L},\pm}\in[u_{\min},u_{\max}]$ for all $i$,
	and the CFL condition \eqref{eq:1d_av_pp_dt} holds.
	Then the scheme
	\eqref{eq:1d_limited_average}, equipped with the limiting parameters
	\eqref{eq:1d_convection_limiting_parameter} and
	\eqref{eq:1d_diffusion_limiting_parameter}, is MPP, i.e., $\bar{u}_i^{\texttt{Lim}}\in[u_{\min},u_{\max}]$.
	Moreover, the convection discretization remains locally conservative, and the diffusion discretization remains doubly conservative.
\end{proposition}

\begin{remark}
	If $\theta^{\texttt{c}}_{\xr}=1$, $\theta^{\texttt{d}}_{i}=1$, then $\hat{f}_{\xr}^{\texttt{Lim}} = \hat{f}_{\xr}$, $\hat{a}_i^{\texttt{Lim}} = \hat{a}_i^{\texttt{H}}$,
	and the original high-order cell-average update in the AF method is recovered exactly.
\end{remark}

\subsubsection{Post-processing for point value}\label{sec:1d_limiting_point}
Since there is no conservation requirement for this kind of DoF, this paper enforces the point value as
\begin{equation*}
	u_{\xr}^{\texttt{Lim}} = \min\{\max\{u_{\xr}^{\texttt{H}}, u_{\min}\}, u_{\max}\},
\end{equation*}
where $u_{\xr}^{\texttt{H}}$ is obtained by the original high-order AF method.

\subsubsection{Suppressing oscillations}\label{sec:1d_so}
The MPP property alone is insufficient to suppress oscillations, therefore, an additional blending procedure is introduced.
Following \cite{Duan_2025_Active_SJSC}, a smoothness indicator is computed as
\begin{equation*}
	\tilde{s}_{i} = \dfrac{\abs{\bar{u}_{i-1} - 2\bar{u}_{i} + \bar{u}_{i+1}}}{\abs{\bar{u}_{i-1}} + 2\abs{\bar{u}_{i}} + \abs{\bar{u}_{i+1}} + 10^{-13}},
\end{equation*}
with the smoothing $s_{i} = \frac16\left(\tilde{s}_{i-1} + 4\tilde{s}_{i} + \tilde{s}_{i+1}\right)$.
For the cell average, the blending coefficient for the convection part at the cell interface is defined as
\begin{equation*}
	\theta^{\texttt{c,s}}_{i+\frac12} = \exp\left( - \gamma s_{\xr} \right) \in (0,1], \quad
	s_{\xr} = \frac12\left( s_{i} + s_{i+1} \right),
\end{equation*}
and the final numerical flux is given by
\begin{equation*}
	\hat{f}_{\xr} = \theta^{\texttt{c,s}}_{i+\frac12} \hat{f}_{\xr}^{\texttt{H}} + (1-\theta^{\texttt{c,s}}_{i+\frac12}) \hat{f}_{\xr}^{\texttt{L}}.
\end{equation*}
For the diffusion part, the blending coefficient is
\begin{equation*}
	\theta^{\texttt{d,s}}_{i} = \exp\left( - \gamma s_{i} \right) \in (0,1],
\end{equation*}
and the final numerical potential is given by
\begin{equation*}
	\hat{a}_{i} = \theta^{\texttt{d,s}}_{i} \hat{a}_{i}^{\texttt{H}} + (1-\theta^{\texttt{d,s}}_{i}) \hat{a}_{i}^{\texttt{L}}.
\end{equation*}
For the point value, the high-order AF right-hand side is blended with the low-order LLF one using $\theta^{\texttt{c,s}}_{i+\frac12}$, as there is no conservation requirement (see \cite{Duan_2025_Active_SJSC} for more details).
The numerical tests will show that such smoothness-indicator-based limiting can suppress oscillations well.
Note that such a blending is applied before the MPP limiting.

\subsection{Multi-D MPP AF method}\label{sec:3d_pp}
Only the MPP limiting procedure for the 3D case is presented here, since the 2D formulation can be obtained similarly.
Let $\mu_{\ell} = \frac{\Delta t^n}{\Delta x_{\ell}}$, $\nu_{\ell} = \frac{\Delta t^n}{\Delta x_{\ell}^2}$, $\ell=1,2,3$.
Consider the following low-order scheme
\begin{equation}\label{eq:multid_av_pp}
	\bar{u}_{\bm{i}}^{\texttt{L}}
	=
	\bar{u}_{\bm{i}}^{n}
	-
	\sum_{\ell=1}^3
	\mu_{\ell}
	\left(
	\hat{f}_{\ell,\bm{i}+\frac12\bm{e}_{\ell}}^{\texttt{L}}
	-
	\hat{f}_{\ell,\bm{i}-\frac12\bm{e}_{\ell}}^{\texttt{L}}
	\right)
	+
	\sum_{\ell=1}^3
	\mu_{\ell}
	\left(
	\hat{a}_{\ell,\bm{i}+\frac12\bm{e}_{\ell}}^{\texttt{L}}
	-
	\hat{a}_{\ell,\bm{i}-\frac12\bm{e}_{\ell}}^{\texttt{L}}
	\right),
\end{equation}
where
\begin{align*}
	&\hat{f}_{\ell,\bm{i}+\frac12\bm{e}_{\ell}}^{\texttt{L}}
	=
	f_{\ell}^{\texttt{LLF}}
	\left(
	\bar{u}_{\bm{i}}^n,
	\bar{u}_{\bm{i}+\bm{e}_{\ell}}^n
	\right)
	=
	\frac12
	\left(
	f_{\ell}(\bar{u}_{\bm{i}}^n)
	+
	f_{\ell}(\bar{u}_{\bm{i}+\bm{e}_{\ell}}^n)
	\right)
	-
	\frac{\alpha_{\ell,\bm{i}+\frac12\bm{e}_{\ell}}^n}{2}
	\left(
	\bar{u}_{\bm{i}+\bm{e}_{\ell}}^n
	-
	\bar{u}_{\bm{i}}^n
	\right),
	\\
	&\hat{a}_{\ell,\bm{i}+\frac12\bm{e}_{\ell}}^{\texttt{L}}
	=
	\frac{1}{\Delta x_{\ell}}
	\left(
	a(\bar{u}_{\bm{i}+\bm{e}_{\ell}}^n)
	-
	a(\bar{u}_{\bm{i}}^n)
	\right).
\end{align*}
The local viscosity coefficient for the convection part is defined by
\begin{equation*}
	\alpha_{\ell,\bm{i}+\frac12\bm{e}_{\ell}}^n
	=
	\max_{
		u\in
		\left[
		\min\left\{
		\bar{u}_{\bm{i}}^n,
		\bar{u}_{\bm{i}+\bm{e}_{\ell}}^n
		\right\},~
		\max\left\{
		\bar{u}_{\bm{i}}^n,
		\bar{u}_{\bm{i}+\bm{e}_{\ell}}^n
		\right\}
		\right]
	}
	\abs{f_{\ell}'(u)}.
\end{equation*}
For the diffusion part, define $\kappa_{\max} = \max_{u\in[u_{\min},u_{\max}]}\kappa(u)$.
The low-order scheme can be decomposed as
\begin{equation*}
	\bar{u}_{\bm{i}}^{\texttt{L}}
	=
	\frac16
	\sum_{\ell=1}^3
	\left(
	u_{\ell,\bm{i}}^{\texttt{L},+}
	+
	u_{\ell,\bm{i}}^{\texttt{L},-}
	\right),
\end{equation*}
where
\begin{subequations}
	\begin{align}
		u_{\ell,\bm{i}}^{\texttt{L},+}
		={}&
		\bar{u}_{\bm{i}}^n
		-
		6\mu_{\ell}
		\left(
		\hat{f}_{\ell,\bm{i}+\frac12\bm{e}_{\ell}}^{\texttt{L}}
		-
		f_{\ell}(\bar{u}_{\bm{i}}^n)
		\right)
		+
		6\mu_{\ell}
		\hat{a}_{\ell,\bm{i}+\frac12\bm{e}_{\ell}}^{\texttt{L}},
		\label{eq:multid_low_directional_plus}
		\\
		u_{\ell,\bm{i}}^{\texttt{L},-}
		={}&
		\bar{u}_{\bm{i}}^n
		-
		6\mu_{\ell}
		\left(
		f_{\ell}(\bar{u}_{\bm{i}}^n)
		-
		\hat{f}_{\ell,\bm{i}-\frac12\bm{e}_{\ell}}^{\texttt{L}}
		\right)
		-
		6\mu_{\ell}
		\hat{a}_{\ell,\bm{i}-\frac12\bm{e}_{\ell}}^{\texttt{L}}.
		\label{eq:multid_low_directional_minus}
	\end{align}
\end{subequations}
A sufficient CFL condition for
$u_{\ell,\bm{i}}^{\texttt{L},\pm}\in[u_{\min},u_{\max}]$ is
\begin{equation}\label{eq:multid_av_pp_dt}
	\Delta t^n
	\leqslant
	\frac{C_{\texttt{CFL}}}{
		\displaystyle
		\max_{\bm{i}, \ell=1,2,3}
		\left(
		\frac{
			\alpha_{\ell,\bm{i}+\frac12\bm{e}_{\ell}}^n
		}{
			\Delta x_{\ell}
		}
		+
		\frac{\kappa_{\max}}{\Delta x_{\ell}^2}
		\right)
	},
	\quad
	C_{\texttt{CFL}}\leqslant\frac16.
\end{equation}
The proof is similar to that of Lemma~\ref{lem:1d_mpp_scheme}.

\subsubsection{Monolithic limiting for cell averages}

As in the 1D case, the high-order correction of the convection part is
associated with each cell interface, whereas the high-order correction
of the diffusion part is associated with each cell-centered numerical
potential. The two types of correction are limited simultaneously using
the same low-order states.


Define $\Delta\hat{f}_{\ell,\bm{i}+\frac12\bm{e}_{\ell}}
=
\hat{f}_{\ell,\bm{i}+\frac12\bm{e}_{\ell}}^{\texttt{H}}
-
\hat{f}_{\ell,\bm{i}+\frac12\bm{e}_{\ell}}^{\texttt{L}}$,
and $\Delta\hat{a}_{\ell,\bm{i}}
=
\hat{a}_{\ell,\bm{i}}^{\texttt{H}}
-
\hat{a}_{\ell,\bm{i}}^{\texttt{L}}$,
where the high-order numerical flux $\hat{f}_{\ell,\bm{i}+\frac12\bm{e}_{\ell}}^{\texttt{H}}$ and numerical potential come from the AF method,
and the low-order numerical flux $\hat{f}_{\ell,\bm{i}+\frac12\bm{e}_{\ell}}^{\texttt{L}}$ and numerical potential are given by the LLF flux and $\hat{a}_{\ell,\bm{i}}^{\texttt{L}} = a(\bar{u}_{\bm{i}}^{n})$, respectively.
The limited convection flux and cell-centered numerical potential are
given by
\begin{align*}
	&\hat{f}_{\ell,\bm{i}+\frac12\bm{e}_{\ell}}^{\texttt{Lim}}
	=
	\hat{f}_{\ell,\bm{i}+\frac12\bm{e}_{\ell}}^{\texttt{L}}
	+
	\theta_{\ell,\bm{i}+\frac12\bm{e}_{\ell}}^{\texttt{c}}
	\Delta\hat{f}_{\ell,\bm{i}+\frac12\bm{e}_{\ell}},
	\\
	&\hat{a}_{\ell,\bm{i}}^{\texttt{Lim}}
	=
	\hat{a}_{\ell,\bm{i}}^{\texttt{L}}
	+
	\theta_{\ell,\bm{i}}^{\texttt{d}}
	\Delta\hat{a}_{\ell,\bm{i}},
\end{align*}
where $\theta_{\ell,\bm{i}+\frac12\bm{e}_{\ell}}^{\texttt{c}},
\theta_{\ell,\bm{i}}^{\texttt{d}}
\in[0,1]$.
Consequently, the limited cell-average update is
\begin{equation}
	\bar{u}_{\bm{i}}^{\texttt{Lim}}
	=
	\bar{u}_{\bm{i}}^n
	-
	\sum_{\ell=1}^3
	\mu_{\ell}
	\left(
	\hat{f}_{\ell,\bm{i}+\frac12\bm{e}_{\ell}}^{\texttt{Lim}}
	-
	\hat{f}_{\ell,\bm{i}-\frac12\bm{e}_{\ell}}^{\texttt{Lim}}
	\right)
	+
	\sum_{\ell=1}^3
	\nu_{\ell}
	\left(
	\hat{a}_{\ell,\bm{i}+\bm{e}_{\ell}}^{\texttt{Lim}}
	-
	2\hat{a}_{\ell,\bm{i}}^{\texttt{Lim}}
	+
	\hat{a}_{\ell,\bm{i}-\bm{e}_{\ell}}^{\texttt{Lim}}
	\right),
	\label{eq:multid_limited_average}
\end{equation}
which can be decomposed as $\bar{u}_{\bm{i}}^{\texttt{Lim}}
=
\frac16
\sum_{\ell=1}^3
\left(
u_{\ell,\bm{i}}^{\texttt{Lim},+}
+
u_{\ell,\bm{i}}^{\texttt{Lim},-}
\right)$,
with
\begin{align*}
	u_{\ell,\bm{i}}^{\texttt{Lim},+}
	={}&
	\frac13
	\Bigg[
	\left(
	u_{\ell,\bm{i}}^{\texttt{L},+}
	-
	18\mu_{\ell}
	\theta_{\ell,\bm{i}+\frac12\bm{e}_{\ell}}^{\texttt{c}}
	\Delta\hat{f}_{\ell,\bm{i}+\frac12\bm{e}_{\ell}}
	\right)
	\nonumber\\
	&\quad
	+
	\left(
	u_{\ell,\bm{i}}^{\texttt{L},+}
	+
	18\nu_{\ell}
	\theta_{\ell,\bm{i}+\bm{e}_{\ell}}^{\texttt{d}}
	\Delta\hat{a}_{\ell,\bm{i}+\bm{e}_{\ell}}
	\right)
	+
	\left(
	u_{\ell,\bm{i}}^{\texttt{L},+}
	-
	18\nu_{\ell}
	\theta_{\ell,\bm{i}}^{\texttt{d}}
	\Delta\hat{a}_{\ell,\bm{i}}
	\right)
	\Bigg], \\
	u_{\ell,\bm{i}}^{\texttt{Lim},-}
	={}&
	\frac13
	\Bigg[
	\left(
	u_{\ell,\bm{i}}^{\texttt{L},-}
	+
	18\mu_{\ell}
	\theta_{\ell,\bm{i}-\frac12\bm{e}_{\ell}}^{\texttt{c}}
	\Delta\hat{f}_{\ell,\bm{i}-\frac12\bm{e}_{\ell}}
	\right)
	\nonumber\\
	&\quad
	+
	\left(
	u_{\ell,\bm{i}}^{\texttt{L},-}
	+
	18\nu_{\ell}
	\theta_{\ell,\bm{i}-\bm{e}_{\ell}}^{\texttt{d}}
	\Delta\hat{a}_{\ell,\bm{i}-\bm{e}_{\ell}}
	\right)
	+
	\left(
	u_{\ell,\bm{i}}^{\texttt{L},-}
	-
	18\nu_{\ell}
	\theta_{\ell,\bm{i}}^{\texttt{d}}
	\Delta\hat{a}_{\ell,\bm{i}}
	\right)
	\Bigg].
\end{align*}
Using the function $\Theta$ defined in the 1D case, the limiting parameter for the convection flux is chosen as
\begin{align}
	\theta_{\ell,\bm{i}+\frac12\bm{e}_{\ell}}^{\texttt{c}}
	=
	\min\Big\{
	&
	\Theta\left(
	u_{\ell,\bm{i}}^{\texttt{L},+},
	u_{\ell,\bm{i}}^{\texttt{L},+}
	-
	18\mu_{\ell}
	\Delta\hat{f}_{\ell,\bm{i}+\frac12\bm{e}_{\ell}}
	\right),
	\nonumber\\
	&
	\Theta\left(
	u_{\ell,\bm{i}+\bm{e}_{\ell}}^{\texttt{L},-},
	u_{\ell,\bm{i}+\bm{e}_{\ell}}^{\texttt{L},-}
	+
	18\mu_{\ell}
	\Delta\hat{f}_{\ell,\bm{i}+\frac12\bm{e}_{\ell}}
	\right)
	\Big\}.
	\label{eq:multid_convection_limiting_parameter}
\end{align}
For the diffusion part, the parameter $\theta_{\ell,\bm{i}}^{\texttt{d}}$ is chosen as
\begin{align}
	\theta_{\ell,\bm{i}}^{\texttt{d}}
	=
	\min\Big\{
	&
	\Theta\left(
	u_{\ell,\bm{i}-\bm{e}_{\ell}}^{\texttt{L},+},
	u_{\ell,\bm{i}-\bm{e}_{\ell}}^{\texttt{L},+}
	+
	18\nu_{\ell}
	\Delta\hat{a}_{\ell,\bm{i}}
	\right),
	\nonumber\\
	&
	\Theta\left(
	u_{\ell,\bm{i}}^{\texttt{L},+},
	u_{\ell,\bm{i}}^{\texttt{L},+}
	-
	18\nu_{\ell}
	\Delta\hat{a}_{\ell,\bm{i}}
	\right),
	\nonumber\\
	&
	\Theta\left(
	u_{\ell,\bm{i}}^{\texttt{L},-},
	u_{\ell,\bm{i}}^{\texttt{L},-}
	-
	18\nu_{\ell}
	\Delta\hat{a}_{\ell,\bm{i}}
	\right),
	\nonumber\\
	&
	\Theta\left(
	u_{\ell,\bm{i}+\bm{e}_{\ell}}^{\texttt{L},-},
	u_{\ell,\bm{i}+\bm{e}_{\ell}}^{\texttt{L},-}
	+
	18\nu_{\ell}
	\Delta\hat{a}_{\ell,\bm{i}}
	\right)
	\Big\}.
	\label{eq:multid_diffusion_limiting_parameter}
\end{align}
The convection and diffusion limiting parameters are evaluated from the same low-order intermediate states in one monolithic step.
The result is summarized as follows.
\begin{proposition}\label{prop:multid_joint_mpp_limiter}
	Suppose that the low-order intermediate states defined by
	\eqref{eq:multid_low_directional_plus} and
	\eqref{eq:multid_low_directional_minus} satisfy $u_{\ell,\bm{i}}^{\texttt{L},\pm}\in[u_{\min},u_{\max}]$, $\ell=1,2,3$,
	and that the CFL condition \eqref{eq:multid_av_pp_dt} holds.
	Then the scheme \eqref{eq:multid_limited_average}, equipped with the
	limiting parameters
	\eqref{eq:multid_convection_limiting_parameter} and
	\eqref{eq:multid_diffusion_limiting_parameter}, is MPP, i.e., $\bar{u}_{\bm{i}}^{\texttt{Lim}}\in[u_{\min},u_{\max}]$.
	Moreover, the convection discretization remains locally
	conservative, and the diffusion discretization remains doubly
	conservative in each direction.
\end{proposition}

\subsubsection{Post-processing for point value}
The point value is directly enforced as
\begin{equation*}
	u_{\bm{\zeta}}^{\texttt{Lim}} = \min\{\max\{u_{\bm{\zeta}}^{\texttt{H}}, u_{\min}\}, u_{\max}\}.
\end{equation*}

\subsubsection{Suppressing oscillations}
The oscillation-control procedure in Section~\ref{sec:1d_so} can be applied direction by direction, so that its multidimensional details are omitted for brevity.

\section{Numerical experiments}\label{sec:results}
This section presents numerical experiments to assess the accuracy and MPP property of the proposed method and to compare it with the LDG method.
In the 2D and 3D plots, the numerical solutions are visualized on a refined mesh with half the mesh size, where the values at the grid points are the cell averages or point values on the original mesh.
For the convection--diffusion cases, the CFL number is chosen as $0.33/d$ for $d=1,2,3$.
The smoothness-indicator-based blending is only activated in Examples~\ref{ex:1d_strongly_degenerate} and \ref{ex:2d_strongly_degenerate}.

\begin{example}[Smooth linear problems]\label{ex:accuracy}
	To examine the convergence rates of the proposed AF method numerically, we solve \eqref{eq:cde_conservative} with $f_\ell=cu$ and $a(u) = \kappa_0u$, where $c=1$ and $\kappa_0=0.01$.
	For each $d=1,2,3$, take the exact solution
	\begin{equation*}
		u(\bm{x},t) = 1 + \exp(-4d\pi^2 \kappa_0 t)\sin(2\pi \sum_{\ell=1}^d (x_\ell - ct)),
	\end{equation*}
	and all the tests are computed until $T=0.2$ on a periodic computational domain $[0,1]^d$.
	
	The maximum CFL numbers for stability are experimentally verified, which are consistent with the analysis in Table~\ref{tab:1d_stability}.
	The numerical errors and convergence rates are shown in Figure~\ref{fig:1d_af_accuracy}.
	One observes that the cell average and point value are both $4$th-order accurate.
	
	\begin{figure}[htb!]
		\centering
		\includegraphics[width=0.32\linewidth]{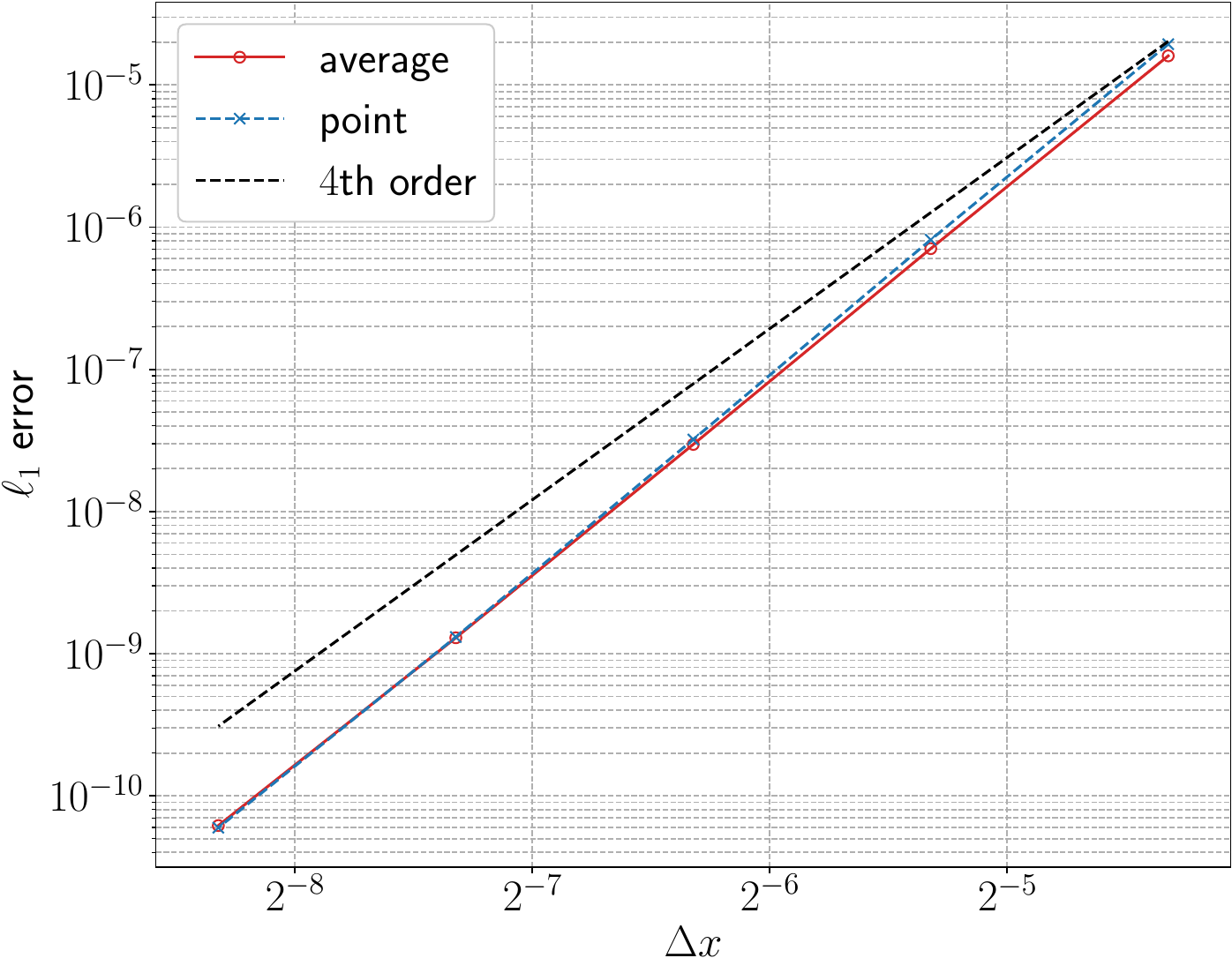}
		\includegraphics[width=0.32\linewidth]{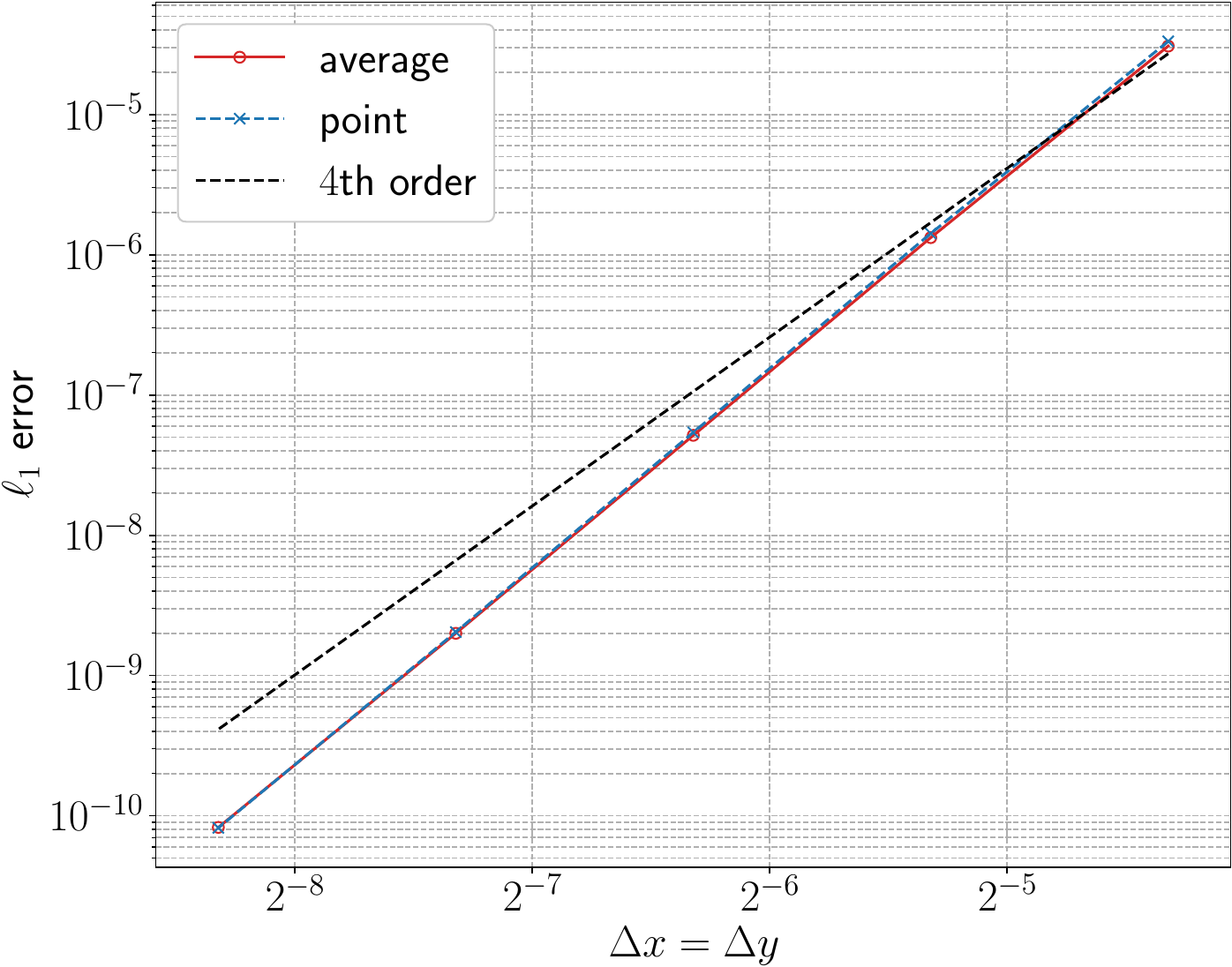}
		\includegraphics[width=0.32\linewidth]{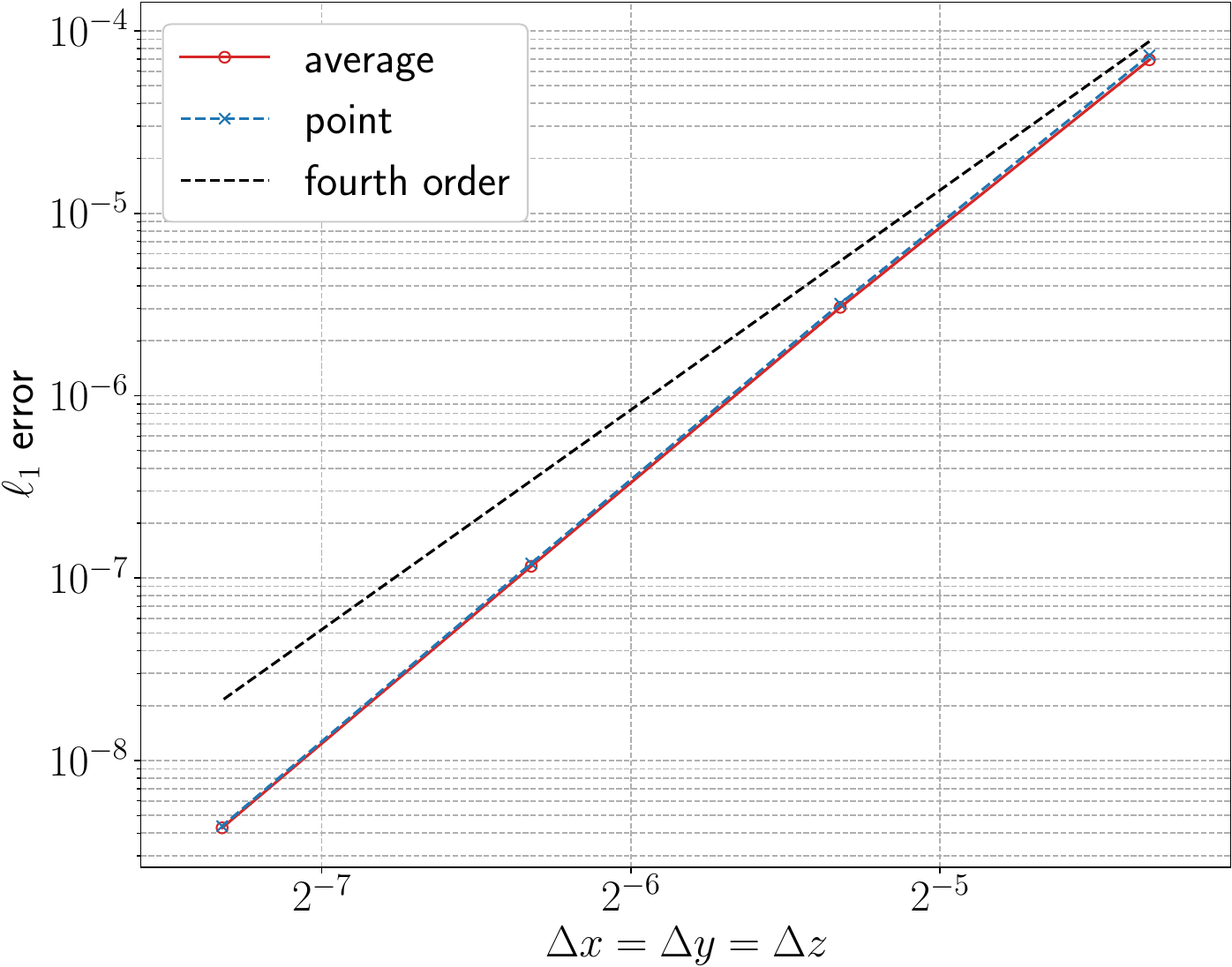}
		\caption{Example~\ref{ex:accuracy}. The errors and convergence rates of the proposed AF method.
			From left to right: 1D, 2D, and 3D cases.}
		\label{fig:1d_af_accuracy}
	\end{figure}
	
	Figure~\ref{fig:1d_comparison} compares the proposed AF method with the $P^3$ LDG method in 1D.
	Both methods are used with their own maximal stable CFL number, i.e., $0.33$ versus $0.012$ as shown in Table~\ref{tab:1d_stability}.
	To achieve the same level of accuracy, the AF method uses fewer DoFs and also a lower estimated total FLOP count, defined as the multiplication of the estimated FLOPs per cell, number of cells, and inverse of the maximal stable time step size, where the common factors shared by the two methods are omitted.
	
	\begin{figure}[htb!]
		\centering
		\includegraphics[width=0.45\linewidth]{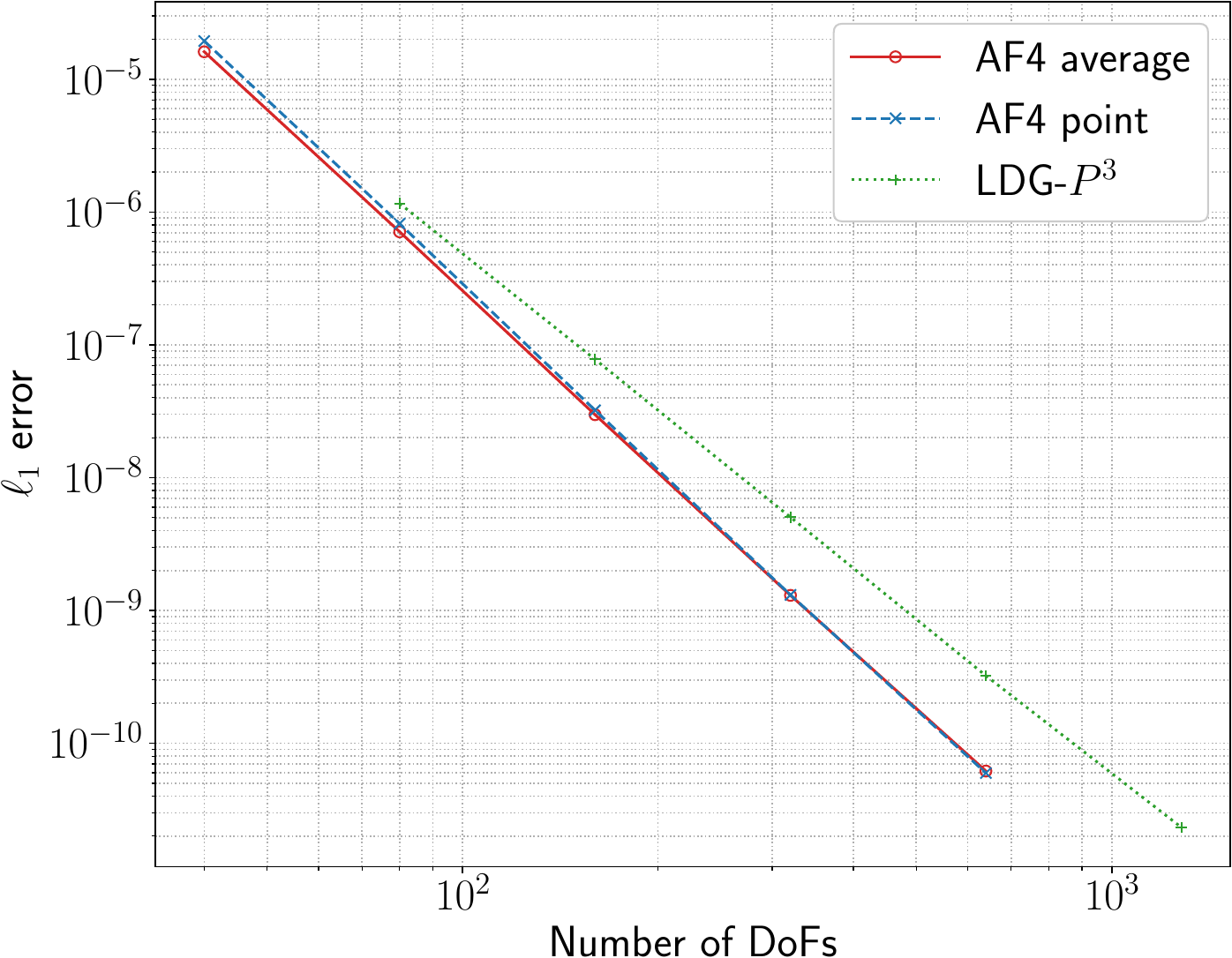}
		\includegraphics[width=0.45\linewidth]{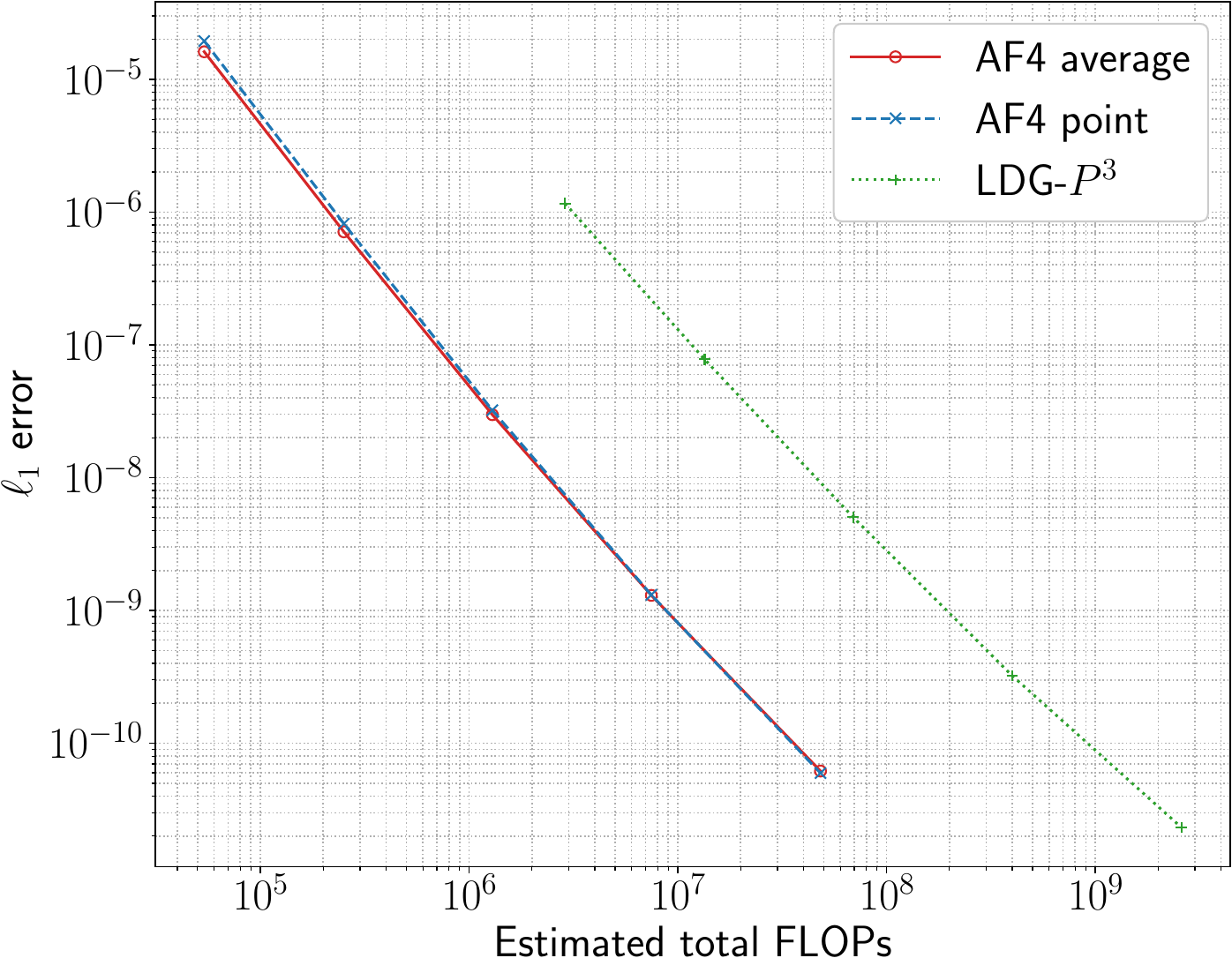}
		\caption{The comparison between the proposed AF and $P^3$ LDG methods for the 1D case.
			From left to right: the $\ell_1$ error versus the number of DoFs and estimated total FLOPs, respectively.}
		\label{fig:1d_comparison}
	\end{figure}
	
\end{example}

\begin{example}[Accuracy test with limiting]\label{ex:accuracy_limiting}
	To examine the effects of the MPP limiting, consider the advection equation $u_t + u_x = 0$ with the exact solution $u(x,t) = 1 + \sin^4(2\pi(x-t))$, thus the minimum and maximum are $u_{\min} = 1$ and $u_{\max} = 2$.
	This example is computed until $T=1$ on a periodic computational domain $[0,1]$.
	
	Table~\ref{tab:1d_accuracy_limiting} shows a comparison of the discrete minimum/maximum (considering all DoFs) without and with the proposed MPP limiting.
	One can clearly observe the violation of the discrete maximum principle without the MPP limiting, and the proposed limiting helps to maintain the discrete maximum principle up to roundoff precision.
	
	\begin{table}[htb!]
		\centering
		\begin{tabular}{ccccc}
			\toprule
			& $(u_h)_{\min}-u_{\min}$ & $(u_h)_{\max}-u_{\max}$ & $(u_h)_{\min}-u_{\min}$ & $(u_h)_{\max}-u_{\max}$ \\
			\cmidrule(lr){2-3} \cmidrule(lr){4-5}
			$N$ & \multicolumn{2}{c}{without limiting} & \multicolumn{2}{c}{with limiting} \\
			\midrule
			$20$ & $ \bm{-1.444\times 10^{-2}} $ & $-1.916\times 10^{-2}$ & $\bm{8.253\times 10^{-7}}$ & $-1.949\times 10^{-2}$ \\
			$40$ & $\bm{-1.334\times 10^{-3}}$ & $-1.374\times 10^{-3}$ & $\bm{2.740\times 10^{-6}}$ & $-1.437\times 10^{-3}$ \\
			$80$ & $\bm{-6.348\times 10^{-5}}$ & $-6.602\times 10^{-5}$ & $\bm{1.605\times 10^{-6}}$ & $-9.223\times 10^{-5}$ \\
			$160$ & $\bm{-2.129\times 10^{-6}}$ & $-2.581\times 10^{-6}$ & $\bm{1.110\times 10^{-7}}$ & $-1.138\times 10^{-5}$ \\
			$320$ & $\bm{-7.845\times 10^{-8}}$ & $-9.890\times 10^{-8}$ & $\bm{2.680\times 10^{-9}}$ & $-1.830\times 10^{-6}$ \\
			\bottomrule
		\end{tabular}
		\caption{Comparison of the discrete minimum/maximum without and with the MPP limiting in this paper.}
		\label{tab:1d_accuracy_limiting}
	\end{table}
	
	The $\ell_1$ errors and the corresponding convergence rates are presented in Figure~\ref{fig:accuracy_limiting}.
	It is seen from the upper left figure that the MPP limiting proposed in this paper retains the observed $4$th-order convergence rate in this test thanks to the monolithic formulation, whereas the one in \cite{Duan_2025_Active_SJSC} (upper right) reduces the convergence order.
	For the lower panel, the 1D convection--diffusion case in Example~\ref{ex:accuracy} with $c = 1$, $\kappa_0=0.01$ is solved with and without the MPP limiting.
	One observes that the proposed limiting does not affect the convergence rate.
	
	\begin{figure}[htb!]
		\centering
		\includegraphics[width=0.45\linewidth]{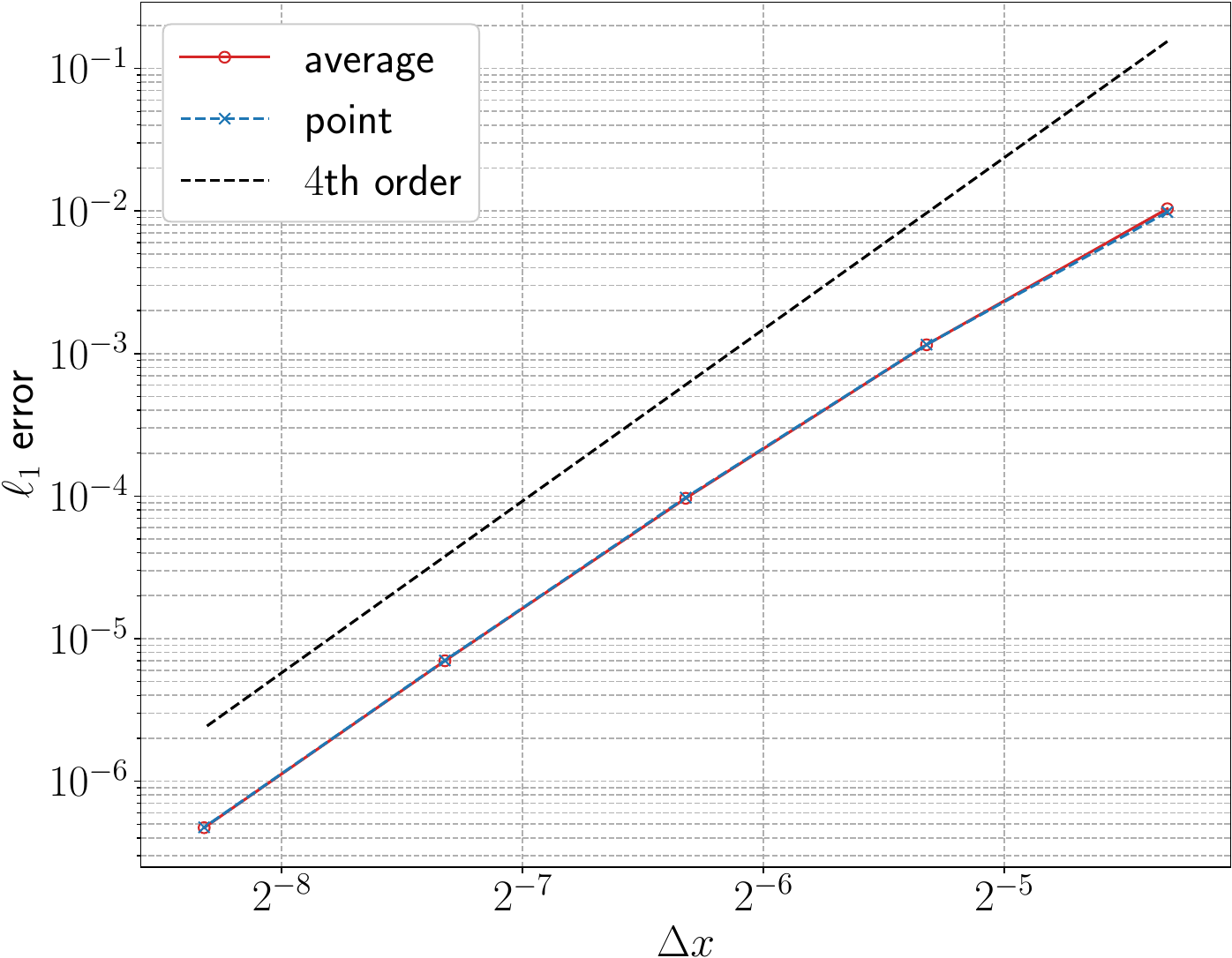}
		~
		\includegraphics[width=0.45\linewidth]{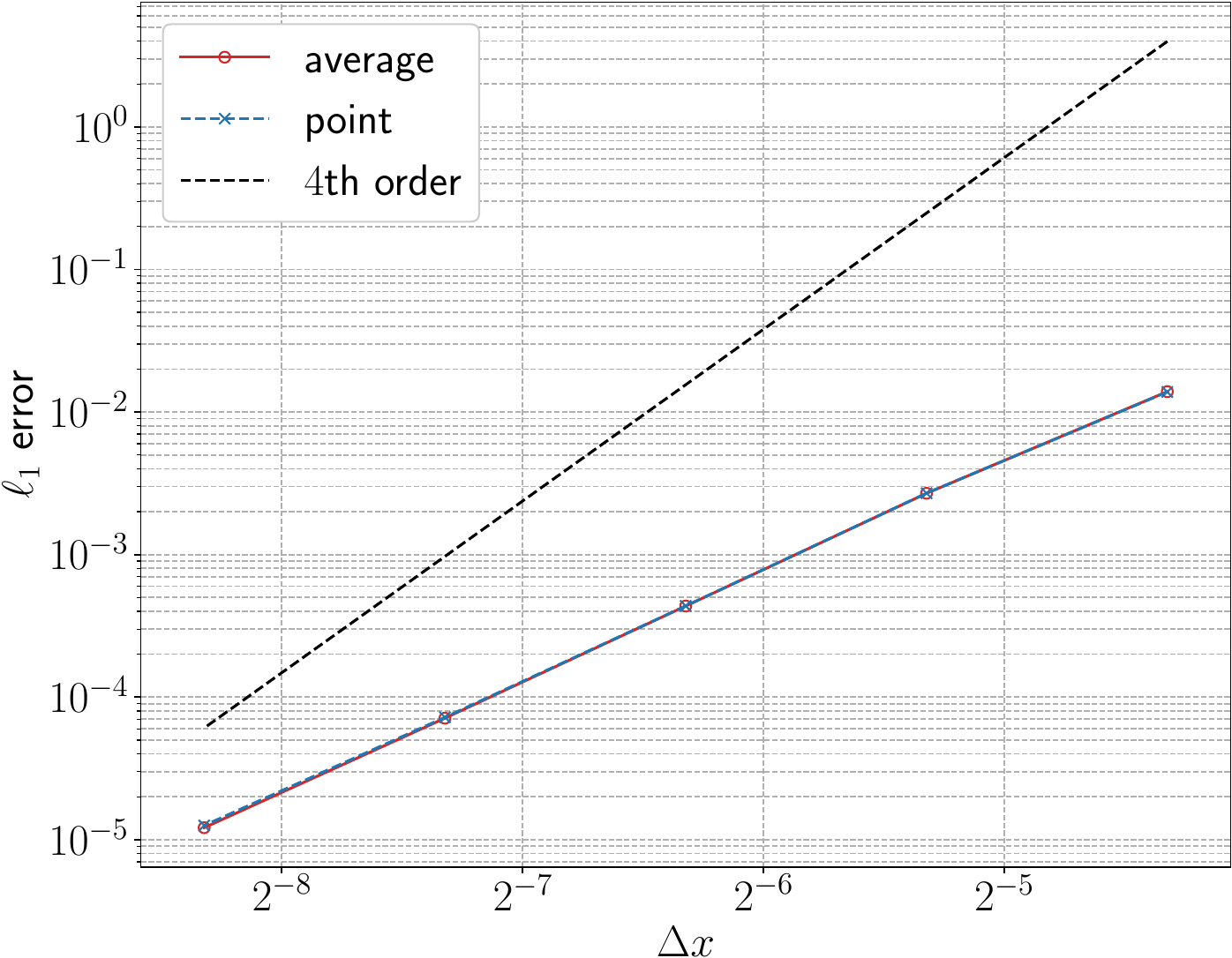}
		
		\vspace{2pt}
		\includegraphics[width=0.45\linewidth]{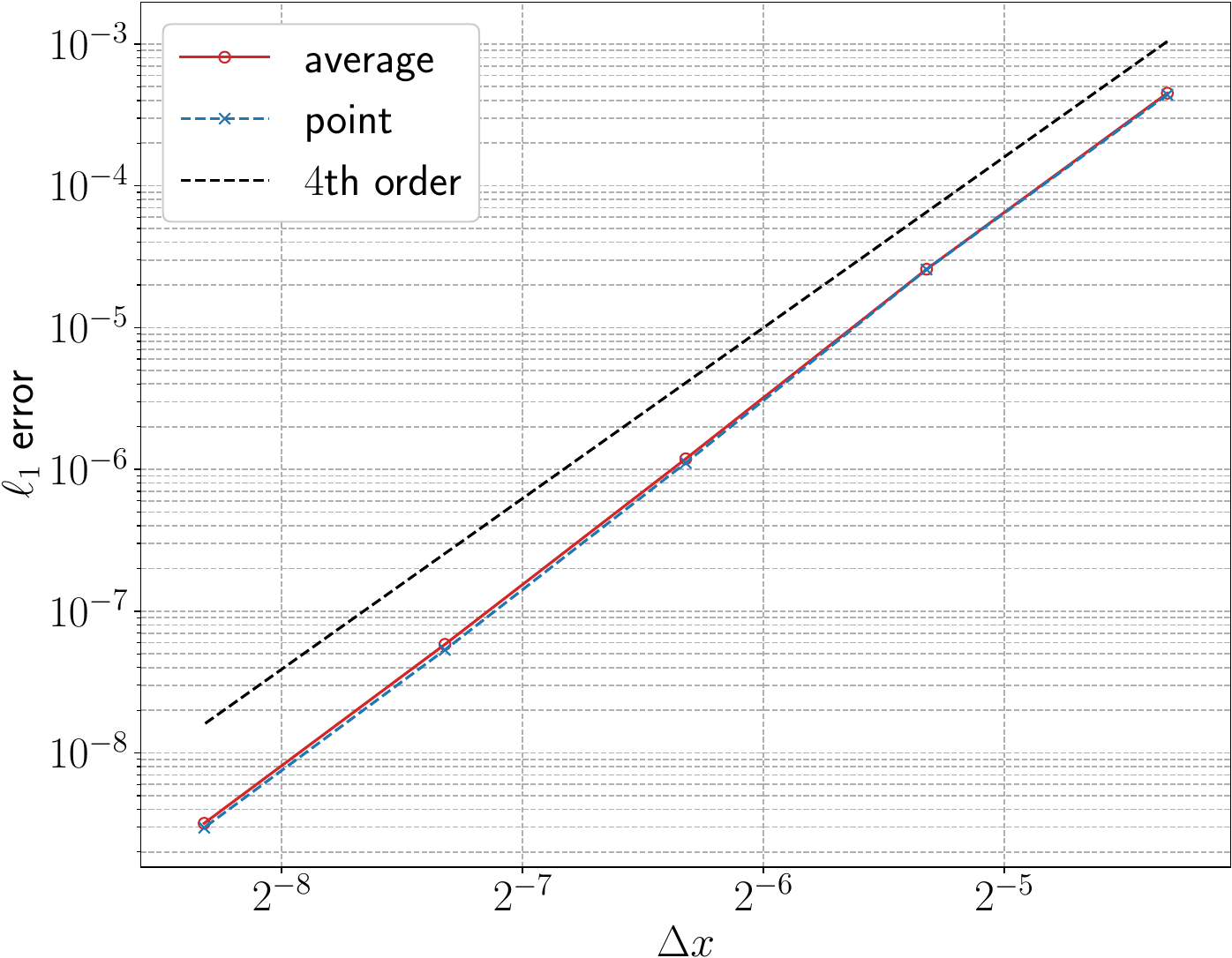}
		~
		\includegraphics[width=0.45\linewidth]{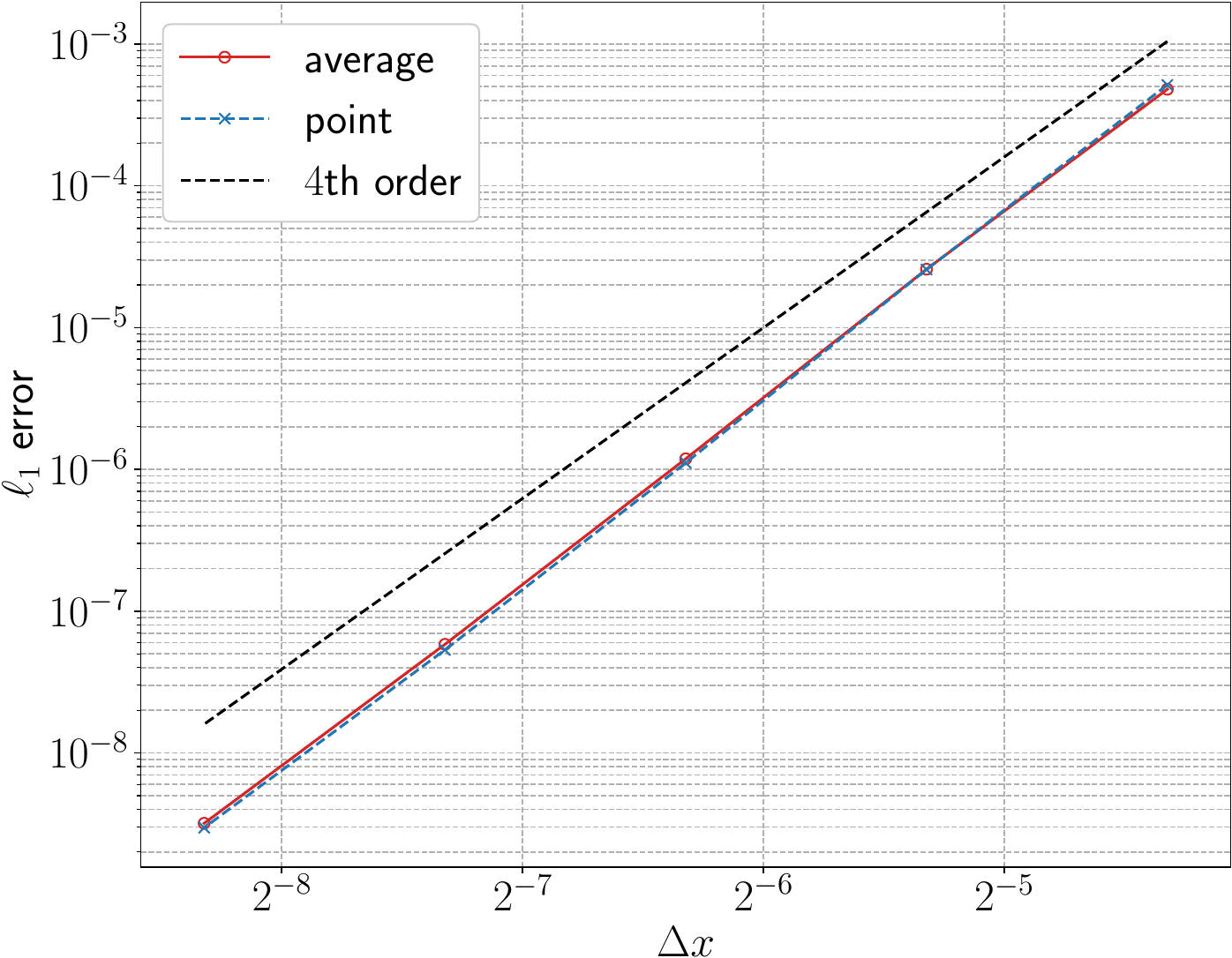}
		\caption{Example~\ref{ex:accuracy_limiting}.
			The $\ell_1$ errors and convergence rates of the proposed AF method.
			Upper left: $c = 1$, $\kappa_0=0$ with the MPP limiting in this paper,
			upper right: $c = 1$, $\kappa_0=0$ with the MPP limiting in \cite{Duan_2025_Active_SJSC},
			lower left: $c = 1$, $\kappa_0=0.01$ with the MPP limiting in this paper,
			lower right: $c = 1$, $\kappa_0=0.01$ without limiting.
		}
		\label{fig:accuracy_limiting}
	\end{figure}
\end{example}

In the following two tests, the proposed AF method is applied to the porous medium equation (PME),
\begin{equation}\label{eq:pme}
	u_t = \Delta (u^m),
\end{equation}
where $m \geqslant 2$ is an integer.
This equation can be used to model the flow of an isentropic gas through a porous medium, where $u \geqslant 0$ is the density of the gas and the gas pressure is proportional to $u^{m-1}$.
At $u = 0$, the equation degenerates, leading to interesting behavior, such as finite propagation speed and sharp fronts.
The diffusion coefficient is $m u^{m-1}$ in each direction.
Therefore, the non-negativity of $u$ is important both for the physical admissibility and the well-posedness of the PME.
The existing high-order methods for the PME include \cite{Zhang_2009_Numerical_JSC,Liu_2011_High_SJSC,Jiang_2021_High_JSC,Xu_2025_High_JCP,Yu_2026_Lax_MC}.

\begin{example}[Barenblatt solution]\label{ex:barenblatt}
	The Barenblatt solution \cite{Barenblatt_1952_some_PMiM} is a symmetric self-similar solution with compact support and is given by
	\begin{equation*}
		B_{m}(\bx,t) =  t^{-\alpha}\max\left[0,~
		1-\dfrac{\alpha(m-1)}{2dm}\dfrac{\norm{\bx}_2^2}{t^{2\alpha/d}}\right] ^{\frac{1}{m-1}},
	\end{equation*}
	where $\alpha=(m - 1 + 2/d)^{-1}$.
	In this example, the initial condition is taken as the solution at $t=1$, and the test is computed until $T = 2$.
	The computational domain is $[-6, 6]$, $[-8, 8]\times[-8,8]$, and $[-10, 10]\times[-10,10]\times[-10,10]$, for $d=1,2,3$, respectively.
	
	The numerical solutions with the MPP limiting on uniform meshes with $N=80$ cells in each direction are shown in Figures~\ref{fig:1d_barenblatt}, \ref{fig:2d_barenblatt}, and \ref{fig:3d_barenblatt} for $d=1,2,3$, respectively,
	and the corresponding cutlines are presented in Figures~\ref{fig:2d_barenblatt_cutline} and \ref{fig:3d_barenblatt_cutline}.
	The solution tends to move outwards, and the wave front is sharper with larger $m$.
	The proposed AF method can accurately capture the evolution and the sharp fronts without oscillations.
	
	\begin{figure}[htb!]
		\centering
		\includegraphics[width=\linewidth]{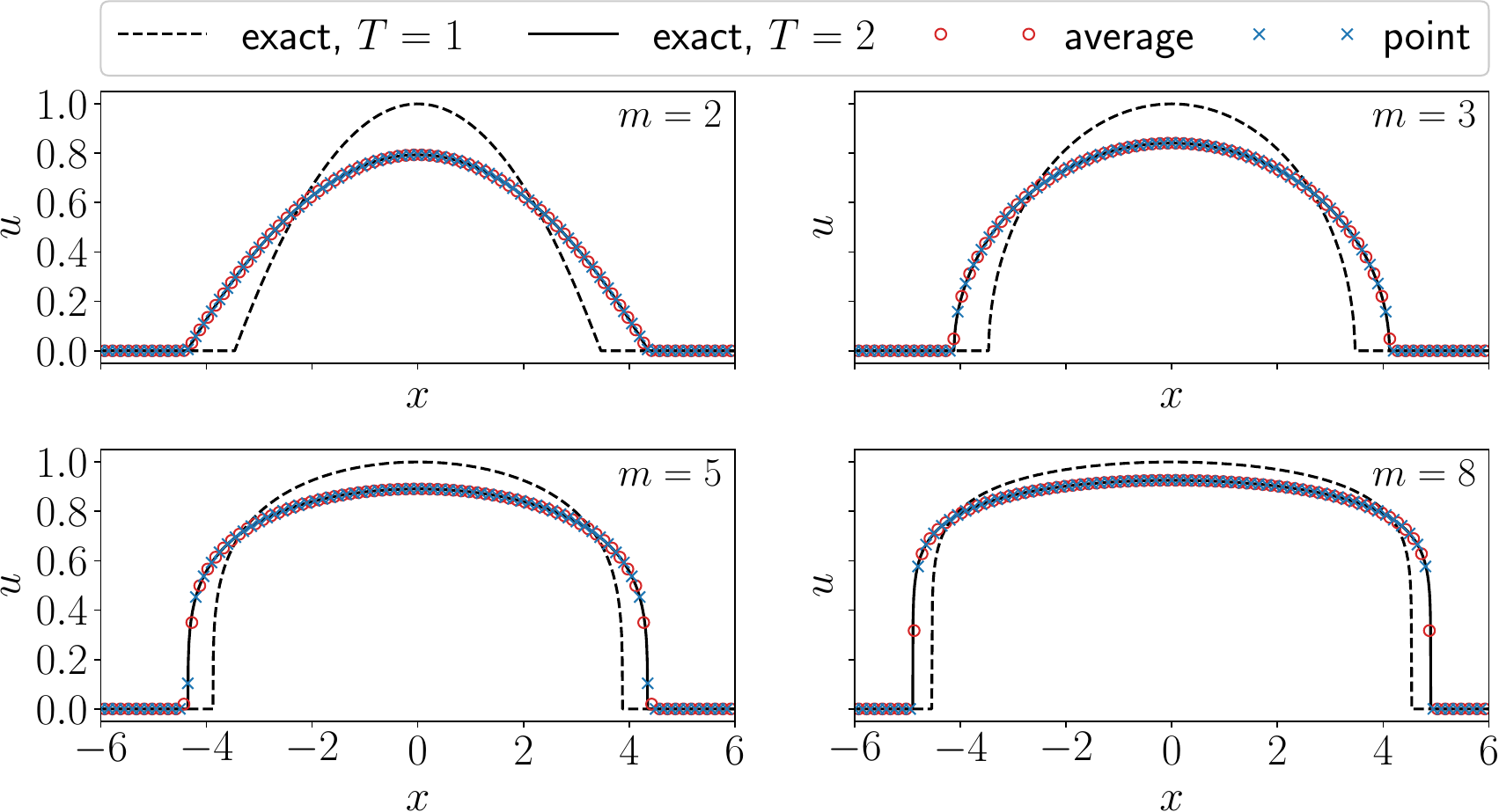}
		\caption{1D case of Example~\ref{ex:barenblatt}. The exact solutions, and numerical solutions obtained by the proposed AF method with the MPP limiting.}
		\label{fig:1d_barenblatt}
	\end{figure}
	
	\begin{figure}[htb!]
		\centering
		\includegraphics[width=0.45\linewidth]{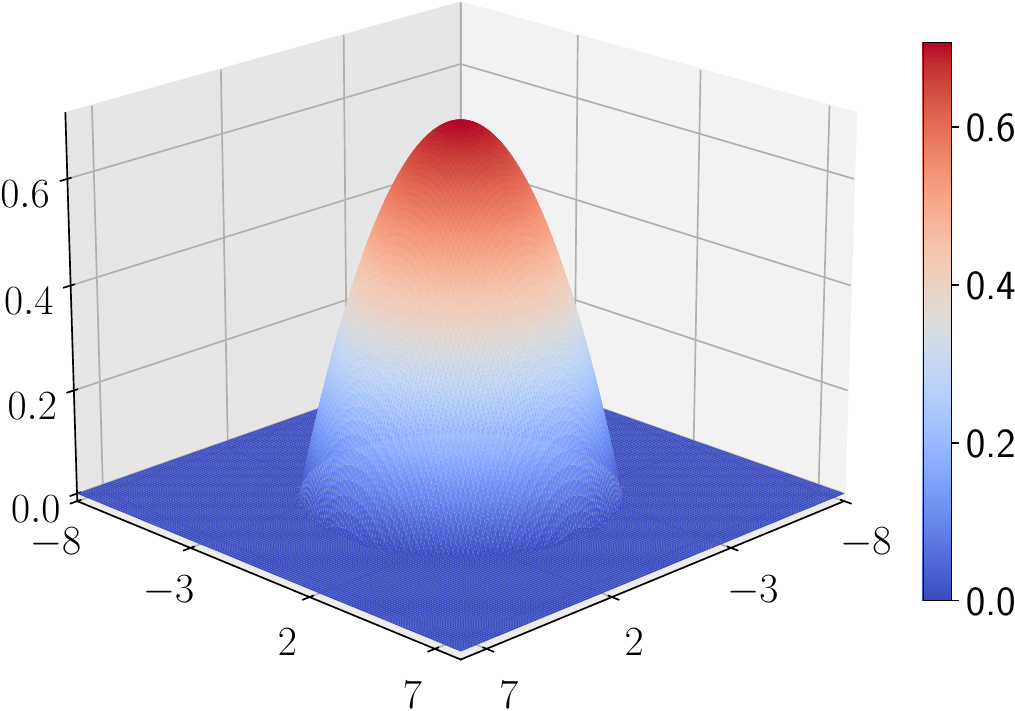}
		~
		\includegraphics[width=0.45\linewidth]{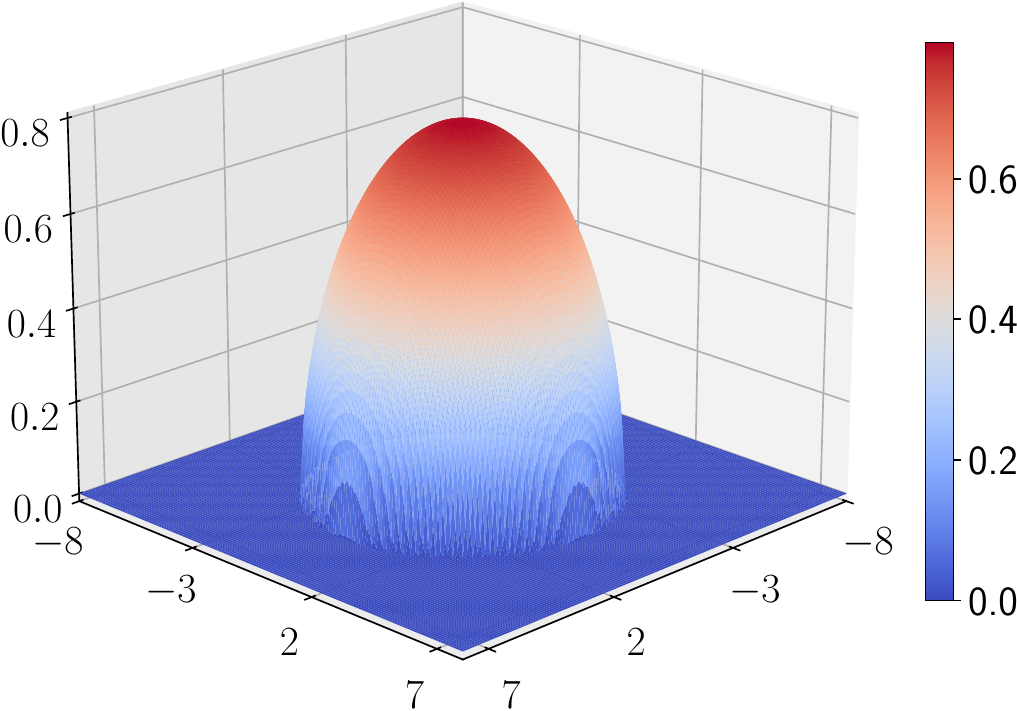}
		
		\includegraphics[width=0.45\linewidth]{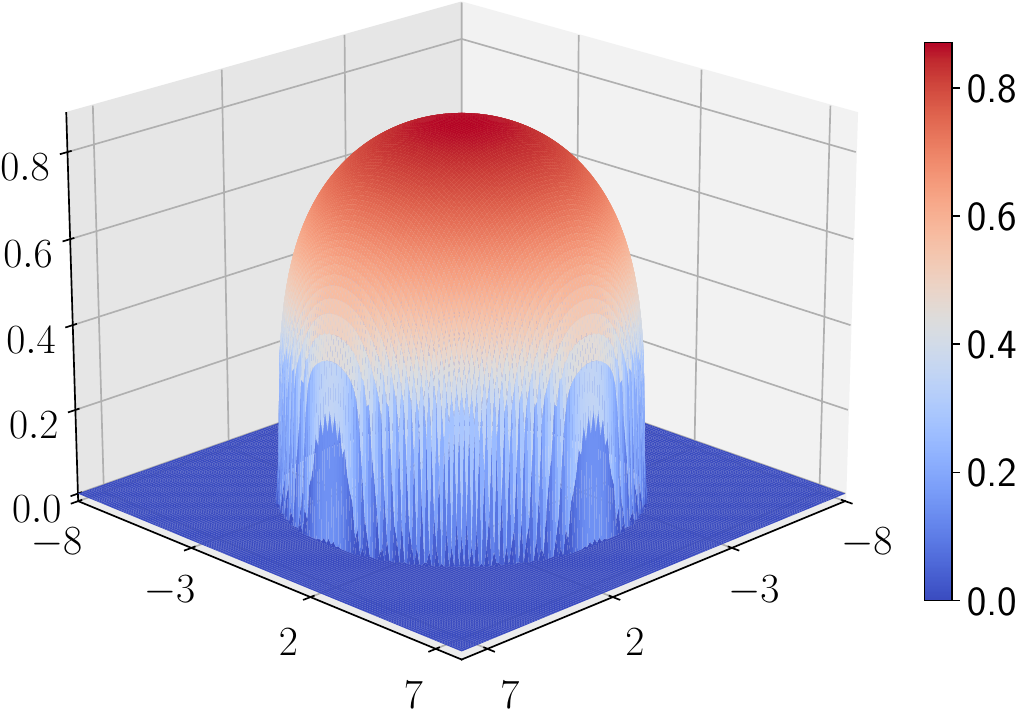}
		~
		\includegraphics[width=0.45\linewidth]{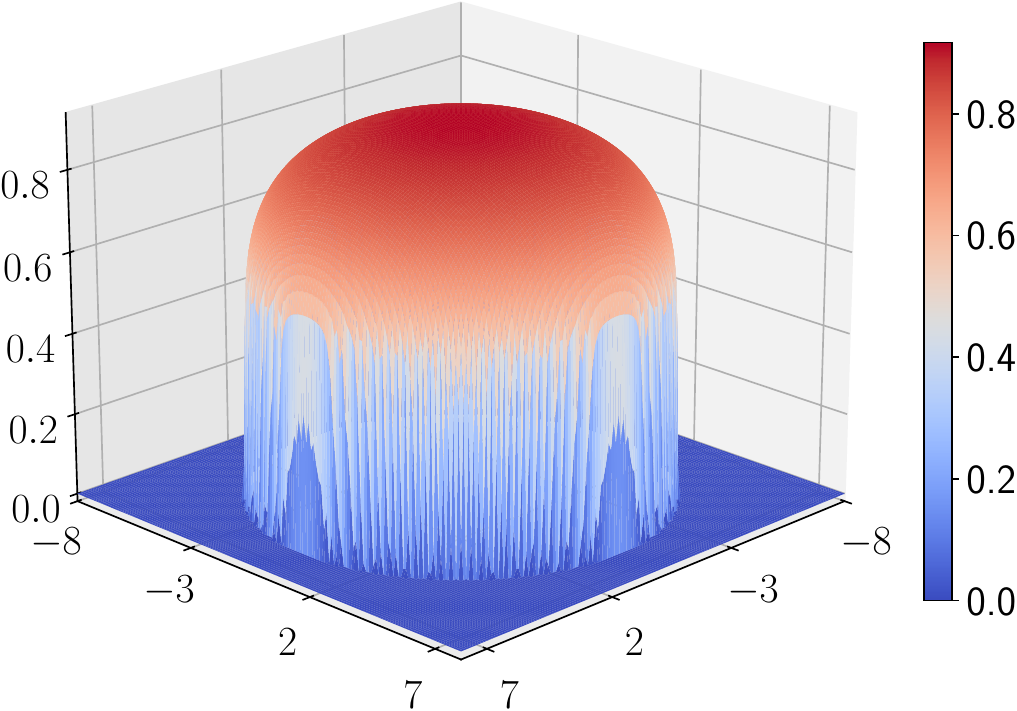}
		\caption{2D case of Example~\ref{ex:barenblatt}.
			The numerical solutions with $m=2,3,5,8$ obtained by the proposed AF method with the MPP limiting.}
		\label{fig:2d_barenblatt}
	\end{figure}
	
	\begin{figure}[htb!]
		\centering
		\includegraphics[width=\linewidth]{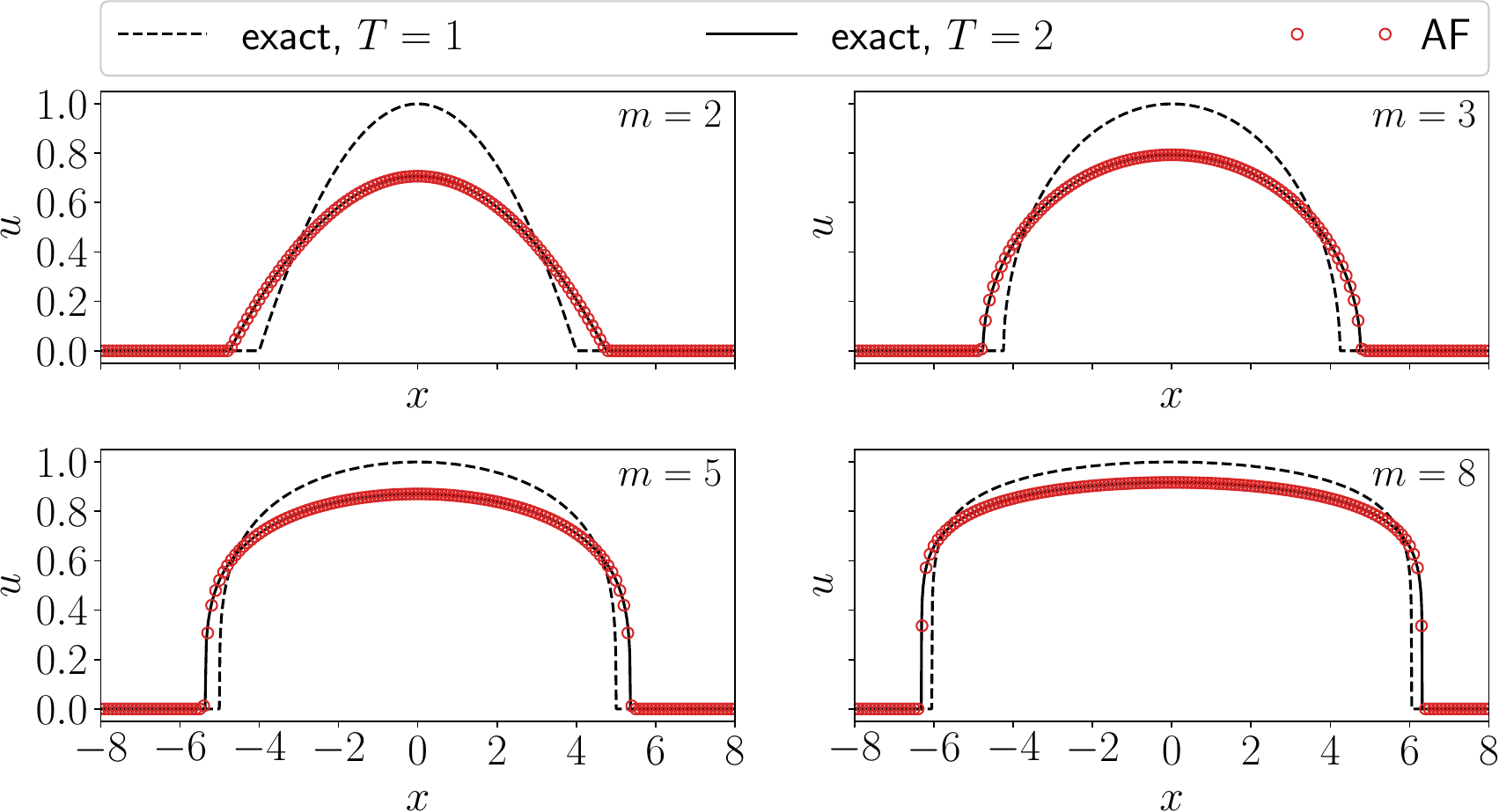}
		\caption{2D case of Example~\ref{ex:barenblatt}. The cutline along $y=0$ of the exact and numerical solutions with the MPP limiting.}
		\label{fig:2d_barenblatt_cutline}
	\end{figure}
	
	\begin{figure}[htb!]
		\centering
		\includegraphics[width=0.45\linewidth]{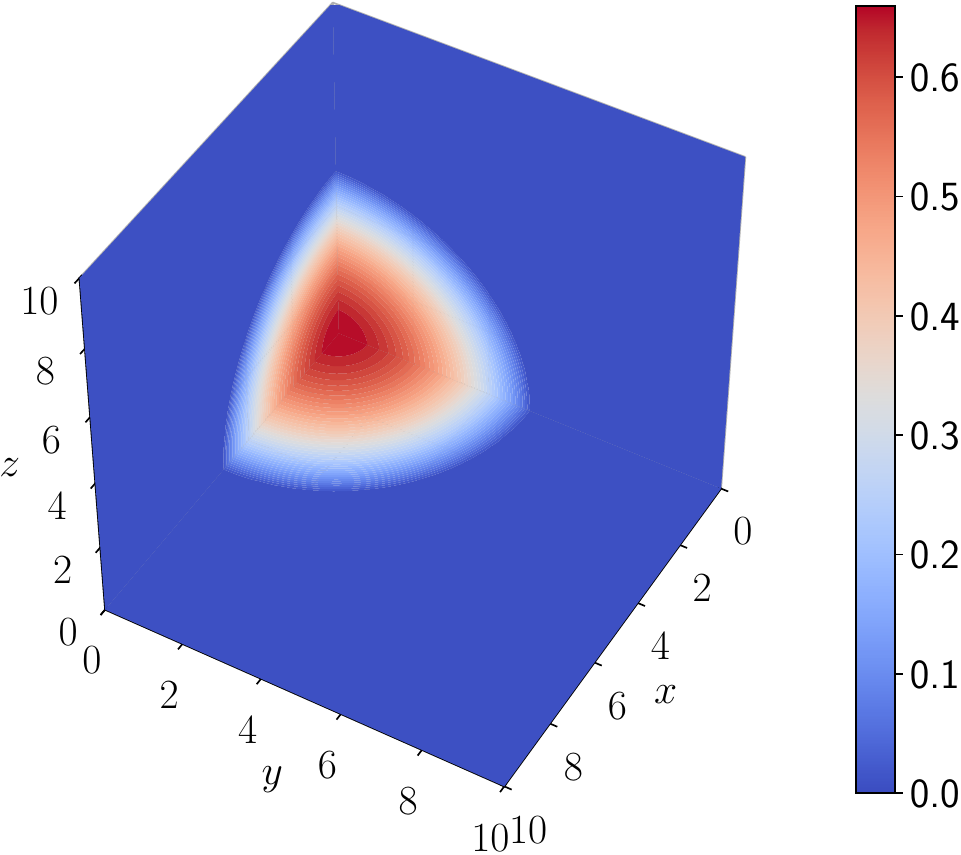}
		~
		\includegraphics[width=0.45\linewidth]{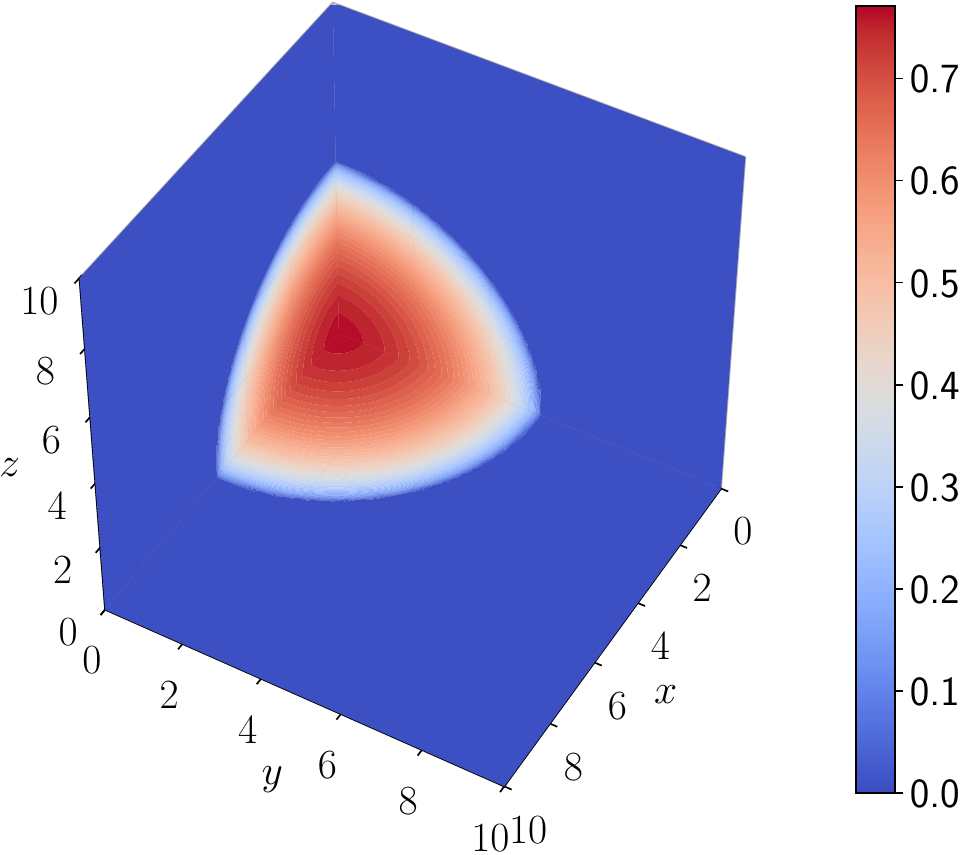}
		
		\includegraphics[width=0.45\linewidth]{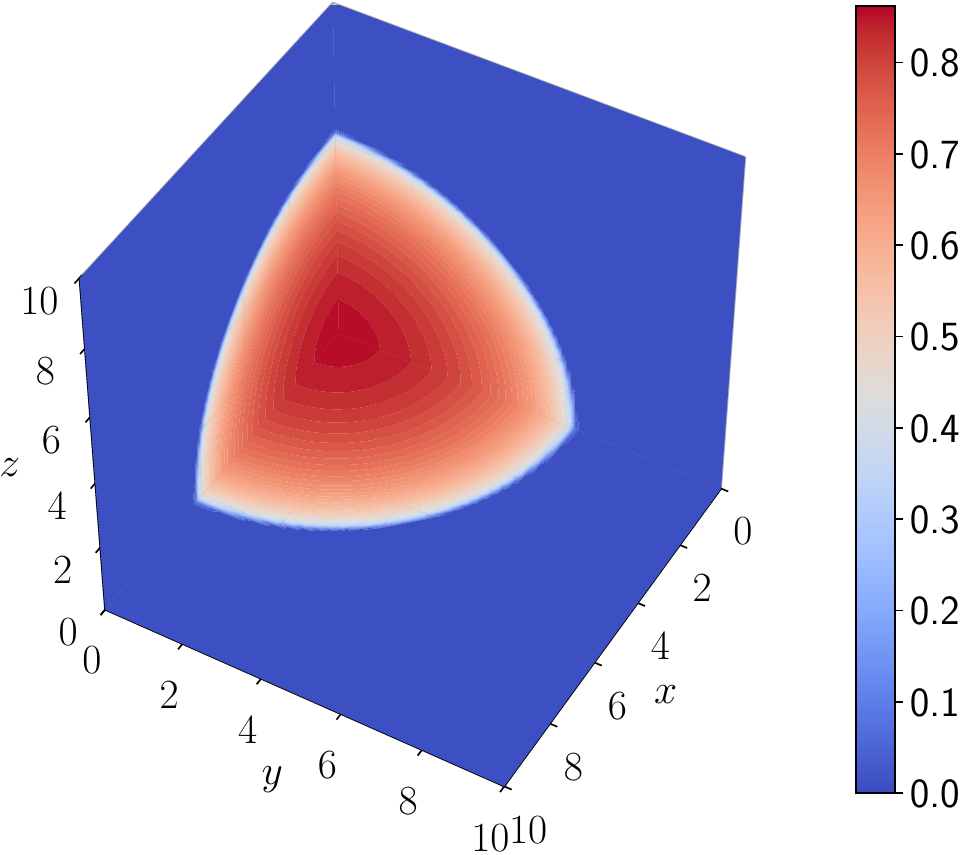}
		~
		\includegraphics[width=0.45\linewidth]{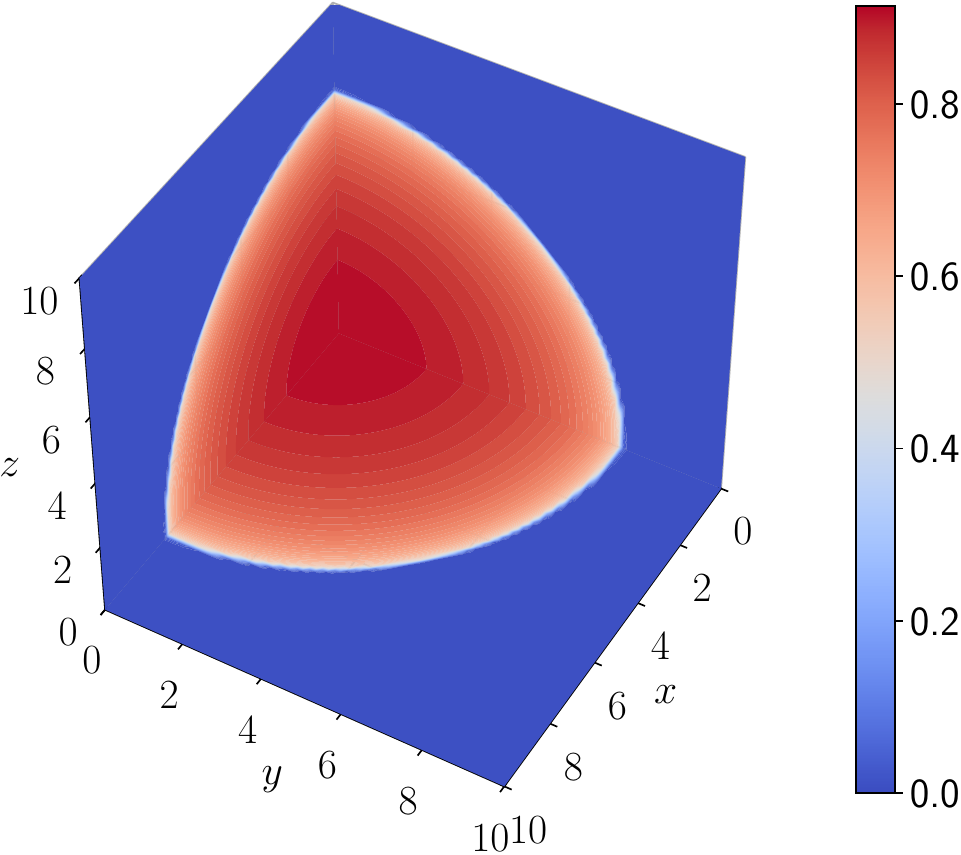}
		\caption{3D case of Example~\ref{ex:barenblatt}.
			The numerical solutions with $m=2,3,5,8$ obtained by the proposed AF method with the MPP limiting.}
		\label{fig:3d_barenblatt}
	\end{figure}
	
	\begin{figure}[htb!]
		\centering
		\includegraphics[width=\linewidth]{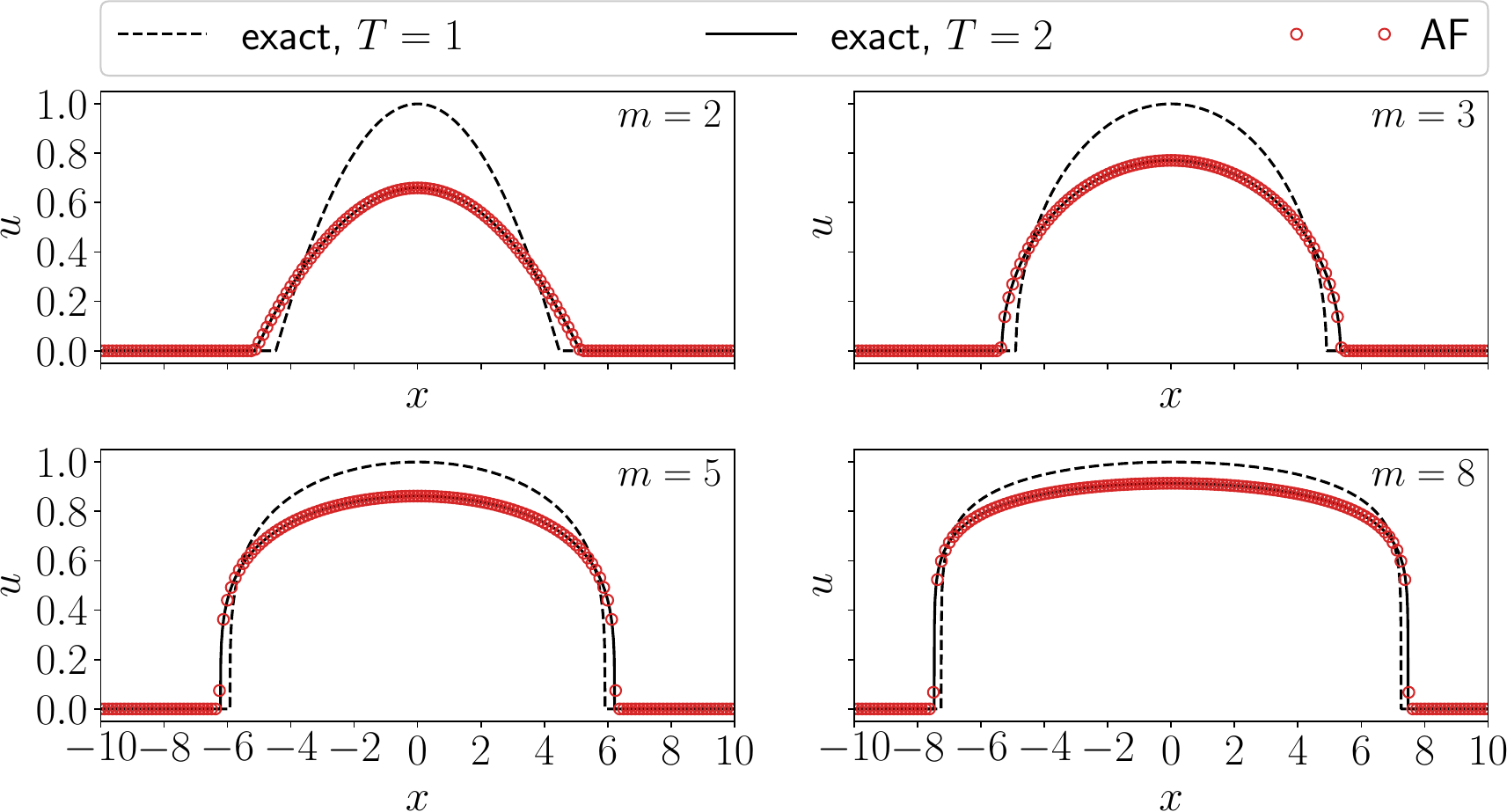}
		\caption{3D case of Example~\ref{ex:barenblatt}. The cutline along $y=z=0$ of the exact and numerical solutions with the MPP limiting.}
		\label{fig:3d_barenblatt_cutline}
	\end{figure}
	
	To examine the efficacy of the MPP limiting, the enlarged views of the 1D numerical solutions are shown in Figure~\ref{fig:barenblatt_enlarge}.
	The numerical solutions maintain non-negativity with the help of the limiting; otherwise, the undershoots would violate non-negativity.
	For the 2D and 3D cases, the numerical solutions without the limiting generate negative solutions, see the minimum in Table~\ref{tab:barenblatt_minimum},
	thus the MPP limiting is essential for preserving the discrete maximum principle.
	
	\begin{figure}[htb!]
		\centering
		\includegraphics[width=0.475\linewidth]{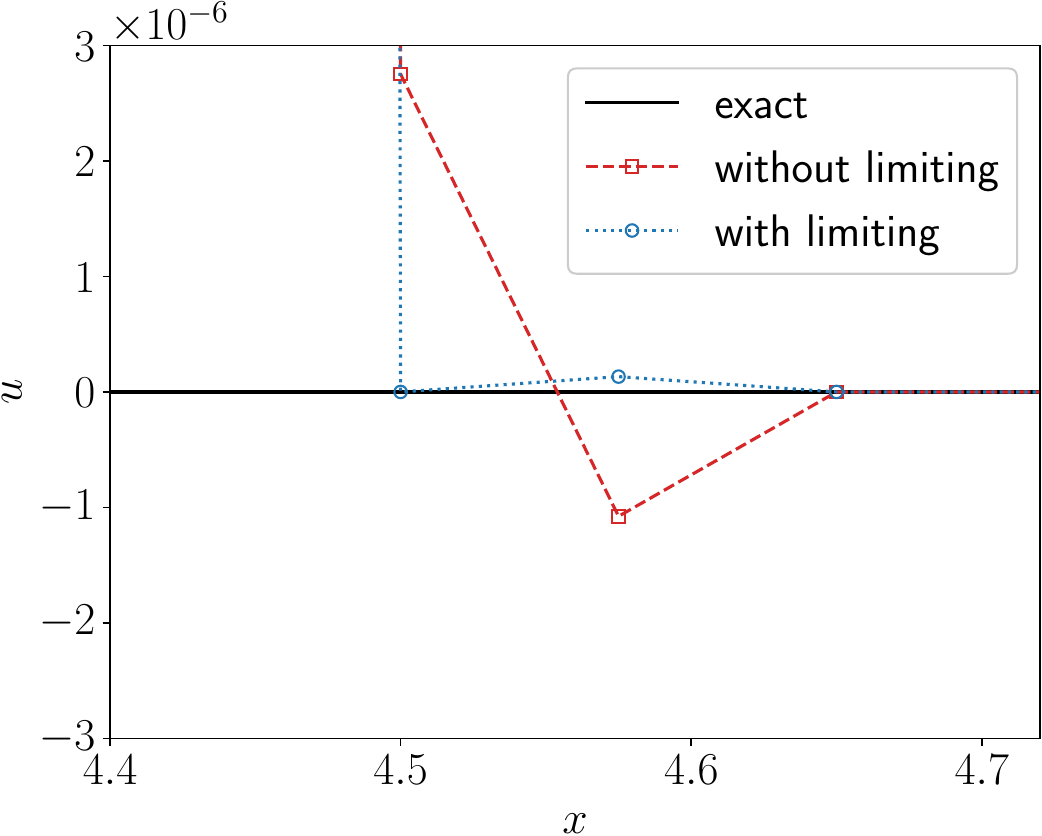}\hfill
		\includegraphics[width=0.51\linewidth]{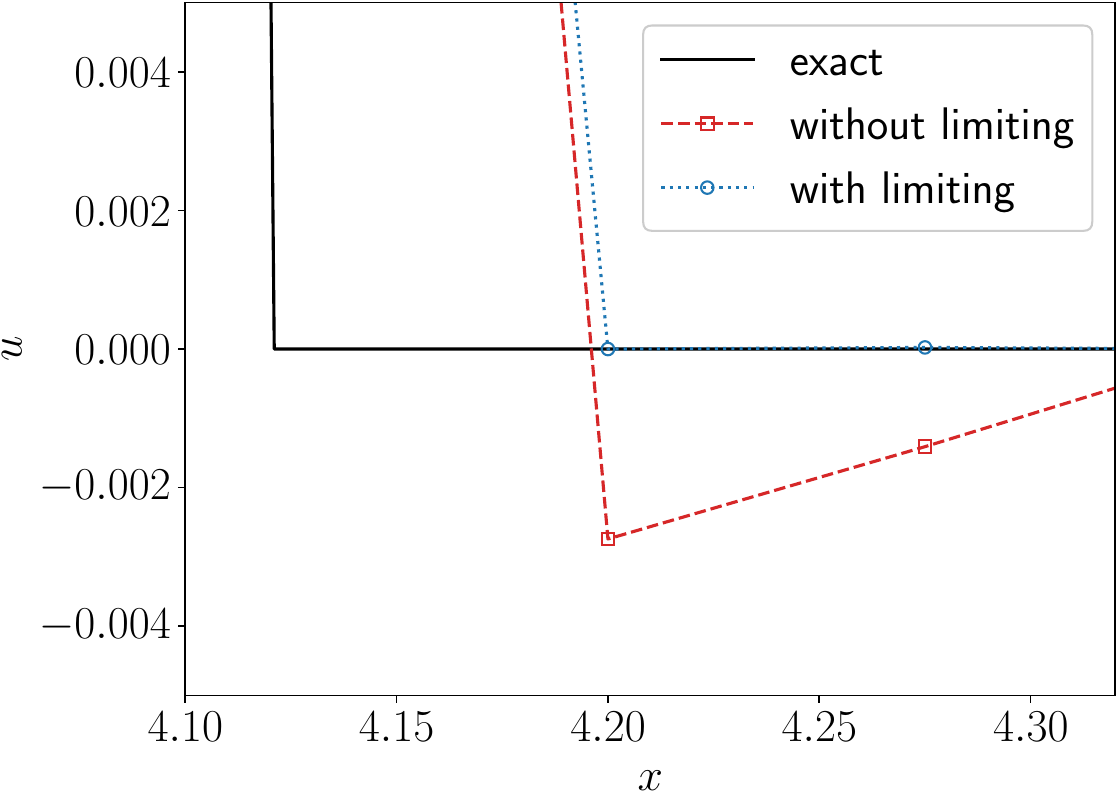}\hfill
		
		\includegraphics[width=0.47\linewidth]{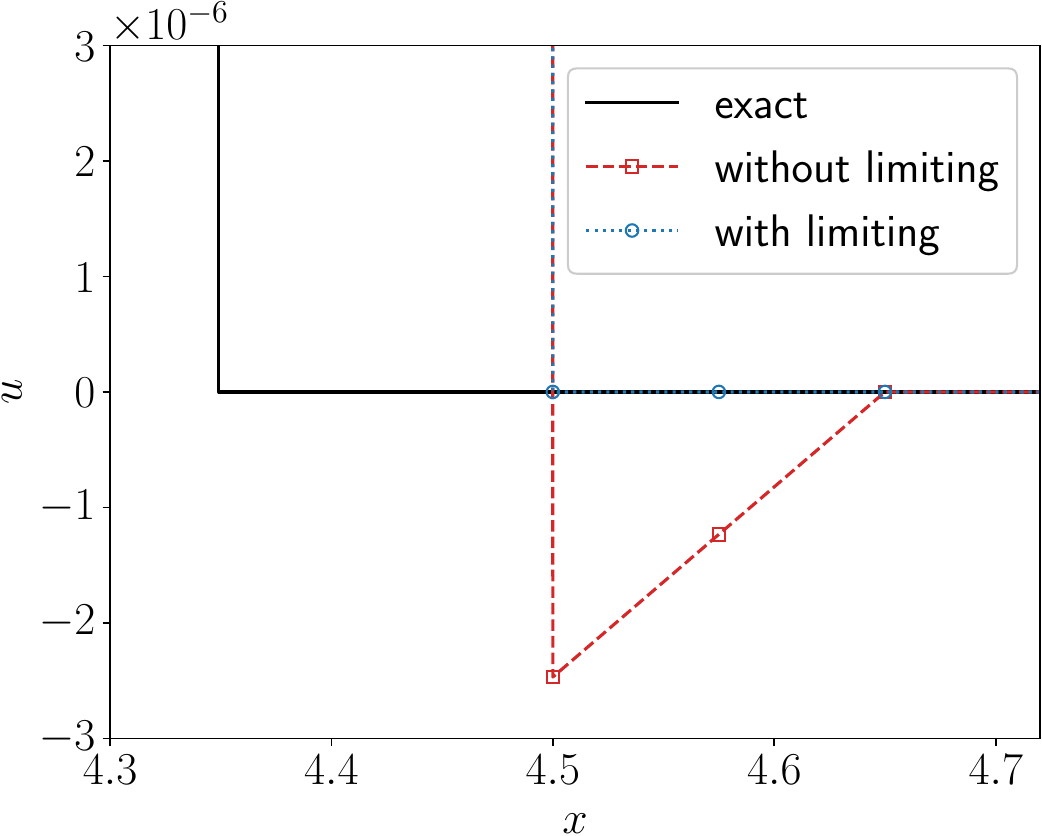}\hfill
		\includegraphics[width=0.495\linewidth]{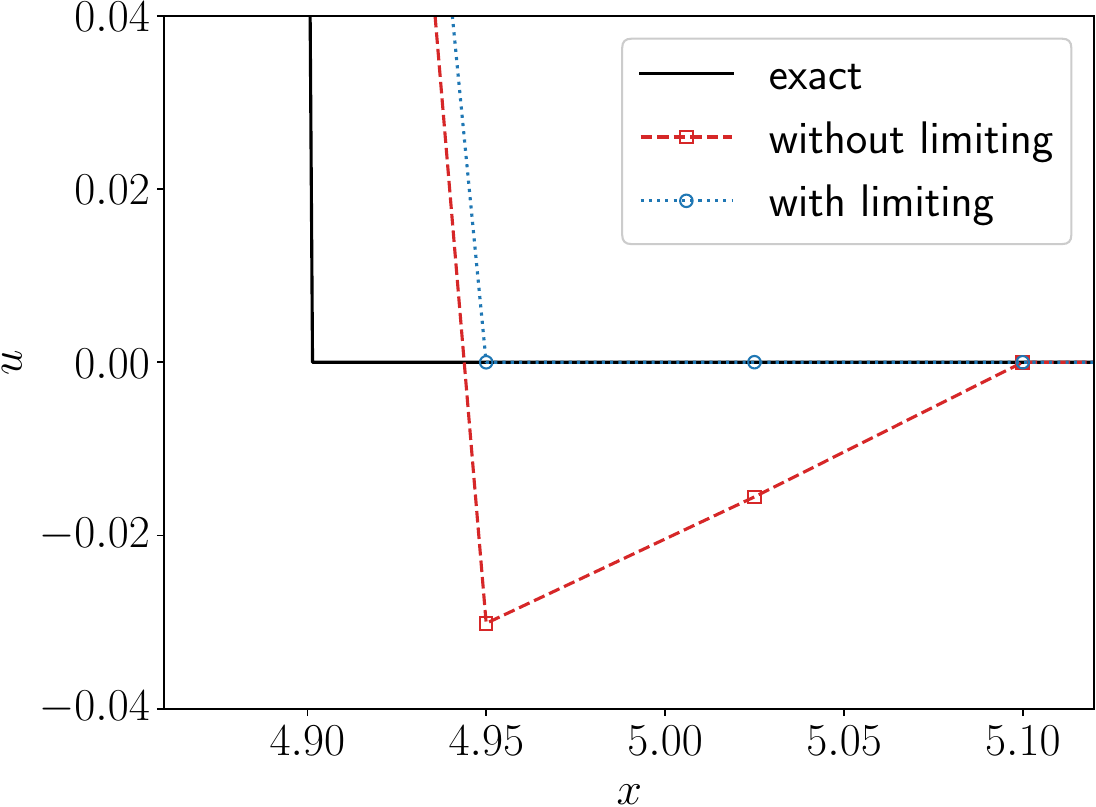}
		\caption{1D case of Example~\ref{ex:barenblatt}. The enlarged views of the numerical solutions without and with the MPP limiting.
			From left to right, top to bottom: $m=2, 3, 5, 8$.}
		\label{fig:barenblatt_enlarge}
	\end{figure}
	
	\begin{table}[htb!]
		\centering
		\begin{tabular}{ccccc}
			\toprule
			& \multicolumn{2}{c}{2D} & \multicolumn{2}{c}{3D} \\
			\cmidrule(lr){2-3} \cmidrule(lr){4-5}
			$m$ & without limiting & with limiting & without limiting & with limiting \\
			\midrule
			$2$ & $-6.048\times 10^{-4}$ & $0.000$ & $-6.645\times 10^{-4}$ & $0.000$ \\
			$3$ & $-4.736\times 10^{-3}$ & $0.000$ & $-4.997\times 10^{-3}$ & $0.000$ \\
			$5$ & $-1.621\times 10^{-2}$ & $0.000$ & $-1.639\times 10^{-2}$ & $0.000$ \\
			$8$ & $-2.897\times 10^{-2}$ & $0.000$ & $-2.893\times 10^{-2}$ & $0.000$ \\
			\bottomrule
		\end{tabular}
		\caption{2D and 3D cases of Example~\ref{ex:barenblatt}. Comparison of the discrete minimum without and with the MPP limiting.}
		\label{tab:barenblatt_minimum}
	\end{table}

\end{example}

\begin{example}[Collision of two hills]\label{ex:two_hill}
	This is a 2D test taken from \cite{Liu_2011_High_SJSC} with $m=2$.
	The initial condition is
	\begin{equation*}
		u(\bx,0) = \begin{cases}
			\exp\left(-\frac{1}{6 - (x-2)^2 - (y+2)^2}\right), &\text{if}~(x-2)^2 + (y+2)^2 < 6, \\
			\exp\left(-\frac{1}{6 - (x+2)^2 - (y-2)^2}\right), &\text{if}~(x+2)^2 + (y-2)^2 < 6, \\
			0, &\text{otherwise}.
		\end{cases}
	\end{equation*}
	The computational domain is $[-10,10]\times[-10,10]$ with periodic boundary conditions.
	
	The evolution of the numerical solution with the MPP limiting is shown in Figure~\ref{fig:two_hill}, obtained with $100\times100$ cells.
	The initial two hills diffuse and merge, which is comparable to the behavior shown in \cite{Liu_2011_High_SJSC}.
	If the MPP limiting is not activated, negative solutions appear during the simulation.
	
	\begin{figure}[htb!]
		\centering
		\includegraphics[width=0.32\linewidth]{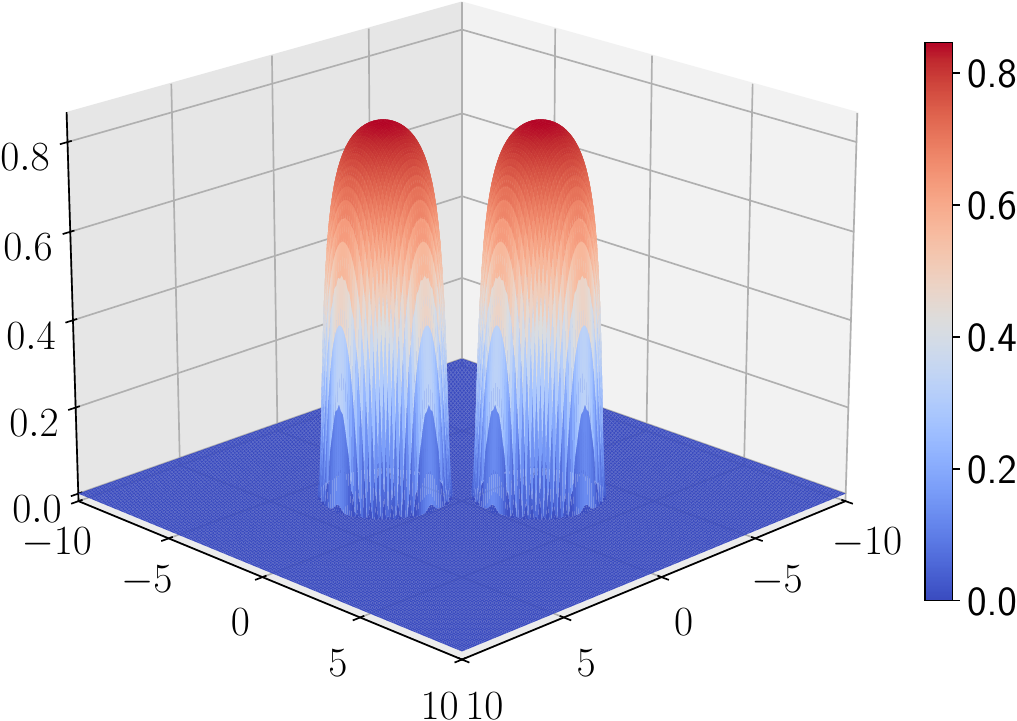}
		\includegraphics[width=0.32\linewidth]{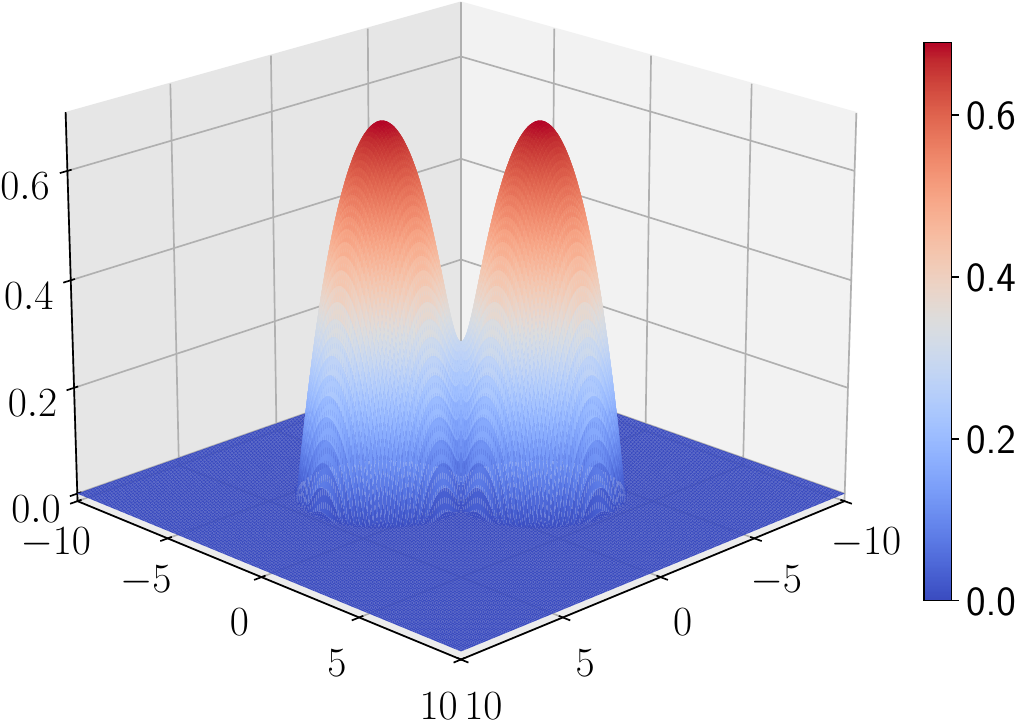}
		\includegraphics[width=0.32\linewidth]{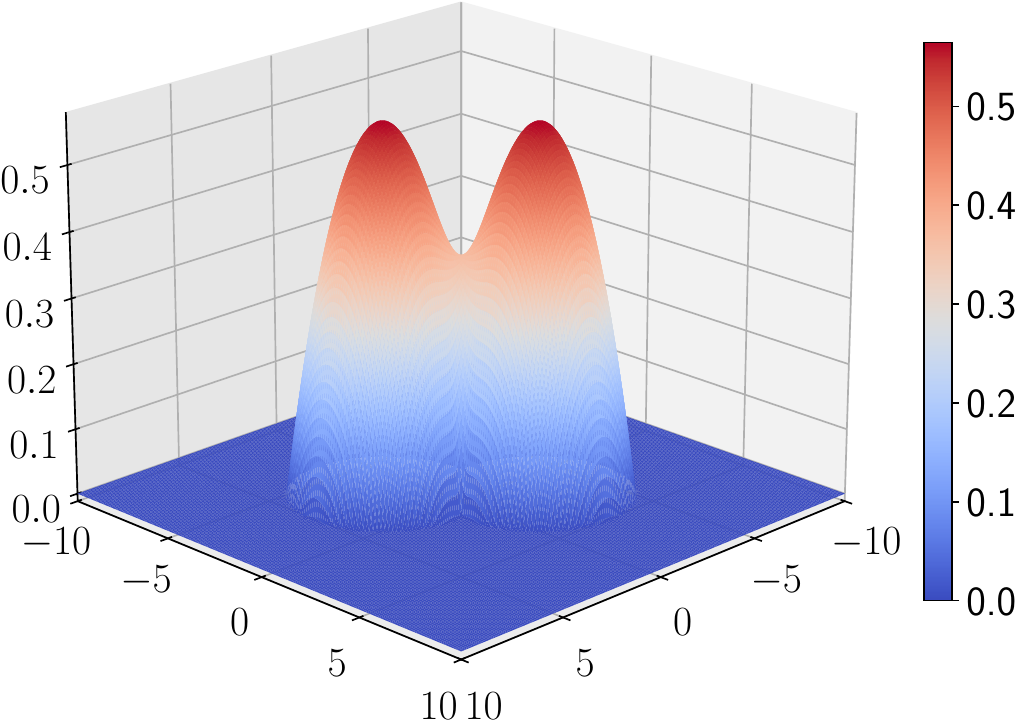}
		
		\includegraphics[width=0.32\linewidth]{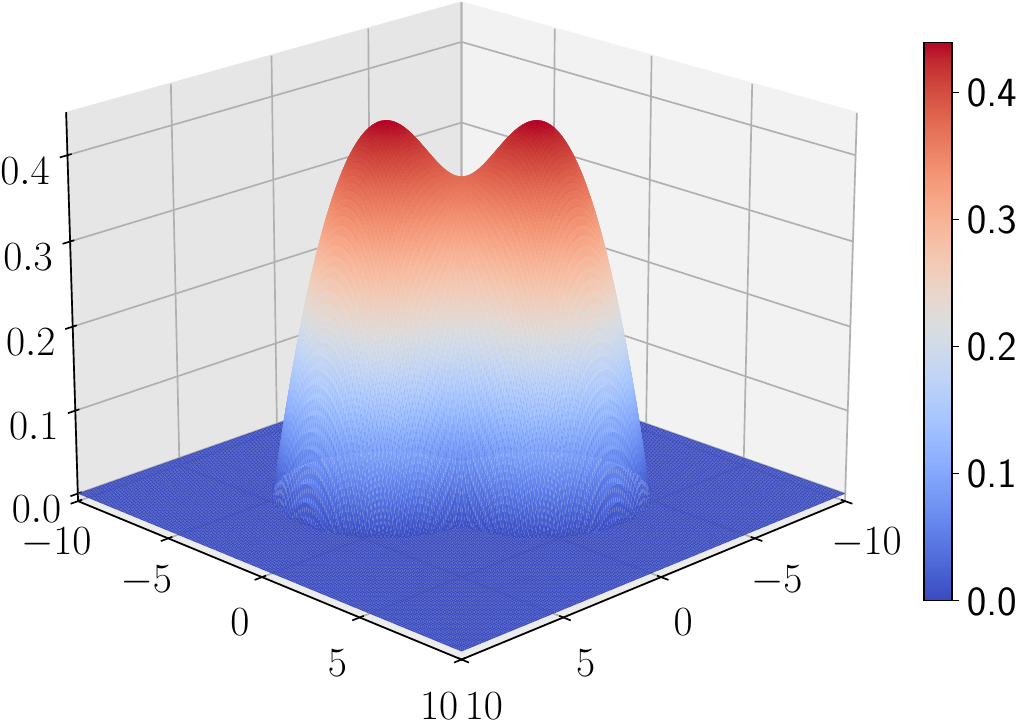}
		\includegraphics[width=0.32\linewidth]{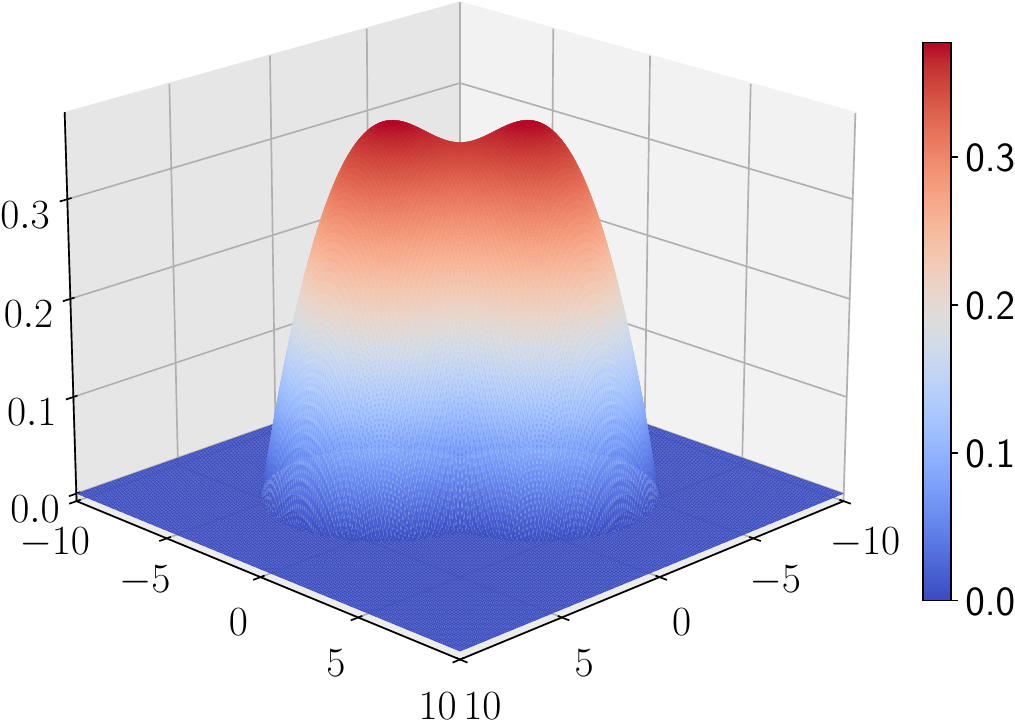}
		\includegraphics[width=0.32\linewidth]{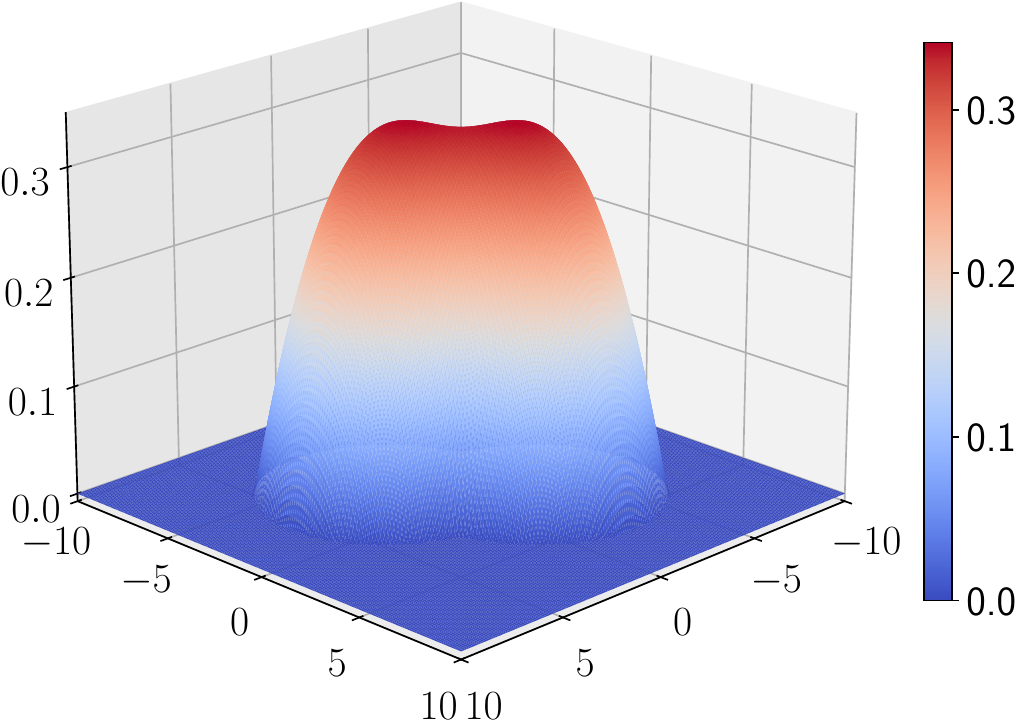}
		\caption{Example~\ref{ex:two_hill}.
			The numerical solutions obtained by the proposed AF method with the MPP limiting at $T=0$, $0.5$, $1$, $2$, $3$, $4$ (from left to right, top to bottom).}
		\label{fig:two_hill}
	\end{figure}	
\end{example}

\begin{example}[Evje and Karlsen's example]\label{ex:1d_degenerate_Evje_Karlsen}
	This example solves \eqref{eq:1d_cde_dc} with $f(u) = u^2/4$ and
	\begin{equation*}
		a(u) =
		\begin{cases}
			0, &\text{if}\ u \leqslant 0.5, \\
			5(2u-1)^2/4, &\text{if}\ 0.5 < u < 0.6, \\
			u - 11/20, &\text{otherwise}, \\
		\end{cases}
	\end{equation*}
	then the corresponding diffusion coefficient is 
	\begin{equation*}
		\kappa(u) =
		\begin{cases}
			0, &\text{if}\ u \leqslant 0.5, \\
			5(2u-1), &\text{if}\ 0.5 < u < 0.6, \\
			1, &\text{otherwise}. \\
		\end{cases}
	\end{equation*}
	The initial condition is 
	\begin{equation*}
		u(x,0) =
		\begin{cases}
			1, &\text{if}\ x\in[-1, -0.4), \\
			-2.5x, &\text{if}\ x\in[-0.4, 0), \\
			0, &\text{otherwise},
		\end{cases}
	\end{equation*}
	and the computational domain is $[-1, 1]$ with Dirichlet boundary conditions $u(-1,t) = 1$, $u(1,t) = 0$.
	
	Figure~\ref{fig:1d_degenerate_Evje_Karlsen} shows the numerical solutions obtained by the proposed AF method with the MPP limiting and $200$ cells at $T=1$.
	The proposed doubly conservative AF method accurately captures the wave front.
	In contrast, the method without double conservation fails to approximate the reference solution, demonstrating the numerical importance of double conservation in this degenerate problem.
	
	\begin{figure}[htb!]
		\centering
		\includegraphics[width=0.45\textwidth]{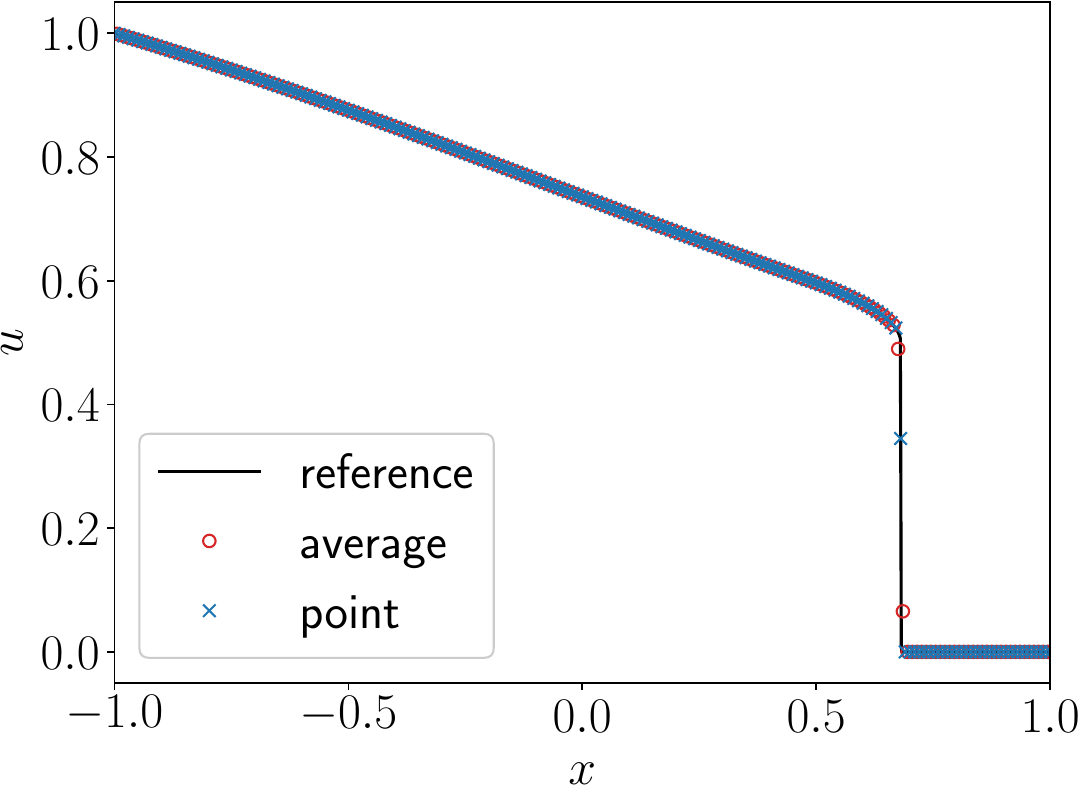}
		\quad
		\includegraphics[width=0.45\textwidth]{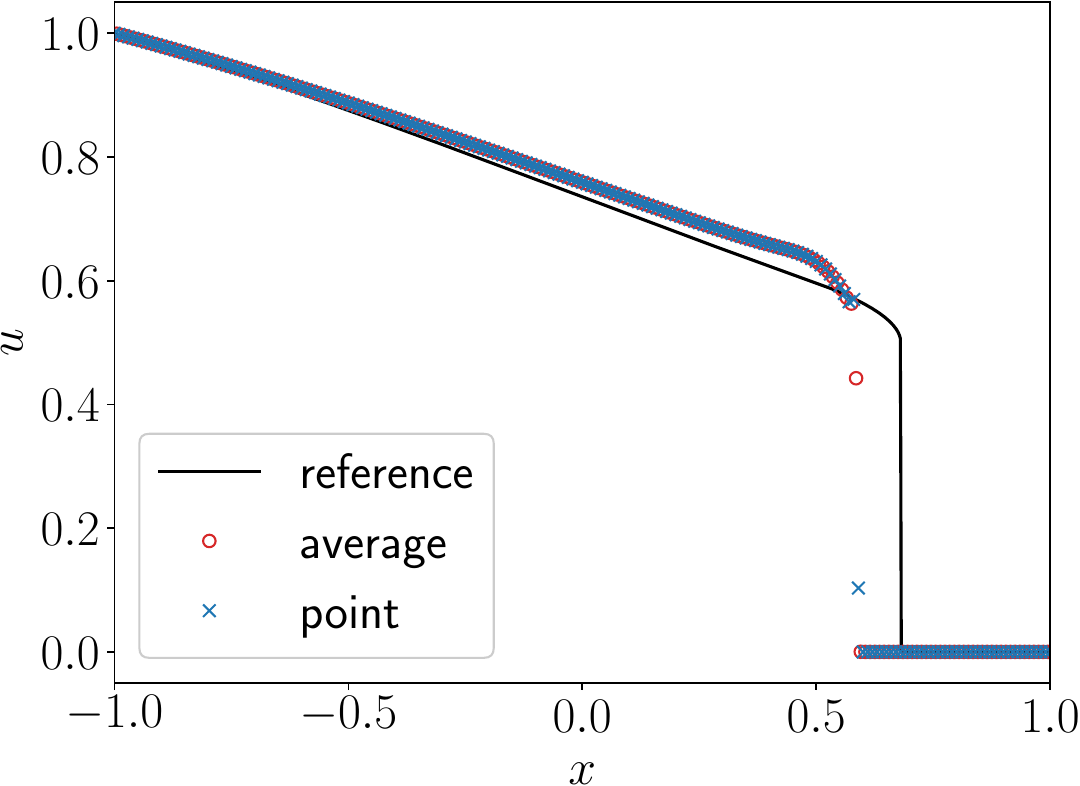}
		\caption{Evje and Karlsen's example~\ref{ex:1d_degenerate_Evje_Karlsen}.
			Left: doubly conservative AF method, right: the AF method without double conservation.
			The numerical solutions are obtained by the proposed AF methods with the MPP limiting and $200$ cells and the reference solution is obtained by a first-order scheme with $2000$ cells.}
		\label{fig:1d_degenerate_Evje_Karlsen}
	\end{figure}
	
\end{example}

\begin{example}[1D Buckley Leverett]\label{ex:1d_buckley_leverett}
	This example solves \eqref{eq:1d_cde_dc} with
	\begin{equation*}
		a(u) =
		\begin{cases}
			0, &\text{if}\ u < 0, \\
			\epsilon (2u^2 - \frac43 u^3), &\text{if}\ 0\leqslant u \leqslant 1, \\
			\frac23 \epsilon, &\text{otherwise}, \\
		\end{cases}
	\end{equation*}
	where $\epsilon=0.01$.
	The computational domain is $[0,1]$.
	Two cases of the flux function are considered.
	The first case adopts the flux function without the gravitational effect
	$f(u) = \frac{u^2}{u^2 + (1-u)^2}$,
	while the second with the gravitational effect uses $f(u) = \frac{u^2}{u^2 + (1-u)^2} \left( 1 - 5(1-u)^2 \right)$.
	For each case, an initial boundary value problem is considered with the following initial condition,
	
	\begin{equation*}
		u(x,0) =
		\begin{cases}
			1-3x, &\text{if}\ x\in[0, \frac13), \\
			0, &\text{otherwise},
		\end{cases}
	\end{equation*}
	with the Dirichlet boundary condition $u(0,t) = 1$, $u(1,t)=0$.
	A Riemann problem is also solved with the initial data
	\begin{equation*}
		u(x,0) =
		\begin{cases}
			0, &\text{if}\ x\in[0, 1-\frac{1}{\sqrt{2}}), \\
			1, &\text{otherwise}.
		\end{cases}
	\end{equation*}
	
	Figure~\ref{fig:1d_buckley_leverett} shows the numerical solutions obtained by the proposed AF method with the MPP limiting and $100$ cells at $T=1$.
	The results are comparable to those in \cite{Jiang_2021_High_JSC,Xu_2025_High_JCP} without oscillations.
	
	\begin{figure}[htb!]
		\centering
		\begin{subfigure}{0.48\textwidth}
			\centering
			\includegraphics[width=\textwidth]{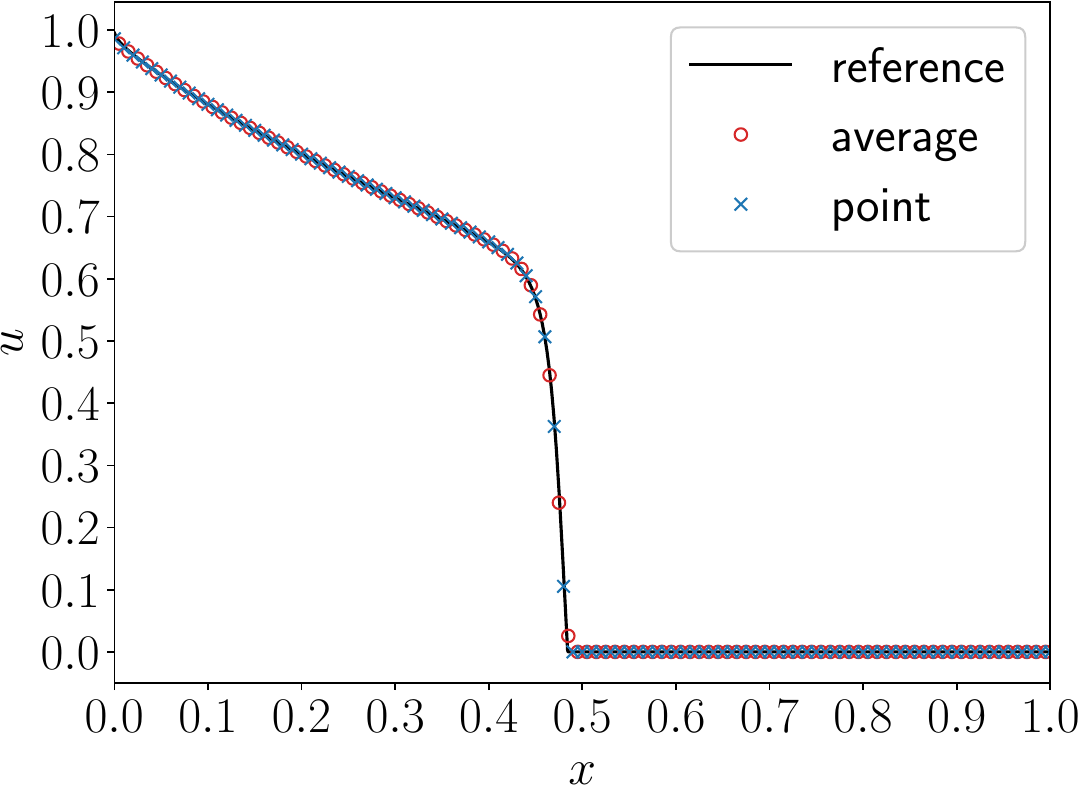}
			\caption{The IBVP without gravity}
		\end{subfigure}
		~
		\begin{subfigure}{0.48\textwidth}
			\centering
			\includegraphics[width=\textwidth]{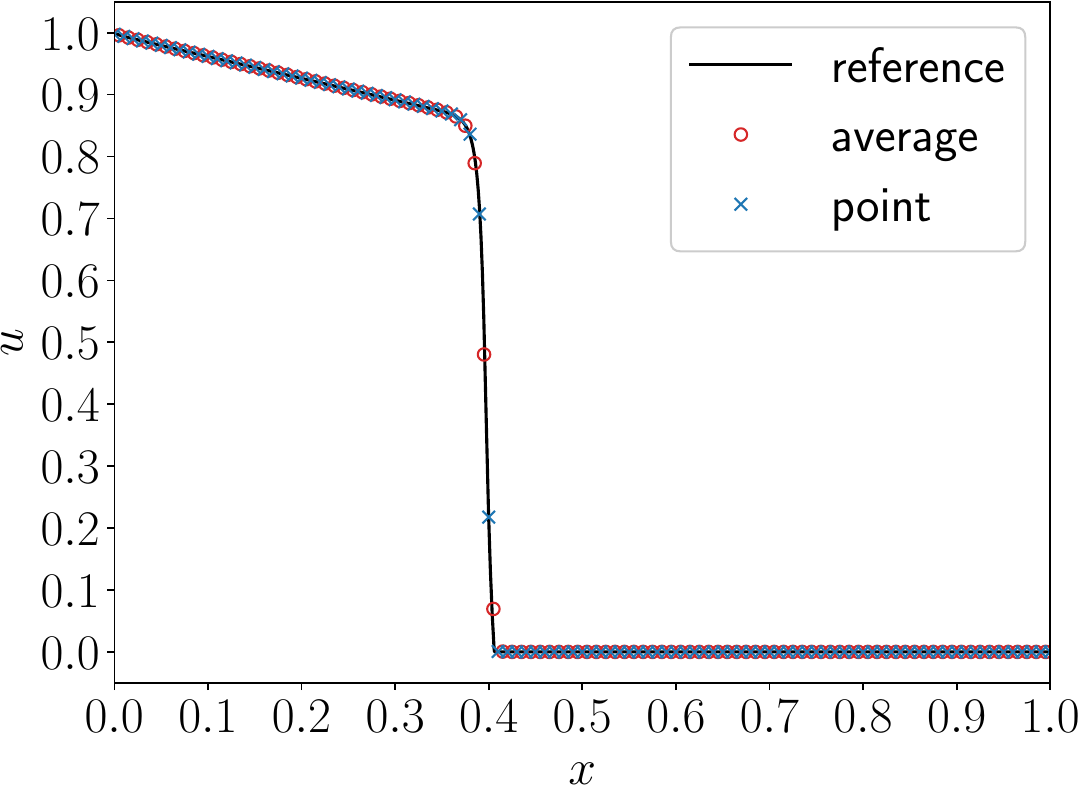}
			\caption{The IBVP with gravity}
		\end{subfigure}
		
		\vspace{10pt}
		\begin{subfigure}{0.48\textwidth}
			\centering
			\includegraphics[width=\textwidth]{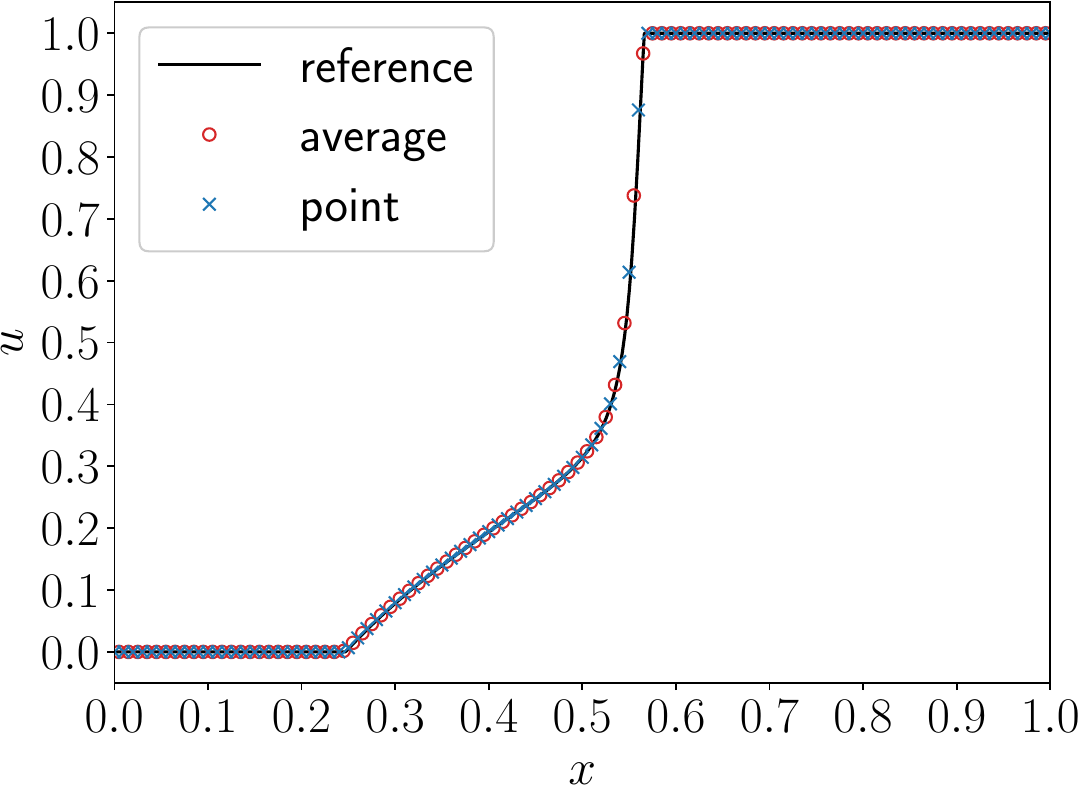}
			\caption{The RP without gravity}
		\end{subfigure}
		~
		\begin{subfigure}{0.48\textwidth}
			\centering
			\includegraphics[width=\textwidth]{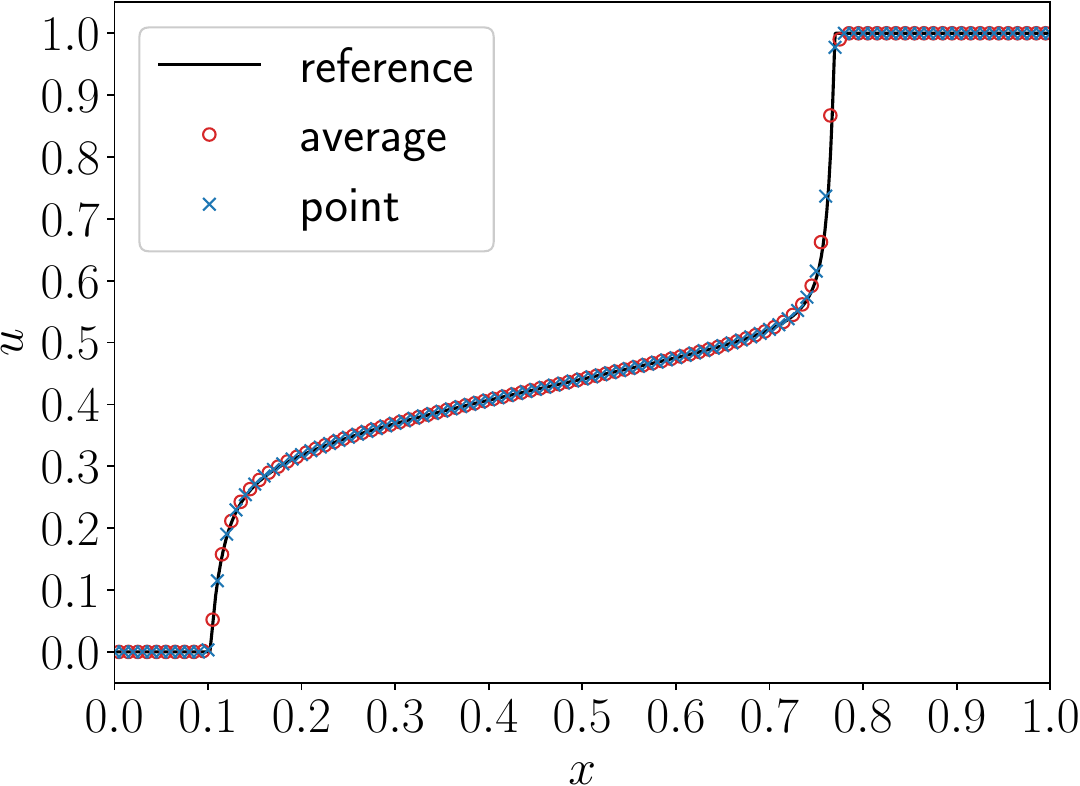}
			\caption{The RP with gravity}
		\end{subfigure}
		\caption{Example~\ref{ex:1d_buckley_leverett}.
			The numerical solutions obtained by the proposed AF method with the MPP limiting and the reference solutions obtained by the first-order scheme with $2000$ cells.}
		\label{fig:1d_buckley_leverett}
	\end{figure}
\end{example}

\begin{example}[2D Buckley Leverett equation]\label{ex:2d_buckley_leverett}
	This example solves $u_t + f(u)_x + g(u)_y = \epsilon(u_{xx} + u_{yy})$, with $\epsilon=0.01$ and
	\begin{equation*}
		f(u) = \frac{u^2}{u^2 + (1-u)^2},
		\quad
		g(u) = \frac{(1 - 5(1-u)^2)u^2}{u^2 + (1-u)^2}.
	\end{equation*}
	The initial value is $u(\bx,0) = 1$ if $x^2+y^2 < 1/2$, and $0$ otherwise.
	The computational domain is $[-1.5, 1.5]\times[-1.5,1.5]$ with outflow boundary conditions and the problem is solved until $T=0.5$.
	
	The numerical solutions with the MPP limiting are shown in Figure~\ref{fig:2d_buckley_leverett}, obtained with the proposed AF method and $100\times100$ cells.
	The numerical solutions with the limiting agree well with the results in the literature \cite{Xu_2025_High_JCP}.
	One also observes that the limiting is important for suppressing oscillations and capturing the correct solution.
	
	\begin{figure}[htb!]
		\begin{subfigure}[t]{0.4\textwidth}
			\centering
			\includegraphics[width=\textwidth]{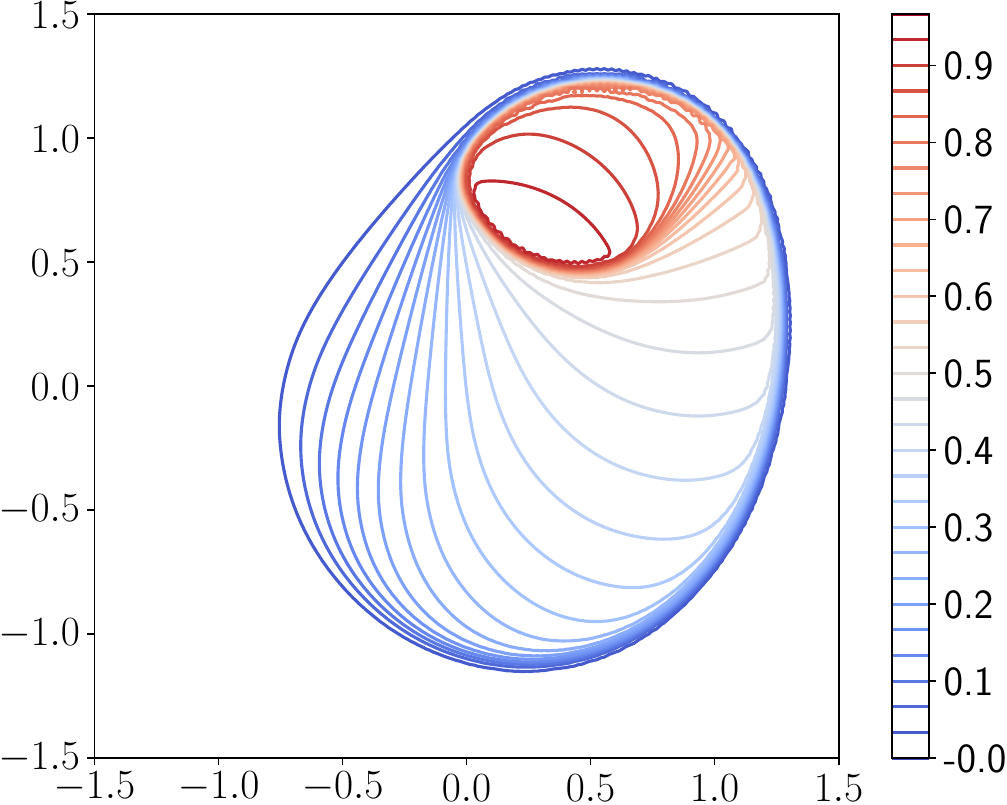}
		\end{subfigure}
		~
		\begin{subfigure}[t]{0.58\textwidth}
			\centering
			\includegraphics[width=\textwidth]{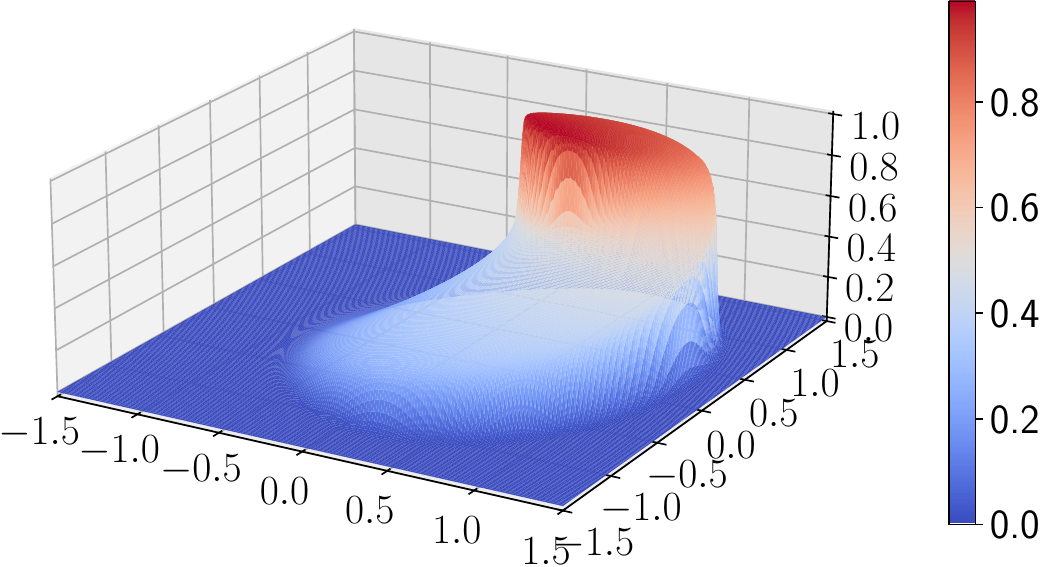}
		\end{subfigure}
		
		\vspace{2pt}
		\begin{subfigure}[t]{0.4\textwidth}
			\centering
			\includegraphics[width=\textwidth]{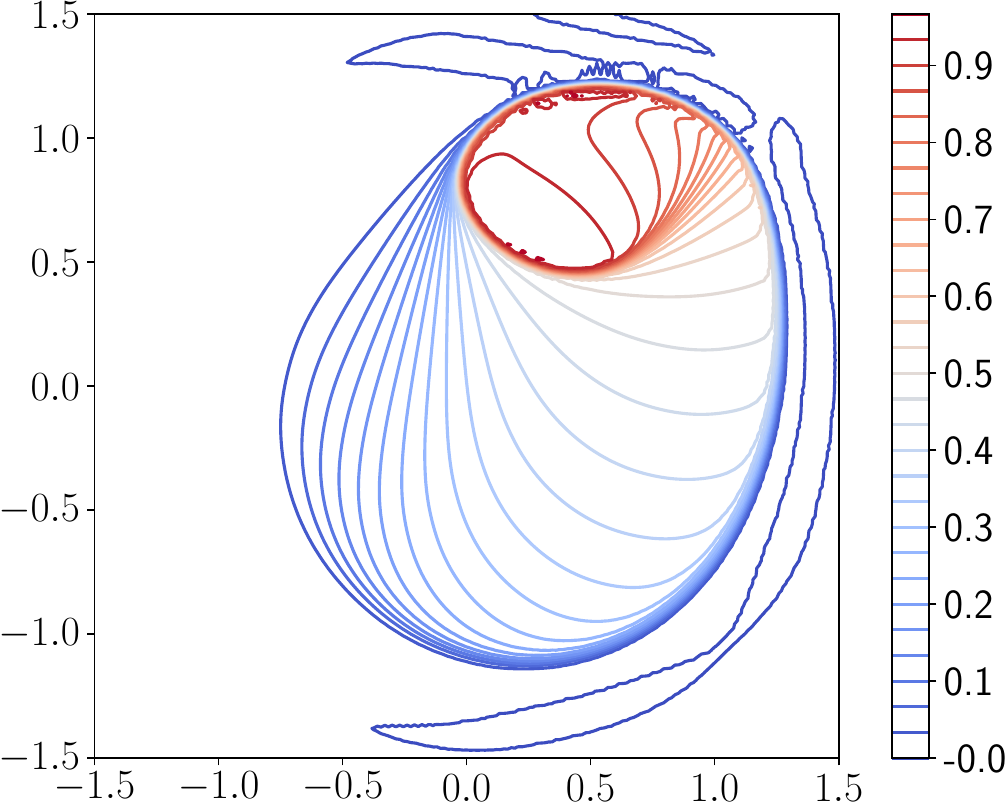}
		\end{subfigure}
		~
		\begin{subfigure}[t]{0.58\textwidth}
			\centering
			\includegraphics[width=\textwidth]{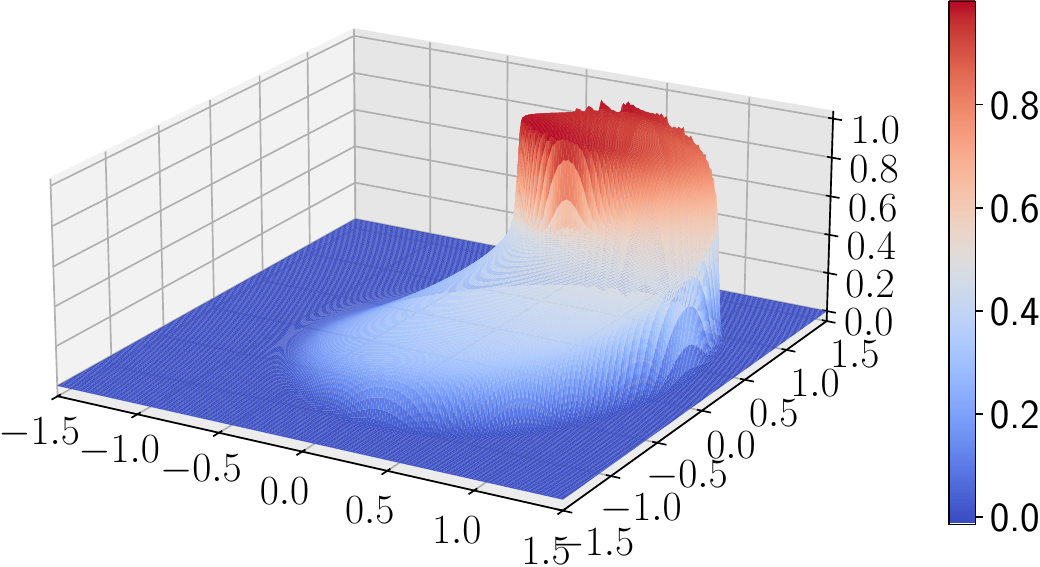}
		\end{subfigure}
		\caption{Example~\ref{ex:2d_buckley_leverett}.
			The numerical solutions obtained by the proposed AF method with (upper) and without (lower) the MPP limiting.}
		\label{fig:2d_buckley_leverett}
	\end{figure}
\end{example}

\begin{example}[1D strongly degenerate case]\label{ex:1d_strongly_degenerate}
	This example solves \eqref{eq:1d_cde_dc} with $f(u) = u^2$ and
	\begin{equation*}
		a(u) =
		\begin{cases}
			\epsilon(u+0.25), &\text{if}\ u \leqslant -0.25, \\
			\epsilon(u-0.25), &\text{if}\ u \geqslant 0.25, \\
			0, &\text{otherwise}, \\
		\end{cases}
	\end{equation*}
	where $\epsilon=0.1$.
	The initial condition is
	\begin{equation*}
		u(x,0) =
		\begin{cases}
			1, &\text{if}\ \abs{x+\frac{1}{\sqrt{2}}} < 0.4, \\
			-1, &\text{if}\ \abs{x-\frac{1}{\sqrt{2}}} < 0.4, \\
			0, &\text{otherwise}.
		\end{cases}
	\end{equation*}
	The computational domain is $[-2, 2]$ with outflow boundary conditions.
	
	Figure~\ref{fig:1d_strongly_degenerate_cases} shows the numerical solutions obtained by the proposed AF method with and without the limiting on a mesh of $100$ cells at $T=0.7$.
	To suppress oscillations, the blending based on the smoothness indicator with $\gamma=1$ is activated.
	The proposed method captures the sharp interface as well as the kinks, which match those in \cite{Jiang_2021_High_JSC,Xu_2025_High_JCP} very well.
	The smoothness indicator-based blending is effective in suppressing oscillations.
	This example also provides numerical evidence for the importance of the double conservation for the cell average, otherwise the numerical solution deviates from the reference near the sharp interface at $x=0$ in the right panel.
	
	\begin{figure}[htb!]
		\centering
		\includegraphics[width=0.32\textwidth]{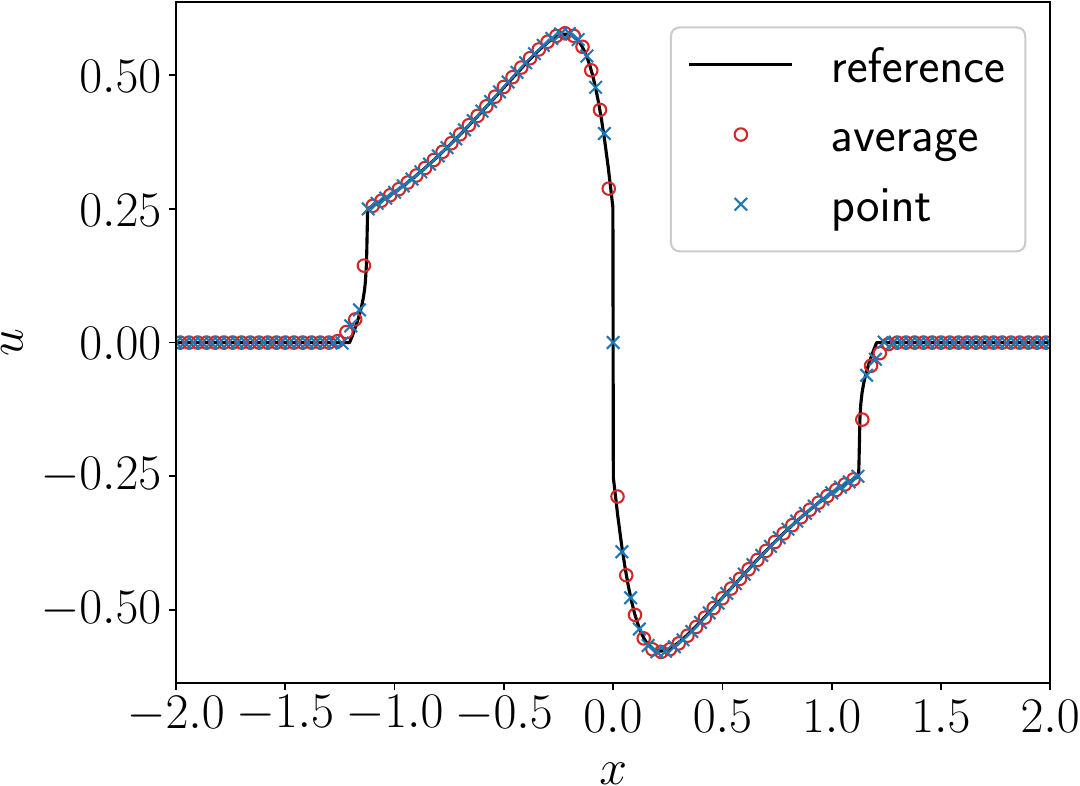}
		\includegraphics[width=0.32\textwidth]{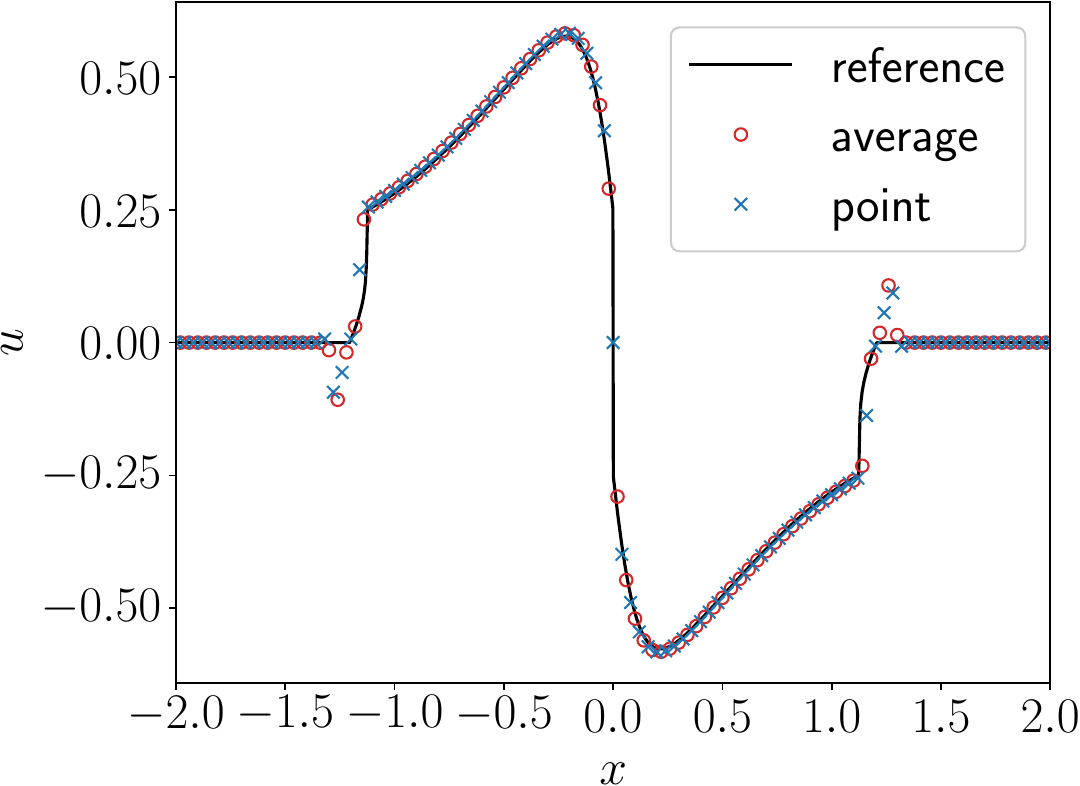}
		\includegraphics[width=0.32\textwidth]{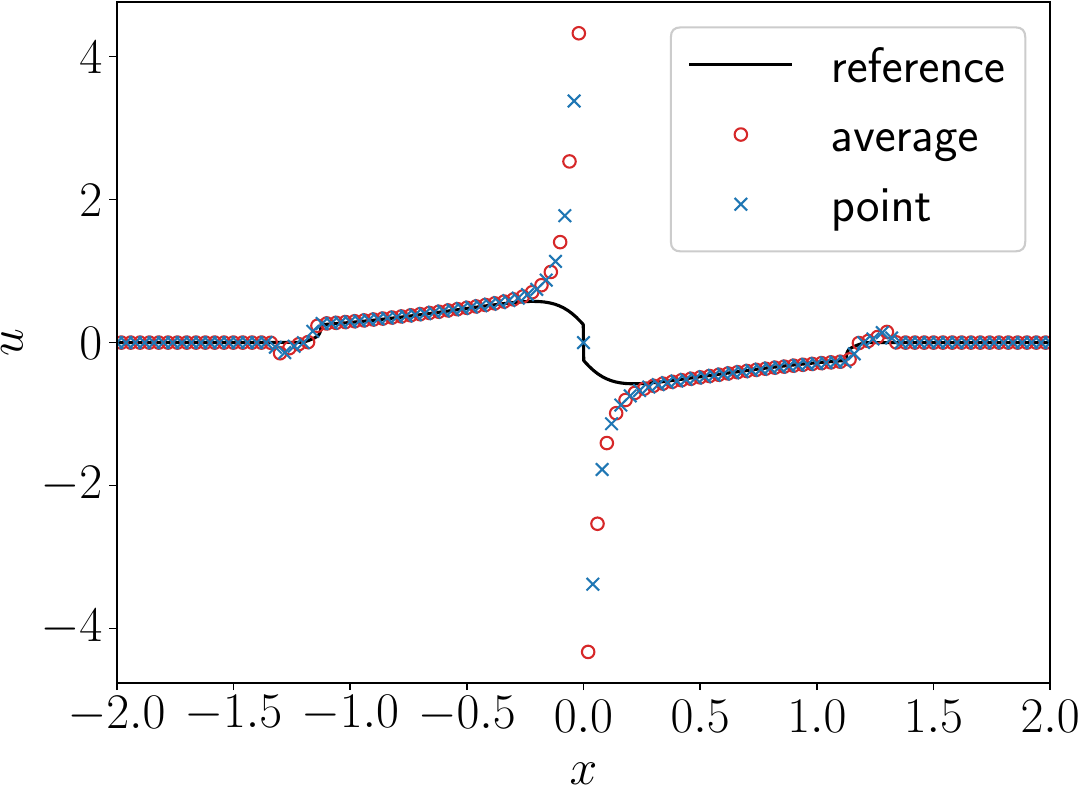}
		\caption{Example~\ref{ex:1d_strongly_degenerate}.
			The numerical solutions obtained by the AF methods and the reference solution obtained by the first-order scheme with $2000$ cells.
			Left: the proposed doubly conservative AF method with the smoothness-indicator-based blending with $\gamma=1$,
			middle: the same as the left but without the smoothness-indicator-based blending,
			right: the AF method \eqref{eq:1d_af_nodc} without double conservation and limiting.}
		\label{fig:1d_strongly_degenerate_cases}
	\end{figure}
\end{example}

\begin{example}[2D strongly degenerate case]\label{ex:2d_strongly_degenerate}
	This example solves $u_t + (u^2)_x + (u^2)_y = \epsilon(a(u)_{xx} + a(u)_{yy})$ with $\epsilon=0.1$ and
	\begin{equation*}
		a(u) =
		\begin{cases}
			u - \frac14, &\text{if}\ u > \frac14, \\
			u + \frac14, &\text{if}\ u < -\frac14, \\
			0, &\text{otherwise}. \\
		\end{cases}
	\end{equation*}
	The initial condition is
	\begin{equation*}
		u(\bx,0) =
		\begin{cases}
			1, &\text{if}\ (x+\frac{1}{2})^2 + (y+\frac{1}{2})^2 < \frac{4}{25}, \\
			-1, &\text{if}\ (x-\frac{1}{2})^2 + (y-\frac{1}{2})^2 < \frac{4}{25}, \\
			0, &\text{otherwise},
		\end{cases}
	\end{equation*}
	and the computational domain is $[-1.5, 1.5]\times[-1.5,1.5]$ with outflow boundary conditions.
	
	The numerical solution with the MPP limiting at $T=0.5$ is shown in Figure~\ref{fig:2d_strongly_degenerate}, obtained with $100\times100$ cells.
	As in the last example, the smoothness-indicator-based blending is activated only in the upper row with $\gamma=1$, which effectively reduces the oscillations.
	The proposed AF method can capture the sharp interface along the diagonal $y=-x$, and the results agree well with those in \cite{Jiang_2021_High_JSC,Xu_2025_High_JCP}.
	
	\begin{figure}[htb!]
		\centering
		\includegraphics[width=0.44\textwidth]{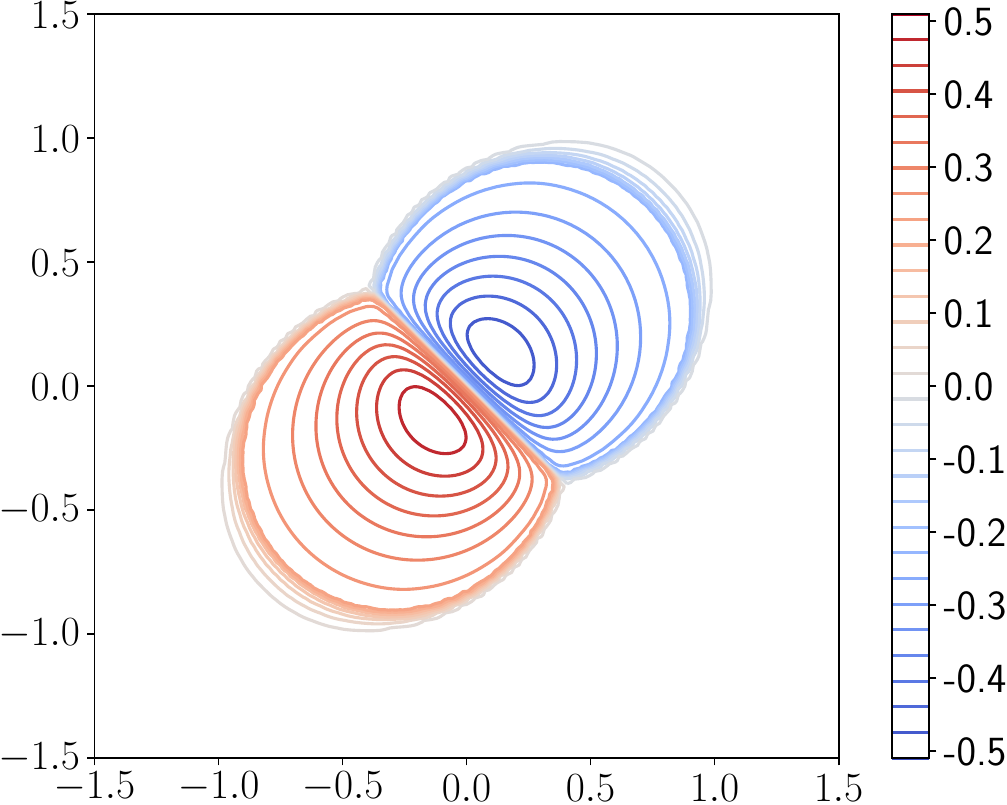}
		\includegraphics[width=0.54\textwidth]{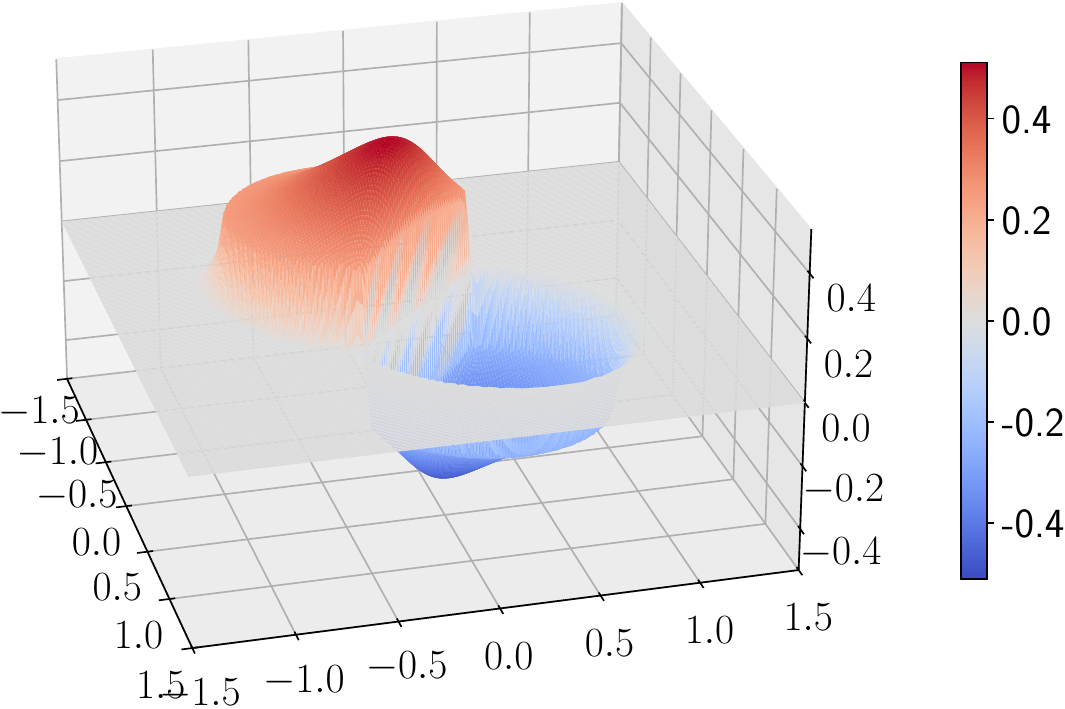}
		
		\vspace{2pt}
		\includegraphics[width=0.44\textwidth]{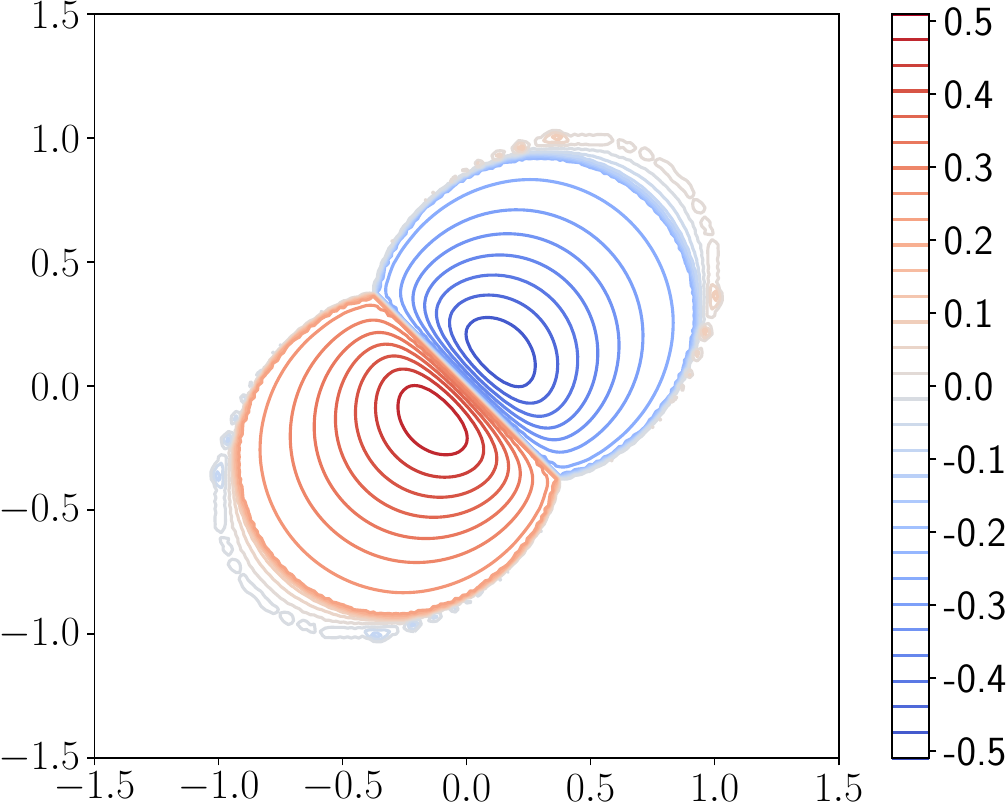}
		\includegraphics[width=0.54\textwidth]{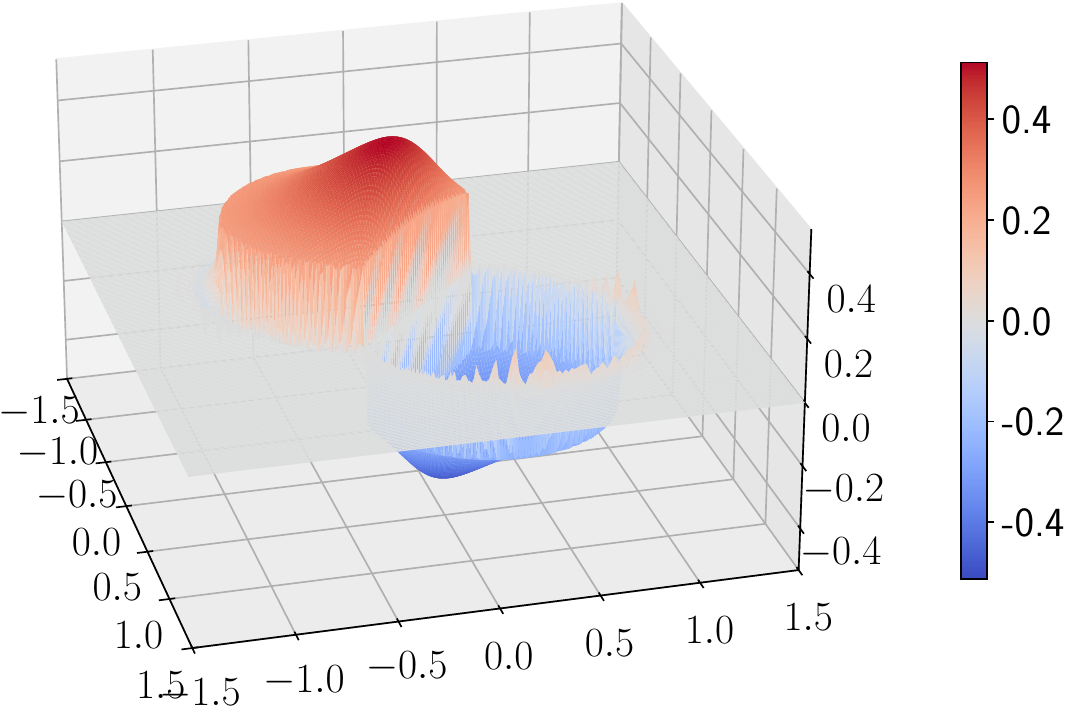}
		\caption{Example~\ref{ex:2d_strongly_degenerate}.
			The numerical solutions obtained by the proposed doubly conservative AF method with (upper) and without (lower) the smoothness-indicator-based limiting with $\gamma=1$.}
		\label{fig:2d_strongly_degenerate}
	\end{figure}
\end{example}

\section{Conclusion}\label{sec:conclusion}

We have developed a compact fourth-order active flux (AF) method for scalar degenerate convection--diffusion equations on Cartesian meshes.
The method uses the same set of evolved unknowns as the standard third-order AF method.
Fourth-order convection discretization is obtained by incorporating a downwind point value into the biased point-value stencil, while diffusion is treated by directly discretizing the diffusion potential with compact fourth-order operators.
This construction avoids auxiliary gradient variables, hyperbolic reformulations, pseudo-time stepping, or inner iterations in the existing AF methods for diffusion problem.
The diffusion contribution to the cell-average evolution is written in a doubly conservative form to deal with degeneracy.
To enforce the discrete maximum principle without destroying the underlying conservation and double conservation form, a monolithic convex limiter is applied simultaneously to convection fluxes and diffusion potentials.
Together with the projection of point values onto the admissible interval, the resulting method preserves the prescribed bounds under the derived CFL condition.
An optional smoothness-indicator-based blending provides additional oscillation control near discontinuities.
For the 1D linear problem, Fourier analysis confirms fourth-order spatial accuracy, and von Neumann analysis shows that the proposed AF method has larger stable CFL limits than the fourth-order $P^3$ LDG method for the explicit RK methods considered.
The reported 1D comparison also shows that the AF method requires fewer DoFs and a lower estimated FLOP count at comparable accuracy.
Numerical experiments demonstrate the observed fourth-order convergence, the MPP property, and the ability of the method to resolve shocks and sharp wave fronts in strongly degenerate problems.
They also provide numerical evidence that the double conservation is important for dealing with the tested strongly degenerate cases.
Extensions to systems such as the Navier--Stokes equations, broader performance comparisons based on measured computational cost, and implicit or IMEX time integration remain subjects for future work.

\section*{Acknowledgements}
This work was supported by 1+1+1 CUHK-CUHK(SZ)-GDSTC Joint Collaboration Fund (No. 2025A0505000079) and University Development Fund (No. UDF01004388).
The author would also like to thank Christian Klingenberg in W\"urzburg, Wasilij Barsukow in Bordeaux, and Kailiang Wu in Shenzhen for helpful discussions during the preparation of this work.


\appendix

\section{A compact discretization for the diffusion term without double conservation}\label{app:nodc}
Consider the equation
\begin{equation}\label{eq:1d_cde_nodc}
	u_t + f(u)_x = (\kappa(u)u_x)_x, 
\end{equation}
which can be related to \eqref{eq:1d_cde_dc} by
\begin{equation*}
	a(u) = \int^u \kappa(s) \dd s,\quad
	a(u)_x = \kappa(u)u_x,\quad
	\kappa(u)\geqslant 0.
\end{equation*}
One can derive the following AF method for the diffusion part
\begin{subequations}\label{eq:1d_diff_nodc}
	\begin{align}
		\frac{\dd \bar{u}_i }{\dd t} &= \frac{1}{\Delta x}\left( \kappa(u_{\xr})(u_x)_{\xr} - \kappa(u_{\xl})(u_x)_{\xl} \right), \label{eq:1d_diff_nodc_av_u} \\
		(u_x)_{\xr} &= D^{\texttt{4th}}_{c}(u)_{\xr}, \label{eq:1d_diff_nodc_pnt_p} \\
		\frac{\dd u_{\xr} }{\dd t} &= \tilde{D}^{2,\texttt{4th}} (\kappa, u)_{\xr}, \label{eq:1d_diff_nodc_pnt_u}
	\end{align}
\end{subequations}
where $D^{\texttt{4th}}_{c}(u)_{\xr}$ and $\tilde{D}^{2,\texttt{4th}} (\kappa, u)_{\xr}$ approximate $(u_x)_{\xr}$ and $(\kappa(u)u_x)_x$, respectively.
A natural approach is to use central finite differences because no upwind mechanism is available for diffusion.
For $(\kappa(u)u_x)_x$, one can apply the operator \eqref{eq:1d_diff_central_4th} twice, i.e., taking
\begin{equation}\label{eq:1d_diff2_twice}
	\tilde{D}^{2,\texttt{4th}} (\kappa, u)_{\xr} = D^{\texttt{4th}}_{c} \left(\kappa(u) D^{\texttt{4th}}_{c}(u)\right)_{\xr}.
\end{equation}
However, such a discretization leads to a larger spatial stencil than the original AF method.
For the constant-coefficient case $a(u) = \kappa_0 u$, equation~\eqref{eq:1d_diff_nodc_pnt_u} becomes
\begin{align*}
	\frac{\dd u_{\xr} }{\dd t}
	=&\ \frac{\kappa_0}{4\Delta x^2}\left( u_{i-\frac32} - 4\bar{u}_{i-1}  + 8u_{i-\frac12} + 4\bar{u}_{i} - 18u_{\xr} \right. \nonumber\\
	&\left. +\ 4\bar{u}_{i+1} + 8u_{i+\frac32} - 4\bar{u}_{i+2} + u_{i+\frac52} \right).
\end{align*}
The spatial stencil spans the cells $I_{i-2}$ to $I_{i+2}$, extending beyond the immediately neighboring cells.


To maintain the compact stencil of the AF method, which depends on only the direct neighboring cells, it is possible to adapt the $4$th-order compact finite difference \cite{Kamakoti_2010_High_SJSC} to our AF method.
The discretization reads
\begin{equation}\label{eq:1d_diff_compact_kappa_u}
	\tilde{D}^{2,\texttt{4th}} (\kappa, u)_{\xr} = \dfrac{4}{\Delta x^2}\sum_{m=-2}^{2}\sum_{n=-2}^{2} b_{mn}\kappa\left(u_{i+\frac{m+1}{2}}\right) u_{i+\frac{n+1}{2}}.
\end{equation}
Here the coefficients $b_{mn}$ are obtained from the order conditions, sparsity requirement, and symmetry.
More precisely, $b_{mn}$ is the $(m,n)$th entry of the matrix $B$ defined as follows,
\begin{equation*}
	B = \begin{bmatrix}
		0 & 0 \\ 0 & M \\
	\end{bmatrix}
	- \begin{bmatrix}
		M & 0 \\ 0 & 0 \\
	\end{bmatrix}, \quad
	M =
	\begin{bmatrix}
		\frac{1}{8} & -\frac{1}{6} & \frac{1}{24} & 0 \\[3pt]
		-\frac{1}{6} & -\frac{3}{8} & \frac{2}{3} & -\frac{1}{8} \\[3pt]
		\frac{1}{8} & -\frac{2}{3} & \frac{3}{8} & \frac{1}{6} \\[3pt]
		0 & -\frac{1}{24} & \frac{1}{6} & -\frac{1}{8} \\
	\end{bmatrix}.
\end{equation*}
Note that the spatial stencil of \eqref{eq:1d_diff_compact_kappa_u} uses the DoFs in the cell $I_{i-1}$ to $I_{i+1}$.
Moreover, for constant diffusion coefficient $\kappa_0$, the scheme \eqref{eq:1d_diff_compact_kappa_u} reduces to
\begin{equation}\label{eq:1d_diff_nodc_constant}
	\frac{\dd u_{\xr} }{\dd t} = \frac{4\kappa_0}{\Delta x^2}\left(- \frac{1}{12}u_{i-\frac12} + \frac{4}{3}u_{i} - \frac{5}{2}u_{i+\frac12} + \frac{4}{3}u_{i+1} - \frac{1}{12}u_{i+\frac32}\right).
\end{equation}
Finally, the compact 1D AF method for solving \eqref{eq:1d_cde_nodc} is
\begin{align}
	&\frac{\dd \bar{u}_i }{\dd t} + \frac{1}{\Delta x}\left( f_{\xr} - f_{\xl} \right) = \frac{1}{\Delta x}\left( \kappa(u_{\xr})(u_x)_{\xr} - \kappa(u_{\xl})(u_x)_{\xl} \right), \nonumber\\
	&\frac{\dd u_{\xr} }{\dd t} + D^{\texttt{4th}}_{+} (f^{+})_{\xr} + D^{\texttt{4th}}_{-}(f^{-})_{\xr} = \tilde{D}^{2,\texttt{4th}} (\kappa, u)_{\xr}, \nonumber\\
	&(u_x)_{\xr} = D^{\texttt{4th}}_{c}(u)_{\xr}, \label{eq:1d_af_nodc}
\end{align}
with $\tilde{D}^{2,\texttt{4th}}$ defined in \eqref{eq:1d_diff_compact_kappa_u}.
The above scheme maintains the compact stencil and aligns with \eqref{eq:1d_af_dc} for the constant case $a(u) = \kappa_0 u$.

\end{document}